\documentclass[12pt]{article}
\usepackage{amsmath, amsthm, amssymb, stmaryrd, fullpage}
\usepackage[all]{xy}

\numberwithin{equation}{subsection}
\newtheorem{theorem}[equation]{Theorem}
\newtheorem{lemma}[equation]{Lemma}
\newtheorem{conj}[equation]{Conjecture}
\newtheorem{cor}[equation]{Corollary}
\newtheorem{prop}[equation]{Proposition}
\theoremstyle{definition}
\newtheorem{question}[equation]{Question}
\newtheorem{defn}[equation]{Definition}
\newtheorem{remark}[equation]{Remark}
\newtheorem{hypo}[equation]{Hypothesis}
\newtheorem{example}[equation]{Example}
\newtheorem{convention}[equation]{Convention}

\newcommand{\AAA}{\mathbb{A}}
\newcommand{\CC}{\mathbb{C}}
\newcommand{\HH}{\mathbb{H}}
\newcommand{\NN}{\mathbb{N}}
\newcommand{\PP}{\mathbb{P}}
\newcommand{\QQ}{\mathbb{Q}}
\newcommand{\RR}{\mathbb{R}}
\newcommand{\ZZ}{\mathbb{Z}}
\newcommand{\calC}{\mathcal{C}}
\newcommand{\calD}{\mathcal{D}}
\newcommand{\calE}{\mathcal{E}}
\newcommand{\calF}{\mathcal{F}}
\newcommand{\calG}{\mathcal{G}}
\newcommand{\calH}{\mathcal{H}}
\newcommand{\calI}{\mathcal{I}}
\newcommand{\calL}{\mathcal{L}}
\newcommand{\calM}{\mathcal{M}}
\newcommand{\calN}{\mathcal{N}}
\newcommand{\calO}{\mathcal{O}}
\newcommand{\calP}{\mathcal{P}}
\newcommand{\calU}{\mathcal{U}}
\newcommand{\gotha}{\mathfrak{a}}
\newcommand{\gothm}{\mathfrak{m}}
\newcommand{\gotho}{\mathfrak{o}}

\newcommand{\an}{\mathrm{an}}

\newcommand{\dual}{\vee}
\newcommand{\del}{\partial}

\newcommand{\be}{\mathbf{e}}
\newcommand{\bbf}{\mathbf{f}}
\newcommand{\bv}{\mathbf{v}}
\newcommand{\bw}{\mathbf{w}}
\newcommand{\bx}{\mathbf{x}}

\DeclareMathOperator{\ad}{ad}
\DeclareMathOperator{\Ext}{Ext}
\DeclareMathOperator{\Fil}{Fil}
\DeclareMathOperator{\Gal}{Gal}
\DeclareMathOperator{\GL}{GL}
\DeclareMathOperator{\gp}{gp}
\DeclareMathOperator{\Hom}{Hom}
\DeclareMathOperator{\id}{id}
\DeclareMathOperator{\image}{im}
\DeclareMathOperator{\Isoc}{Isoc}
\DeclareMathOperator{\LNM}{LNM}
\DeclareMathOperator{\Maxspec}{Max}
\DeclareMathOperator{\Spec}{Spec}
\DeclareMathOperator{\speci}{sp}
\DeclareMathOperator{\Spf}{Spf}
\DeclareMathOperator{\tr}{tr}
\DeclareMathOperator{\triv}{triv}
\DeclareMathOperator{\ULNM}{ULNM}

\begin{document}

\title{Semistable reduction for overconvergent $F$-isocrystals, I:
Unipotence and logarithmic extensions}
\author{Kiran S. Kedlaya \\ Department of Mathematics, Room 2-165 \\ Massachusetts
Institute of Technology \\ 77 Massachusetts Avenue \\
Cambridge, MA 02139 \\
\texttt{kedlaya@mit.edu}}
\date{January 20, 2007}

\maketitle

\centerline{\textit{Dedicated to Pierre Berthelot}}

\begin{abstract}
Let $X$ be a smooth variety over a field $k$
of characteristic $p>0$, and let $\calE$ be an overconvergent
isocrystal on $X$. 
We establish a criterion for the existence of a ``canonical logarithmic
extension'' of $\calE$ to a smooth compactification $\overline{X}$
of $X$ whose complement is a strict normal crossings divisor. We also
obtain some related results, including a form of Zariski-Nagata purity
for isocrystals.
\end{abstract}

\tableofcontents

\section{Introduction}

This paper is intended as the first in a series in which we pursue a
``semistable reduction'' theorem for overconvergent $F$-isocrystals,
a class of $p$-adic analytic objects associated to schemes of finite type
over a field of characteristic $p>0$.
Such a theorem would have consequences for the theory of rigid
cohomology, in which overconvergent $F$-isocrystals play the role of
coefficient objects of locally constant rank.
In this introduction, we give a high-level description of a complex
analytic model situation and the $p$-adic situation that imitates
it, a bit about intended applications, and the structure of the paper.
For a more detailed description of the questions we will be
considering in subsequent papers, see Section~\ref{sec:finale}.

\subsection{An analogy: complex local systems}

Let $X \hookrightarrow \overline{X}$ be an open immersion of smooth 
varieties over $\CC$, with $\overline{X}$ proper and
$Z = \overline{X} \setminus X$ a strict normal crossings divisor.
(Here and throughout, ``variety'' will be used as shorthand for
``reduced, separated scheme of finite type'' over some field.)
A \emph{$\nabla$-module} on the complex analytic space $X^{\an}$ 
consists of a coherent locally free sheaf $\calE$ of
$\calO_{X^{\an}}$-modules (or equivalently, a holomorphic vector bundle)
equipped with an integrable connection.
The integrability condition means that $\calE$ admits a basis of
horizontal sections on any contractible open subset; these fit together
to form a local system of finite dimensional
$\CC$-vector spaces on $X^{\an}$. (In fact,
the categories of $\nabla$-modules and of local systems of
finite dimensional $\CC$-vector spaces are equivalent, by the easy part
of the Riemann-Hilbert correspondence.)

Suppose that $X$ is connected, so that $\calE$ has some rank $n$
everywhere. Associated to $\calE$ (or rather, from its associated
local system) is a monodromy representation $\rho: \pi_1(X^{\an})
\to \GL_n(\CC)$ of the (topological) fundamental group of $X^{\an}$.
Specifically, given a pointed loop, one analytically continues a basis of
local horizontal sections along the loop, and compares the basis before and
after this parallel transport.

Given a component $D$ of $Z$, one obtains from
$\rho$ a new representation by restriction to the subgroup of $\pi_1(X^{\an})$
generated by some loop winding once around $D$ (with the
correct orientation). Of course this subgroup depends
on the choice of the loop, but that choice acts on the loop by a conjugation
in $\pi_1(X^{\an})$, 
and so does not alter the isomorphism class of the restricted
representation. That restriction is called the \emph{local monodromy
representation} associated to $D$.

The local monodromy representation measures the ``badness'' of the
singularities of the connection along $D$. For instance, if the
connection extends without singularities across $D$, the local
monodromy representation is a trivial representation. More interestingly,
by a theorem of Deligne \cite[Proposition~II.5.2]{deligne},
the local monodromy representation is unipotent (i.e., its semisimplification
is a direct sum of trivial representations) if and only if the $\nabla$-module
extends to a log-$\nabla$-module with logarithmic singularities
and nilpotent residues along $D$;
such an extension is unique if it exists.
(This uniqueness relies crucially on the nilpotent residue condition;
otherwise many distinct extensions are possible.)
In particular, the existence of such a
``canonical logarithmic extension'' (the ``prolongement canonique''
of \cite{deligne})
is determined by a codimension 1
criterion, so its existence on $\overline{X}$ minus a codimension 2
subscheme implies its existence on $\overline{X}$
\cite[Corollaire~II.5.8]{deligne}.

For local systems of ``algebro-geometric origin'', e.g., the $i$-th
relative Betti cohomology of a smooth proper morphism to $X$,
one typically obtains a canonical logarithmic extension after pulling back
along a suitable finite cover of $\overline{X}$. This can be shown
``extrinsically'', using semistable reduction of varieties, but a more
intrinsic approach involves recognizing such local systems as analytic
objects equipped with extra data, namely variations of Hodge
structures. (At this point our discussion, being purely of motivational
nature, will turn unabashedly cursory; see \cite{griffiths} for a
more comprehensive overview.)

A \emph{polarized variation of Hodge structures} on $X$ consists of a local
system of finitely generated $\ZZ$-modules on $X^{\an}$, plus some
additional Hodge-theoretic data which we will not describe here, 
save to mention
the principal example (arising from a theorem of Griffiths): 
the $i$-th cohomology of a family of smooth projective
complex analytic varieties. A basic fact about polarized variations
of Hodge structures is the \emph{monodromy theorem}, due in this form
to Borel \cite[Lemma~4.5]{schmid}: the local monodromy representation
associated to any component of $Z$ is \emph{quasi-unipotent}, i.e.,
becomes unipotent
upon further restriction to a subgroup of finite index.

{}From the monodromy theorem,
one easily deduces the following. 
Given a $\nabla$-module $\calE$ on $X^{\an}$ whose associated 
local system can be obtained from a polarized variation of Hodge structures
(by tensoring over $\ZZ$ with $\CC$),
for any closed point $x$ of $\overline{X}$, one can find
an open neighborhood $U$ of $x$ in $\overline{X}$ and a finite
cover $f: V \to U$ such that $V$ is etale over $U \cap X$, $V$ is smooth,
$f^{-1}(U \cap Z)$ is a strict normal 
crossings divisor on $V$, and $f^* \calE$ extends to a
log-$\nabla$-module on $V$ with logarithmic singularities and
nilpotent residues along $f^{-1}(U \cap Z)$.

It is a bit less clear how to patch things together globally
without further analysis of the local situations,
but using resolution of singularities, one can at least
assert that there is a proper, dominant, generically finite morphism
$f: \overline{Y} \to \overline{X}$ with $\overline{Y}$ smooth
and $f^{-1}(Z)$ a strict normal crossings divisor, such that
$f^* \calE$ extends to a log-$\nabla$-module everywhere on
$\overline{Y}$, with logarithmic singularities and nilpotent
residues along $f^{-1}(Z)$.
We summarize
this situation by saying that $\calE$ 
``admits semistable reduction''. (The reason for this
terminological choice is that
when $\calE$ comes from the cohomology of a family of
varieties, one is guaranteed to have the desired property if the
family pulls back to a semistable family over $\overline{Y}$.)

\subsection{Extension of overconvergent isocrystals}

We now consider a $p$-adic analogue of the situation of the
previous subsection. This will be appropriately vague for an
introduction; see Section~\ref{sec:finale} for a summary
in more precise language.

Let $X \hookrightarrow \overline{X}$ 
be an open immersion of smooth $k$-varieties, for
$k$ a field of characteristic $p>0$, such that
$Z = \overline{X} \setminus X$ is a strict normal crossings divisor.
Let $\calE$ be an isocrystal on $X$ which is overconvergent along $Z$;
this is a positive characteristic analogue of a $\nabla$-module with
some additional convergence conditions, constructed using
$p$-adic rigid analytic geometry.
Although it is not so easy to define a $p$-adic local monodromy group,
one can at least give meaning to the assertion that
``$\calE$ has constant/unipotent local monodromy along $Z$''.
We show (Theorem~\ref{T:extension}) that again 
$\calE$ has unipotent local monodromy if and only if 
$\calE$ admits a ``canonical logarithmic extension'' to $\overline{X}$;
that extension will be a convergent log-isocrystal in the sense of
Shiho \cite{shiho1, shiho2}. This in particular implies a form of Zariski-Nagata purity for isocrystals on smooth varieties

Continuing the analogy, one can then ask whether one can associate to
$\calE$ of ``algebro-geometric origin'' a certain global analytic object
that will ensure that $\calE$ admits a canonical logarithmic extension.
The object that provides this control is a Frobenius structure:
the analogue of the monodromy theorem is that the semisimplified
local monodromy representations, being
equipped with Frobenius structures,
necessarily have finite image when restricted to an inertia subgroup. 
This is ``Crew's conjecture'',
now the
$p$-adic local monodromy theorem of Andr\'e
\cite{andre}, Mebkhout \cite{mebkhout}, and the present author
\cite{me-local}.

Thus one expects that one can pull back $\calE$ along a generically
finite cover and get a canonical
logarithmic extension. Note that this is not
at all a trivial consequence of Theorem~\ref{T:extension}, despite
that the fact of an isocrystal having unipotent monodromy
can be checked in codimension 1! The problem arises because
of wild ramification in positive characteristic: 
the analogue of the local construction
in the complex case produces a singular $\overline{Y}$, to which
Theorem~\ref{T:extension} does not (and should not) apply.
Resolving the resulting singularities
(using an alteration in the manner of de Jong \cite{dejong1}) produces
new components whose local monodromy is not \emph{a priori} under control.
We describe the situation in more detail in Section~\ref{sec:finale}.

It should also be noted that the failure to obtain a canonical logarithmic extension on a finite (not just generically finite) cover is also not merely an artifact of the proof technique. One can exhibit examples of overconvergent isocrystals with Frobenius structure
that cannot admit a canonical logarithmic extension after pullback along
\emph{any} finite cover; obstructions to this can be exhibited using the
Newton polygons of the Frobenius action at various points. We plan to 
include an example of this in a subsequent paper.

\subsection{Applications in rigid cohomology}

In the theory of algebraic de Rham cohomology of varieties over a field
of characteristic zero, the ability to ``compactify coefficients'' makes
it possible to prove various finiteness theorems by passing to smooth
proper varieties. With a semistable reduction theorem for overconvergent
$F$-isocrystals, one would hope to obtain
analogous results in rigid cohomology; we now describe some possible
such results.

Shiho \cite{shiho2} has shown that semistable reduction implies
the finite dimensionality of rigid cohomology with coefficients in
an overconvergent $F$-isocrystal. Although one can also prove this
more directly \cite{me-finite}, Shiho's construction may yield
insight into the relative setting, where a direct argument seems more
difficult.

Nakkajima \cite{nakkajima} has shown that semistable reduction implies the
existence of complexes, constructed from log-crystalline cohomology, that
compute the rigid cohomology of an arbitrary scheme of finite type (not
even separated!) over $k$. These complexes may shed some light on the
rigid weight-monodromy conjecture of Mokrane \cite{mokrane}.

Berthelot (private communication) has suggested that semistable reduction
may be of value in the theory of arithmetic $\mathcal{D}$-modules. In
particular, one currently does not know that the restriction
of a holonomic $\mathcal{D}$-module to a closed subscheme is again
holonomic; possibly this can be proved by ``approximating'' the
$\mathcal{D}$-module with overconvergent log-isocrystals.
Ongoing work of Caro may provide a workaround for this problem, but
we still expect semistable reduction to intervene ultimately.

Some of our side results may have their own relevance. For instance, the
fact that a convergent isocrystal admits an overconvergent structure if
the same is true after restriction to an open dense subset
(Proposition~\ref{P:contagion})
 can be used to prove some results in the direction of Berthelot's conjecture
\cite{ber1}
on overconvergence of direct images of smooth proper morphisms. The
point is that the direct images
one is trying to construct exist in the convergent
category by arguments of Ogus \cite{ogus1}, and
can be shown to
exist ``generically'' (on an open dense subset of the original base)
in the overconvergent category using the techniques of
\cite{me-finite}. We intend to amplify these comments elsewhere.

\subsection{Structure of the paper}

We conclude this introduction with a summary of the structure of the
paper. Note that (unlike the rest of this introduction)
these comments only summarize the structure of the present
paper; the structures of subsequent papers in this series
will be described therein.

In Section~\ref{sec:rigid}, we review some notions from rigid
analytic geometry. In particular, we introduce modules with connection
and log-connection,
as well as Berthelot's notions of tubes and strict neighborhoods,
and define overconvergent isocrystals following Berthelot's
treatment in \cite{ber2}.

In Section~\ref{sec:local}, we analyze modules with connection over
the product of a polyannulus with another space. This amounts to recalling
some results from the local
theory of $p$-adic differential equations. In particular, we 
define the notion of a unipotent $\nabla$-module in this context and analyze
its relationship with log-connections.

In Section~\ref{sec:mono}, we specify what we mean for an isocrystal
on a smooth variety
to have ``constant monodromy'' or ``unipotent monodromy'' along the
boundary in some partial compactification.

In Section~\ref{sec:ext}, we state several results to the effect that
the obstruction to extending an isocrystal over a boundary subvariety
is precisely its failure to have constant monodromy along the subvariety.
Although this sort of result is not really needed for semistable reduction, 
such assertions may be of independent interest.

In Section~\ref{sec:logext}, we state a result to the effect that
the obstruction to the existence of a canonical logarithmic extension
of an isocrystal is precisely its failure to have unipotent monodromy.
Our canonical logarithmic extensions will be convergent log-isocrystals
in the sense of Shiho \cite{shiho2}, and some effort is expended
to relate our construction to his.

In Section~\ref{sec:finale},
we conclude by articulating the questions we intend to address in subsequent papers in this series, fleshing out the discussion initiated in this
introduction.

\subsection*{Acknowledgments}

Thanks to Bernard le Stum, Atsushi Shiho, and Liang Xiao for pointing out
errors in prior versions of this paper.
Thanks to Pierre Berthelot, Johan de Jong, Arthur Ogus,
and Nobuo Tsuzuki for helpful
discussions, and to the Universit\'e de Rennes 1 and the
Institute for Advanced Study for their hospitality. Thanks also
to the organizers of the Semaine Cohomologique in Rennes and
the Miami Winter School for the invitations to present this material there.
Special thanks to the referee for efforts above and beyond the call of
duty, in helping to eradicate a plague of errors
in Section~\ref{sec:local}; in particular, most of the final proof of 
Proposition~\ref{P:generization} was kindly provided by the referee.
The author was partially supported by a National Science Foundation
postdoctoral fellowship, and by NSF grants DMS-0111298
and DMS-0400727.

\section{Rigid analytic setup}
\label{sec:rigid}

In this section, we recall briefly the construction of 
overconvergent isocrystals 
on schemes over a field of positive characteristic.
Our reference for notation and terminology in rigid analytic geometry
is \cite{bgr}; see also \cite{fresnel}.
Also see \cite[Chapter~1]{ber2} for more details on the construction of
isocrystals.

\subsection{Initial notations}
\label{subsec:initial}

We first set some notation and terminology conventions, which will
hold throughout the paper unless otherwise specified.
\begin{convention}
Throughout this paper, let $k$ be an arbitrary
field of characteristic $p>0$. When we speak of a ``$k$-variety'',
we will mean a reduced separated (but not necessarily irreducible)
scheme of finite type over $k$;
any additional modifiers
are to be passed through to the underlying scheme (e.g., connected,
irreducible) or to the structural morphism (e.g., smooth,
affine, proper) as appropriate.
\end{convention}

\begin{convention}
Until further notice (specifically, until Section~\ref{sec:logext}),
let $K$ be a field of characteristic 0
complete with respect to
a nonarchimedean absolute value $|\cdot|: K^* \to \RR^+$,
with residue field $k$. 
Let $\Gamma^*$ denote the divisible closure
of the image of $|\cdot|$.
Let $\gotho = \gotho_K$ denote the ring of integers of $K$.
Any norm or seminorm on a $K$-algebra will be assumed to be compatible
with the given norm on $K$; in particular, any finite extension of $K$
carries a unique such norm, which we also denote by $|\cdot|$.
\end{convention}

\begin{remark}
The fact that $K$ will start the paper 
being any field complete for a nonarchimedean absolute value,
and end the paper being discretely valued, reflects a certain
ambivalence in the $p$-adic cohomological community.
It seems that if one's perspective is informed by crystals or formal-scheme
constructions (like Monsky-Washnitzer cohomology), discretely valued
fields are the ones that arise most naturally, whereas if one's perspective
is informed by $p$-adic analysis, then fields like $\CC_p$
and its spherical completion
also arise naturally. We have decided to split the difference, by
carrying along a general $K$ as far as possible, namely
until we begin to invoke Shiho's papers \cite{shiho1, shiho2}.
\end{remark}

\begin{convention}
When forming an $i$-fold product or fibred product in any category, we use
$\pi_1, \dots, \pi_i$ to denote the projections onto the
respective factors.
\end{convention}

\begin{convention}
When any sort of norm is applied to a matrix, we mean this to be the
maximum of the values of the norm on the individual elements of the matrix,
and not any sort of spectral/operator norm.
\end{convention}

\subsection{Tubes, frames, and strict neighborhoods}

We now set up some of the rigid geometry needed to construct isocrystals,
in order to fix notations.

We start with Raynaud's notion of the ``generic fibre'' of an affine
formal scheme \cite[0.2.2]{ber2}. This construction
provides the ambient rigid spaces inside which we will work.
\begin{defn}
Let $P = \Spf A$ be
an affine formal scheme of finite type over $\gotho_K$,
and put $A_K = A \otimes_{\gotho_K} K$ and $P_K = \Maxspec A_K$. Then
$P_K$ is an affinoid space, called the \emph{generic fibre} of $P$.
The points of $P_K$ correspond to quotients
of $A$ which are integral and  finite flat over $\gotho_K$; under
this interpretation, we get a map $\speci: P_K \to P_k$ by tensoring these
quotients with $k$. This is called the \emph{specialization map}.
For any subvariety $U$ of $P_k$, define the \emph{tube of $U$ (within
$P_K$)},
denoted $]U[_P$, as the inverse image $\speci^{-1}(U)$ within $P_K$;
we drop the subscript $P$ in case it is to be understood.
\end{defn}
\begin{remark}
One could relax the restriction that $P$ be affine; see
Remark~\ref{R:nonaffines} for more discussion.
\end{remark}

\begin{defn} \label{D:partial tubes}
Suppose $X$ is a closed subscheme of $P_k$ cut out by the reductions
of $g_1, \dots, g_n \in \Gamma(P, \calO_P)$. Then
\[
]X[_P = \{x \in P_K: |g_i(x)| < 1 \qquad (i=1, \dots, n)\}.
\]
As in \cite[1.1.8]{ber2}, for $\lambda \in (0,1) \cap \Gamma^*$, put
\[
[X]_{P\lambda} = \{x \in P_K: |g_i(x)| \leq \lambda \qquad (i=1, \dots, n)\}
\]
and
\[
]X[_{P\lambda} = \{x \in P_K: |g_i(x)| < \lambda \qquad (i=1, \dots, n)\};
\]
then each $[X]_{P\lambda}$ is rational, and each of the collections
$\{[X]_{P\lambda}\}$ and $\{]X[_{P\lambda}\}$, for 
$\lambda$ running over a sequence in $(0,1) \cap \Gamma^*$ converging to 1,
forms an admissible covering of $]X[_P$ \cite[Proposition~1.1.9]{ber2}.
Again, we drop the subscript $P$ if it is to be understood.
\end{defn}

We now specify a geometric setup we will be using repeatedly; the terminology
is not standard, but will be rather convenient for us.
\begin{defn}
A \emph{frame} (or \emph{affine frame}) is a tuple $(X,Y,P,i,j)$, in which:
\begin{itemize}
\item
$P$ is an affine formal scheme of finite type over $\gotho_K$;
\item
$Y$ is a $k$-variety and $i: Y \hookrightarrow P_k$ is a closed immersion;
\item
$X$ is a $k$-variety and $j: X \hookrightarrow Y$ is an open immersion;
\item
$P$ is smooth over $\gotho_K$ in a neighborhood of $X$.
\end{itemize}
We say that the frame \emph{encloses} the variety $Y$ and/or
the pair $(X,Y)$.
Given two frames
$F = (X,Y,P,i,j)$ and $F' = (X',Y',P',i',j')$, a \emph{morphism}
$F' \to F$ is a diagram of the form
\begin{equation} \label{eq:functor}
\xymatrix{
X' \ar^{j'}@{^{(}->}[r] \ar_w[d] & Y' \ar_v[d]
\ar^{i'}@{^{(}->}[r] & P' \ar_{u}[d] \\
X \ar^{j}@{^{(}->}[r]  & Y
\ar^{i}@{^{(}->}[r]  & P
}
\end{equation}
in which $u$ is smooth in a neighborhood of $X$.
Define the \emph{product frame} $F \times F'$ as the frame
$(X \times_k X', Y \times_k Y', P \times_{\gotho_K} P', i \times i',
j \times j')$; it is equipped with the obvious projection morphisms
$\pi_1: F \times F' \to F$ and
$\pi_2: F \times F' \to F'$.
\end{defn}

\begin{remark} \label{R:nonaffines}
Berthelot considers also the analogous situation in which $P$
is not necessarily affine. However, since our work here is
entirely ``pre-cohomological'', allowing non-affine $P$
would not really add any
generality, since one can always cover such a $P$ with affines,
work locally, and keep track of glueing maps. (This is basically
what Definition~\ref{defn:isocrys} does.) In fact,
one is forced to do this anyway in order to deal with
varieties which do not lift to characteristic zero.
Thus for simplicity,
we have decided to use only affine frames throughout.
(By contrast, when one passes to cohomological considerations, it is
necessary to consider the case where $P$ is proper in order to invoke
Kiehl's finiteness theorem.)
\end{remark}

We next introduce strict neighborhoods, following \cite[1.2]{ber2}.
\begin{defn}
Let $(X,Y,P,i,j)$ be a frame.
An admissible open subset $V$ of $]Y[_P$ containing $]X[_P$ is a
\emph{strict neighborhood} of $]X[_P$ within $]Y[_P$ if the covering 
$\{V, ]Y \setminus X[_P\}$ of $]Y[_P$ is admissible. (Note that the covering 
$\{]X[_P, ]Y \setminus X[_P\}$ of $]Y[_P$ is typically not admissible.)
\end{defn}

To test locally whether an open set is a strict neighborhood, one may
use the following lemma, which is the variant of 
\cite[Proposition~1.2.2]{ber2} described in
\cite[Remarques~1.2.3(iii)]{ber2}.
\begin{lemma} \label{L:strict neighborhood}
Let $(X,Y,P,i,j)$ be a frame and choose 
$g_1, \dots, g_n \in \Gamma(P,\calO_P)$ whose reductions
cut out $Y \setminus X$ within $Y$.
For $\lambda \in (0,1) \cap \Gamma^*$, put
\begin{align*}
U_\lambda &= ]Y[_P \, \setminus \, ]Y \setminus V(g_1,\dots,g_n)[_{P\lambda} \\
&= \{y \in \,]Y[_P: \max_i \{|g_i(y)|\} \geq \lambda \}
\end{align*}
as in Definition~\ref{D:partial tubes}.
Let $V$ be an admissible open 
subset of $]Y[_P$ containing $]X[_P$. Then $V$ is a strict neighborhood
of $]X[_P$ if and only if for any admissible affinoid $W \subseteq ]Y[_P$,
there exists $\lambda_0 \in (0,1) \cap \Gamma^*$ such that
for all $\lambda \in [\lambda_0, 1) \cap \Gamma^*$, 
$U_\lambda \cap W \subseteq V$. 
\end{lemma}

A key tool in the construction and study of isocrystals
is Berthelot's ``strong fibration theorem'' 
\cite[Th\'eor\`eme~1.3.7]{ber2}, which constructs analogues of
tubular neighborhoods (of a closed subset) in ordinary
topology.
\begin{prop}[Strong fibration theorem] \label{P:strong fibration}
Let $F' \to F$ be a morphism of frames as in
\eqref{eq:functor} with
$X' = X$ and $w = \id_X$. Let $\overline{X}$ be the closure of $X$
in $P'_Y = P' \times_P Y$, and suppose that $\overline{X} \to Y$ is proper
(e.g., if $P' \to P$ is proper). Let $\calI' \subset \calO_{P'}$ be the
defining ideal of $Y'$ within $P'$, and let $\overline{\calI}'$ be the
defining ideal of $Y'$ within $P'_Y$; suppose further that
there exist sections $t_1, \dots, t_d \in \Gamma(P',\calI')$ whose
reductions induce a basis of the conormal sheaf $\overline{\calI}'/
(\overline{\calI}')^2$ on $X$. Put
\[
P'' = P \times_{\gotho_K} \widehat{\AAA^d_{\gotho_K}}
= \Spf \calO_{P'} \langle t_1, \dots, t_d \rangle;
\]
then the morphism $\phi: P' \to P''$ defined by $t_1, \dots, t_d$ is
an isomorphism on $X$, and induces an isomorphism of some strict neighborhood
of $]X[_{P'}$ within $]Y[_{P'}$ with some strict neighborhood of
$]X[_{P''}$ within $]Y[_{P''}$.
\end{prop}

\begin{remark}
The strong fibration theorem is crucial to the independence under pullback
properties of isocrystals (Propositions~\ref{P:pullback0} 
and~\ref{P:pullback1}). It also intervenes in 
the definition
of constant/unipotent monodromy (Subsection~\ref{subsec:monodromy}).
\end{remark}

\subsection{Connections and log-connections}

\begin{convention}
When some construction is made relative to a morphism $f: V \to W$ of rigid
spaces, in case we omit mention of this morphism we take it to be the
structure morphism $f: V \to \Maxspec K$ of a rigid space $V$ over $K$.
\end{convention}

\begin{defn}
For $A$ an affinoid algebra, let $\Omega^1_{A/K}$ denote the module
of continuous 
differentials of $A$ over $K$, as in \cite[Theorem~3.6.1]{fresnel}.
Likewise, for $X$ a rigid space, let $\Omega^1_{X/K}$ denote the sheaf
of continuous differentials on $X$ over $K$; this sheaf is coherent, and is
locally free if $X$ is smooth over $K$ \cite[Theorem~3.6.3]{fresnel}.
If $f: V \to W$ is a morphism of rigid spaces, we define
$\Omega^1_{V/W} = \Omega^1_{V/K} / f^* \Omega^1_{W/K}$.
Write $\Omega^i_{V/W} = \wedge^i_{\calO_V} \Omega^1_{V/W}$.
\end{defn}

\begin{remark} \label{R:polydisc neighborhoods}
Note that if $V$ is smooth over $K$, then for any point $x \in V$,
we can find an affinoid subdomain $W$ of $V$ containing $x$ and some
$t_1, \dots, t_n \in \calO(W)$ such that $dt_1, \dots, dt_n$
freely generate $\Omega^1_{V/K}$ on $W$.
If $x$ is a $K$-rational point, we
can further ensure that $t_1, \dots, t_n$ all vanish at $x$.
For such a choice, we obtain an \'etale map $W \to \AAA^n_K$ defined by
$t_1, \dots, t_n$, sending $x$ to the origin; this map
can be shown (as is done in the proof of
\cite[Proposition~1.3]{grosse-klonne2}) to induce an isomorphism of
an affinoid subdomain of $V$ containing $x$ with some affinoid subdomain
of $\AAA^n_K$ containing the origin.
In particular, we obtain a cofinal set of affinoid subdomains of $V$
containing $x$ of the form
\[
\{y \in V: |t_i(y)| \leq \epsilon \qquad (i=1, \dots, n)\}
\]
for $\epsilon \in (0, +\infty) \cap \Gamma^*$ sufficiently small.
\end{remark}

\begin{defn} \label{defn:nablamod}
Let $f: V \to W$ be a morphism of rigid spaces.
A \emph{$\nabla$-module} on $V$, relative to $W$,
is a coherent sheaf $\calE$ of
$\calO_V$-modules on $V$, equipped with an integrable $f^{-1} \calO_W$-linear
connection $\nabla: \calE \to 
\calE \otimes_{\calO_V} \Omega^1_{V/W}$.
If $V$ is smooth over $K$, then
any $\nabla$-module on $V$ (relative to $\Maxspec K$)
is automatically locally free,
as in \cite[Proposition~2.2.3]{ber2}.
\end{defn}

One can also make a logarithmic analogue of this construction;
we will not use it again in this section, but it will become
crucially important later on.
\begin{defn} \label{defn:logdiff}
Let $f: V \to W$ be a morphism of rigid spaces,
and fix $x_1, \dots, x_m \in \Gamma(V,\calO)$.
Let $\Omega^{1,\log}_{V/W}$ be the coherent sheaf on $V$ given as the 
quotient of
\[
\Omega^1_{V/W} \oplus \calO_V s_1 \oplus \cdots \oplus \calO_V s_m
\]
by the relations $x_i s_i - dx_i$ for
$i=1, \dots, m$. We call $\Omega^{1,\log}_{V/W}$ the \emph{module
of (continuous) logarithmic differentials} with respect to the $x_i$.
\end{defn}
\begin{remark}
A better way to make this definition would be to first define
logarithmic structures on rigid spaces, then define $\Omega^{1,\log}_{V/W}$
to be the module of differentials of $V$ equipped with the log
structure generated by $x_1,\dots,x_m$, relative to $W$.
Rather than do that here, we stick to the  \emph{ad hoc} construction;
however, we will discuss logarithmic structures on schemes and formal
schemes in Section~\ref{sec:logext}.
\end{remark}

\begin{defn} \label{defn:lognablamod} 
With notation as in Definition~\ref{defn:logdiff},
a \emph{log-$\nabla$-module} on $V$ with respect to the $x_i$,
relative to $W$,
is a coherent locally free sheaf $\calE$ of
$\calO$-modules on $V$, equipped with an integrable
$f^{-1} \calO_W$-linear connection $\nabla: \calE \to 
\calE \otimes \Omega^{1,\log}_{V/W}$.
\end{defn}

We will also need the notion of horizontal sections.
\begin{defn}
With notation as in Definition~\ref{defn:lognablamod},
a section $\bv \in \Gamma(V, \calE)$ is said to be
\emph{horizontal relative to $W$} if $\nabla \bv = 0$. Let
$H^0_W(V, \calE)$ denote the set of horizontal sections relative
to $W$;
it is a $\Gamma(W, \calO)$-module. 
\end{defn}

As in the complex analytic setting, a logarithmic connection has a residue
map associated to it.
\begin{defn}
With notation as in Definition~\ref{defn:lognablamod}, note
that over the zero locus $V(x_i)$,
$\nabla$ induces an $\calO$-linear map 
from $\calE$ to $\calE \otimes \calO_V s_i$, after quotienting 
$\calE \otimes \Omega^{1,\log}_{V/W}$ by the image of
$\calE \otimes (\Omega^1_{V/W} \oplus \oplus_{j\neq i} \calO_V s_j)$
and then reducing modulo $x_i$.
Identifying $\calE \otimes \calO_V s_i$ with $\calE$
yields an $\calO$-linear endomorphism of $\calE$ over $V(x_i)$;
we call this the \emph{residue} of $\nabla$ along $V(x_i)$.
\end{defn}

\begin{remark} \label{R:saturated1}
Beware that unlike in the $\nabla$-module case, we built the locally free 
hypothesis into the definition of a log-$\nabla$-module: otherwise
we could take for instance $\calO/t\calO$ on the affine $t$-line
to be a log-$\nabla$-module with respect to $t$. By the same token, a 
log-$\nabla$-submodule $\calF$ of a log-$\nabla$-module $\calE$ need not
have locally free quotient. However, one can get around these issues by
inserting hypotheses about nilpotence of residues; see
Subsection~\ref{subsec:constant unip}.
\end{remark}

\subsection{Convergence of Taylor series}

The construction of overconvergent isocrystals can be described
in terms of a Taylor series associated to a connection; here is
a relevant constraint.

\begin{defn}
Let $X$ be an affinoid space and let $\calE$ be a coherent
$\calO$-module on $X$. For $\eta_1, \dots, \eta_n \in [0,+\infty)$,
we say a multisequence $\{\bv_I\}$
of elements of $\Gamma(X, \calE)$, indexed by $n$-tuples $I
= (i_1, \dots, i_n)$ of nonnegative
integers, is \emph{$(\eta_1,\dots,\eta_n)$-null} if
for any multisequence $\{c_I\}$ of elements of $K$ with
$|c_I| \leq \eta_1^{i_1}\cdots \eta_n^{i_n}$, 
the multisequence $\{c_I \bv_I\}$ converges to
zero in $\Gamma(X, \calE)$ (for the canonical topology induced on this module
from the affinoid topology on $\calO(X)$).
If $\eta_1 = \cdots = \eta_n = \eta$, we simply say the multisequence
is \emph{$\eta$-null}.
Note that it suffices to check the convergence on each element of an
admissible affinoid cover of $X$.
\end{defn}

\begin{defn} \label{D:eta convergent}
Let $h: V \to X$ be a morphism of affinoid spaces,
and suppose that $x_1, \dots, x_m \in \calO(V)$ 
have the property
that $dx_1, \dots, dx_m$ freely generate $\Omega^1_{V/X}$
(so that in particular the morphism $h$ is smooth).
Let $\calE$ be a $\nabla$-module over $V$ relative to $X$;
we may then view $\calE$ as being
equipped with commuting actions of the partial differential operators
$\frac{\del}{\del x_i}$ for $i=1, \dots, m$.
For $\eta \in [0,+\infty)$ and $\bv \in \Gamma(V, \calE)$, we say
$\calE$ (or its connection) is \emph{$\eta$-convergent at $\bv$}
(with respect to $x_1, \dots, x_m$) if the multisequence
\[
\frac{1}{i_1! \cdots i_m!} \frac{\del^{i_1}}{\del x_1^{i_1}}
\cdots \frac{\del^{i_m}}{\del x_m^{i_m}}
\bv
\]
is $\eta$-null; if $\calE$ is $\eta$-convergent at all $\bv 
\in \Gamma(V, \calE)$, we simply say
that $\calE$ is $\eta$-convergent.
\end{defn}

\begin{defn} \label{D:eta-admissible}
With notation as in Definition~\ref{D:eta convergent}, we say that
$x_1, \dots, x_m$ form an \emph{$\eta$-admissible coordinate system}
on $V$ (relative to $X$)
if the trivial $\nabla$-module with $\calE = \calO$ and $\nabla = d$ is
$\eta$-convergent. In this case, by the Leibniz rule, any $\calE$ is
$\eta$-convergent if and only if it is $\eta$-convergent at each of
a set of generators of $\Gamma(V,\calE)$.
\end{defn}

\begin{remark} \label{R:eta convergent}
With notation as in Definition~\ref{D:eta convergent} and~\ref{D:eta-admissible}, suppose that 
$y_1, \dots, y_m \in \Gamma(V,\calO)$ 
form another $\eta$-admissible coordinate system, and suppose that
the $m \times m$ matrix $A$ defined by $A_{ij} = \frac{\del y_i}{\del 
x_j}$ is invertible over $\Gamma(V, \gotho)$. Then
the criterion of $\eta$-convergence with respect to $y_1, \dots, y_m$
is equivalent to the criterion with respect to $x_1, \dots, x_m$.
\end{remark}

\begin{remark} \label{R:eta conv exact}
Retain notation as in Definitions~\ref{D:eta convergent} and~\ref{D:eta-admissible}.
If $0 \to \calE_1 \to \calE \to \calE_2 \to 0$ is a short
exact sequence of $\nabla$-modules, then $\calE$ is $\eta$-convergent
(with respect to a particular $\eta$-admissible coordinate system)
if and only if $\calE_1$ and $\calE_2$ are $\eta$-convergent.
In particular, the $\eta$-convergent $\nabla$-modules
on a given $V$ form an abelian category.
\end{remark}

\subsection{Overconvergent sections}

We recall the ``overconvergent sections'' functor
from \cite[2.1.1]{ber2}.
\begin{defn}
Let $(X,Y,P,i,j)$ be a frame.
For $V' \subset V$ two strict neighborhoods of $]X[$ within $]Y[$, let 
$\alpha_V$ (resp.\ $\alpha_{VV'}$) denote the open immersion of $V$ into
$]Y[$ (resp.\ of $V'$ into $V$). Given an
$\calO_V$-module $\calE$ on $V$, define
\[
j^{\dagger}_V \calE = \lim_{\rightarrow} \alpha_{VV'*} \alpha_{VV'}^* \calE,
\]
the limit taken over strict neighborhoods $V'$ of $]X[$ within $]Y[$
which are contained in $V$.
The functors
$\alpha_{VV'*}$ and $\alpha_{VV'}^*$ induce equivalences of categories
between $j^\dagger_{V} \calO$-modules and $j^\dagger_{V'} \calO$-modules.
The functor $\alpha_{V*} j^\dagger_V \alpha_V^*$ on
$\calO_{]Y[}$-modules 
does not depend on the choice of $V$, so we notate it
simply as $j^\dagger$.
\end{defn}

\begin{remark} \label{R:no dagger}
By \cite[Proposition~2.2.10]{ber2}, any coherent $j^\dagger 
\calO_{]Y[}$-module is the pullback of a coherent $\calO$-module
on a strict neighborhood of $]X[$ in $]Y[$. Moreover, if two such modules are
given, any morphism between them is obtained from a morphism between them
on a strict neighborhood where they are both defined. In practice, then,
we will write down coherent $j^\dagger \calO_{]Y[}$-modules by writing
down coherent $\calO$-modules on strict neighborhoods of $]X[$,
with the understanding that the strict neighborhood is to be shrunk
as needed.
\end{remark}

\begin{defn} \label{D:overconvergent}
Let $(X,Y,P,i,j)$ be a frame.
Let $\delta: P_K \to
P_K \times_K P_K$ be the diagonal, put $j' = \delta \circ j$,
let $\calI 
\subset \calO_{P_K \times P_K}$ be the
ideal of the image of $\delta$, and put $\calP^n = 
\calO_{P_K \times P_K}/\calI^{n+1}$.
Let $\calE$ be a coherent $j^\dagger \calO_{]Y[}$-module equipped with
an integrable $K$-linear connection $\nabla$. 
Then in the usual fashion \cite[2.2.2]{ber2},
the connection gives rise to isomorphisms
\[
\epsilon_n : j^\dagger \calP^n \otimes_{j^\dagger \calO_{]Y[}} \calE
\stackrel{\sim}{\to} \calE \otimes_{j^\dagger \calO_{]Y[}} j^\dagger \calP^n.
\]
We say $\calE$ is \emph{overconvergent along $Y \setminus X$} 
if there exists an
isomorphism $\epsilon: \pi_2^* \calE \stackrel{\sim}{\to} \pi_1^* \calE$
which induces each $\epsilon_n$ by reducing modulo $(j')^\dagger \calI^{n+1}$
and using the canonical identification $\delta^{-1} (j')^\dagger
\cong j^\dagger \delta^{-1}$ of \cite[(2.1.4.4)]{ber2}. If 
$Y \setminus X = \emptyset$, we say instead that $\calE$ is \emph{convergent}.
\end{defn}

\begin{remark}
By \cite[Proposition~2.2.3]{ber2} (as in Remark~\ref{R:no dagger}),
any coherent $j^\dagger \calO_{]Y[}$-module equipped with
an integrable $K$-linear connection $\nabla$ is the pullback of
a $\nabla$-module $\calE$ on some strict
neighborhood of $]X[$ in $]Y[$, and likewise any morphism between
such modules extends to some strict neighborhood of $]X[$ in $]Y[$.
By \cite[Proposition~2.2.6]{ber2}, the connection is overconvergent 
along $Y \setminus X$ if
and only if there exists
$\epsilon: \pi_2^* \calE \stackrel{\sim}{\to} \pi_1^* \calE$ of the desired
form over some strict neighborhood of $]X[_{P^2}$ in $]Y[_{P^2}$.
By abuse of language, we will say that ``$\calE$ is overconvergent
along $Y \setminus X$'' to mean that $j^\dagger \calE$ is overconvergent
along $Y \setminus X$.
\end{remark}

The condition of overconvergence can also be interpreted in terms of
the convergence of the Taylor series associated to the connection,
as follows.
\begin{remark}
Let $(X,Y,P,i,j)$ be a frame
and suppose that
the differentials of $x_1, \dots, x_n \in \Gamma(P, \calO_P)$
generate $\Omega^1_{P/\gotho_K}$ over a neighborhood of $X$.
Then $dx_1, \dots, dx_n$ also generate $\Omega^1_{P_K/K}$ over a strict
neighborhood of $]X[$ in $]Y[$ \cite[Proposition~2.2.13]{ber2}.
\end{remark}

By \cite[Proposition~2.2.13]{ber2},
we have the following. (Note that the statement of
\cite[Proposition~2.2.13]{ber2} only includes the equivalence between
(a) and (b) below; however, the fact that (b) holds for all sufficiently
large $\lambda$ is evident in the proof of
\cite[Proposition~2.2.13]{ber2}.)
\begin{prop} \label{P:taylor}
Let $(X,Y,P,i,j)$ be a frame;
define the sets $[Y]_\eta$ as in Definition~\ref{D:partial tubes}
(using any set of generators).
Suppose further that there exists $g \in \Gamma(P, \calO_P)$ which
cuts out $Y \setminus X$ within $Y$;
define the sets $U_\lambda$ as in Lemma~\ref{L:strict neighborhood}
using $g$.
Suppose further that
the differentials of $x_1, \dots, x_n \in \Gamma(P, \calO_P)$
generate $\Omega^1_{P/\gotho_K}$ over a neighborhood of $X$.
Let $V$ be a strict neighborhood of $]X[$ in $]Y[$,
and let $\calE$ be a $\nabla$-module on $V$. Then the following
conditions are equivalent.
\begin{enumerate}
\item[(a)]
$j^\dagger \calE$ is overconvergent.
\item[(b)]
For each $\eta \in (0,1) \cap \Gamma^*$,
there exists $\lambda \in (0,1) \cap \Gamma^*$ such that
$[Y]_\eta \cap U_\lambda \subseteq V$ and $\calE$ is $\eta$-convergent
with respect to $x_1, \dots, x_n$ over $[Y]_\eta \cap U_\lambda$.
\end{enumerate}
Moreover, if these hold, then for each $\eta \in (0,1) \cap \Gamma^*$,
the conclusion of (b) holds for all $\lambda \in (0,1) \cap \Gamma^*$
sufficiently large.
\end{prop}

\begin{remark}
Note that the trivial $\nabla$-module $\calO_V$ evidently satisfies
the definition of overconvergence given in Definition~\ref{D:overconvergent}.
Hence with conditions as in Proposition~\ref{P:taylor}, for each
$\eta \in (0,1) \cap \Gamma^*$,
$x_1, \dots, x_n$ necessarily form an $\eta$-admissible coordinate system
(in the sense of Definition~\ref{D:eta-admissible})
on $[Y]_\eta \cap U_\lambda$ for all 
$\lambda \in (0,1) \cap \Gamma^*$ sufficiently large. In particular,
the property of $\eta$-convergence may be checked at each element of
a set of generators of $\Gamma([Y]_\eta \cap U_\lambda, \calE)$.
\end{remark}

\begin{remark}
Note that the criterion for overconvergence in Proposition~\ref{P:taylor}
simplifies somewhat in case $Y = P_K$, as in that case $[Y]_\eta = P_K$
for all $\eta \in (0,1)$. This will be a great help as we work with ``small
frames'' in Section~\ref{sec:mono}.
\end{remark}

\begin{remark}
In a previous version of this paper, the restriction that
$Y \setminus X$ must be a divisor in $Y$ was omitted from 
Proposition~\ref{P:taylor}; thanks to Bernard le Stum for pointing this out.
That restriction will be harmless in practice, as we will be able to blow up
in $Y \setminus X$ without disturbing the concept of overconvergence;
see Definition~\ref{defn:invimage} below.
\end{remark}

\subsection{Isocrystals}

Given a morphism of frames as in \eqref{eq:functor},
one obtains a pullback functor $u_K^*$ from the category of
$j^\dagger \calO_{]Y[}$-modules with integrable overconvergent connection
to the analogous category of $(j')^\dagger \calO_{]Y'[}$-modules.
The key consequences of overconvergence are the following two
``homotopy invariance'' results for the pullback functors, which are
\cite[Proposition~2.2.17]{ber2} and
\cite[Th\'eor\`eme~2.3.1]{ber2}, respectively.

\begin{prop} \label{P:pullback0}
Given two morphisms of frames as in \eqref{eq:functor} 
factoring through the same map $Y' \to Y$ with
$u = u_1$ and $u = u_2$, respectively,
there is a canonical isomorphism $\epsilon_{u_1,u_2}$ between the functors
$u_{1K}^*$ and $u_{2K}^*$. Moreover, for any horizontal section $s$,
one has $\epsilon_{u_1,u_2}(u_{1K}^*(s)) = u_{2K}^*(s)$.
\end{prop}
\begin{prop} \label{P:pullback1}
Given a morphism of frames as in \eqref{eq:functor} in which
$X =X', Y=Y'$, and $v$ and $w$ are the identity maps, the functor $u_K^*$ is an
equivalence of categories.
\end{prop}

Using Proposition~\ref{P:pullback1}, one can define a category of
isocrystals. This is done somewhat informally in \cite{ber2}; a more
``crystalline'' presentation is given by ongoing work of le Stum
(see \cite{lestum} for a report, and \cite{lestum2} for further details). 
Here we take a middle road.
\begin{defn} \label{defn:site}
Given an open immersion $i: X \hookrightarrow Y$ of $k$-varieties,
define the site $\calC_{X,Y}$ as follows.
The objects of $\calC_{X,Y}$ 
are tuples $(U, P, j)$, where $U$ is an open subscheme of $Y$ and
$(X \cap U, U, P, i, j)$ form a frame.
A morphism $(U,P,j) \to (U',P',j')$ consists of an inclusion $U \subseteq
U'$ and a morphism $f: P \to P'$ of formal schemes such that $f \circ j$
equals the restriction of $j'$ to $U$. A covering $\{(U_i,P_i,j_i)
\to (U,P,j)\}$ is admissible if $\{U_i \to U\}$ is surjective.
\end{defn}

\begin{defn} \label{defn:isocrys}
With notation as in Definition~\ref{defn:site}, put $Z = Y \setminus X$.
An \emph{isocrystal on $X$ overconvergent along $Z$ (over $K$)}
is a crystal on $\calC_{X,Y}$ of 
coherent locally free $j^\dagger 
\calO$-modules with overconvergent connection: i.e.,
one specifies for each $(U,P,j) \in \calC_{X,Y}$ a coherent locally free
$j^\dagger \calO_{]U[}$-module $\calE_U$ equipped with an integrable connection
overconvergent along $U \cap Z$, and for each morphism
$u: (U,P,j) \to (U',P',j')$ an isomorphism
$\calE_U \stackrel{\sim}{\to} u^* \calE_{U'}$ of modules with connection,
such that the isomorphisms satisfy the obvious cocycle condition.
Let $\Isoc^\dagger(X,Y/K)$ denote the category of these objects.
In case $X=Y$, we call the category $\Isoc(X/K)$ and call its elements
\emph{convergent isocrystals on $X$}.
\end{defn}

\begin{defn}
For $F = (X,Y,P,i,j)$ a frame,
there is an obvious
restriction functor from isocrystals on $X$ overconvergent along 
$Y \setminus X$ to
coherent locally free $j^\dagger \calO_{]Y[}$-modules with 
integrable 
overconvergent connection; this is called the \emph{realization} functor
for the frame $F$.
Using Proposition~\ref{P:pullback1}, one can show that each realization
functor is itself an equivalence of categories. Namely,
to construct an isocrystal with any given realization,
given $j: Y \hookrightarrow P_k$ and $j': U \hookrightarrow P'_k$,
we restrict $j$ to $U$,
pull back along the first projection of $P \times P'$, then
apply Proposition~\ref{P:pullback1} to ``push forward'' along
the second projection. (One can also speak of realizations on frames
enclosing open subvarieties of $X$, but of course those will not
typically be equivalences of categories.)
\end{defn}
\begin{remark}
In fact, carrying the connection around in this construction is superfluous;
as happens for the infinitesimal and crystalline sites, the connection
data is already captured in the structure of a crystal of 
$j^\dagger \calO$-modules. This is the point of view adopted in
\cite{lestum, lestum2}. Another approach is to state the definition in terms
of simplicial schemes, as in \cite[1.3.1]{shiho1}.
\end{remark}

\begin{defn} \label{defn:invimage}
Given a diagram of the form
\[
\xymatrix{
X' \ar^{j'}@{^{(}->}[r] \ar_w[d] & Y' \ar_v[d] \ar[r] & \Spec \gotho_{K'} \ar[d]
\\
X \ar^{j}@{^{(}->}[r]  & Y \ar[r] & \Spec \gotho_{K},
}
\]
one obtains by pullback (as in \cite[2.3.2.2]{ber2}) an inverse image
functor
\[
v^*: \Isoc^\dagger(X,Y/K) \to \Isoc^\dagger(X',Y'/K').
\]
In case $X = X'$, $w = \id_X$, $K = K'$, and $v$ is proper, then $v^*$ is an
equivalence of categories \cite[Th\'eor\`eme~2.3.5]{ber2}.
In particular, if $Y$ itself is proper, then
the category $\Isoc^\dagger(X,Y/K)$ is independent of $Y$;
it is thus denoted $\Isoc^\dagger(X/K)$ and its objects are called
\emph{overconvergent isocrystals on $X$ (over $K$)}. This category
is abelian \cite[Remarques~2.3.3]{ber2}.
\end{defn}

\begin{defn}
Suppose that
\[
\xymatrix{
X' \ar^{j'}@{^{(}->}[r] \ar_w[d] & Y' \ar_v[d] \\
X \ar^{j}@{^{(}->}[r]  & Y 
}
\]
is a commutative diagram of $k$-varieties with $j,j'$ open immersions,
$v$ finite, and $w$ finite \'etale. Then one obtains a pushforward functor
\[
v_*: \Isoc^\dagger(X', Y'/K) \to \Isoc^\dagger(X,Y/K)
\]
from the pushforward along a finite \'etale morphism of rigid spaces.
As shown by Tsuzuki (see \cite[5.1]{tsuzuki-duke}), 
for $\calE, \calF \in \Isoc^\dagger(X',Y'/K)$, we have a canonical bijection
\begin{equation} \label{eq:pushforward}
\Hom(\calE, \calF) \to \Hom(v_* \calE, v_* \calF);
\end{equation}
in addition,
for $\calE \in \Isoc^\dagger(X,Y/K)$ and $\calF \in \Isoc^\dagger(X',Y'/K)$,
one has adjunction and trace morphisms
\[
\calE \stackrel{\ad}{\to} v_* v^* \calE \stackrel{\tr}{\to} \calE, \qquad
\calF \stackrel{\ad}{\to} v^* v_* \calF \stackrel{\tr}{\to} \calF
\]
such that the displayed compositions are multiplication by the degree of $v$.
(Tsuzuki explicitly constructs the first sequence; the second sequence
is obtained from the first by putting $\calE = v_* \calF$
and invoking \eqref{eq:pushforward}.)
\end{defn}

\section{Local monodromy of $p$-adic differential equations}
\label{sec:local}

We next gather some facts about differential modules
on $p$-adic annuli. Various aspects of this theory have been
treated previously, e.g., by
Crew \cite{crewfin}, Tsuzuki \cite{tsuzuki-slope}, de Jong
\cite{dejong2}, and this author \cite{me-full}. 
New features here include the systematic 
presentation in terms of rigid analytic
spaces (which obviates the need to restrict to discretely valued
or even spherically complete
coefficient fields), the treatment of multidimensional annuli,
the consideration of families of annuli (based partly on \cite{me-finite}),
and the introduction of logarithmic singularities.
However, we restrict here to cases of unipotent monodromy; we will
consider ``quasi-unipotent'' differential modules later in the series.

Throughout this section,
we retain the conventions introduced in Subsection~\ref{subsec:initial}.

\subsection{Polyannuli}

\begin{defn}
We say a subinterval $I$ of $[0, +\infty)$ is \emph{aligned} if any
endpoint at which it is closed is either equal to zero or contained
in $\Gamma^*$ (the divisible closure of the image of $|\cdot|$ on $K^*$).
In particular, any open interval is aligned, and any aligned
interval can be written as
the union of a weakly increasing sequence of aligned closed subintervals.
We say $I$ is \emph{quasi-open} if it is open at each 
\emph{nonzero} endpoint, i.e., it is of one of the forms $(a,b)$ or
$[0,b)$; any quasi-open interval is aligned.
\end{defn}

\begin{defn}
For $I$ an aligned subinterval of $[0, +\infty)$, 
we define the \emph{polyannulus} $A^n_K(I)$ as
\[
A^n_K(I) = \{(t_1, \dots, t_n) \in \AAA^n_K: |t_i| \in I \quad
(i=1, \dots, n) \}.
\]
\end{defn}

\begin{convention}
In the notation $A^n_K(I)$, we drop the parentheses
around the interval $I$ if it is being written out explicitly, e.g.,
we write $A^n_K[0,1)$ instead of $A^n_K([0,1))$.
\end{convention}

\begin{remark} \label{R:punctured polydisc}
Note that if $0 \notin I$ and $n > 1$, then $A^n_K(I)$ is not the same
as a punctured polydisc; if $I = J_1 \setminus J_2$, where $J_1$ and
$J_2$ are aligned intervals both containing 0, the latter would be
\[
A^n_K(J_1) \setminus A^n_K(J_2),
\]
which unlike $A^n_K(I)$ is not an affinoid space.
\end{remark}

\begin{defn}
For $X$ an affinoid space and $I$ an aligned subinterval of $[0,+\infty)$, the ring
$\Gamma(X \times A^n_K(I), \calO)$ consists of
Laurent series
\[
\sum_{J \in
\ZZ^n} c_J t^J = \sum_{J = (j_1, \dots, j_n)} c_J t_1^{j_1} \cdots t_n^{j_n}
\]
with coefficients in $\Gamma(X, \calO)$, 
such that $|c_J|_X \rho_1^{j_1}\cdots \rho_n^{j_n} \to 0$ as $J \to \infty$
(that is, $|c_J|_X \rho_1^{j_1} \cdots \rho_n^{j_n}$ exceeds any
particular positive number for only finitely many $J$) for each
$\rho_1, \dots, \rho_n \in I$.
For $R = (r_1, \dots, r_n) \in I^n$,
let $|\cdot|_{X,R}$ denote the function on $\calO(X \times A^n_K(I))$
given by
\[
\left|\sum_J c_J t^J \right|_{X,R}
= \sup_J \{ |c_J|_X r_1^{j_1} \cdots r_n^{j_n}\};
\]
note that the supremum is achieved by at least one, but only
finitely many, tuples $(j_1,\dots,j_n)$. If $R = (r, \dots, r)$,
we also write $|\cdot|_{X,r}$ for $|\cdot|_{X,R}$.
\end{defn}

One has analogues of the maximum modulus principle and the Hadamard 
three circles theorem for $|\cdot|_{X,R}$.
\begin{lemma} \label{L:Hadamard}
Let $X$ be an affinoid space.
\begin{enumerate}
\item[(a)]
For $x \in \Gamma(X \times A^n_K[0,b], \calO)$ with $b \in [0,+\infty) 
\cap \Gamma^*$,
and $R \in [0,b]^n$, we have $|x|_{X,R} \leq |x|_{X,b}$.
\item[(b)]
For $x \in \Gamma(X \times A^n_K(I), \calO)$ with $I$ an aligned
subinterval of $[0,+\infty)$, $A,B \in I^n$, and
$c \in [0,1]$, put $r_i = a_i^{c} b_i^{1-c}$; then
$|x|_{X,R} \leq |x|_{X,A}^c |x|_{X,B}^{1-c}$.
\end{enumerate}
\end{lemma}
\begin{proof}
\begin{enumerate}
\item[(a)]
If $x = \sum c_J t^J \in \Gamma(X \times A^n_K[0,b], \calO)$, 
then $c_J = 0$ unless $j_1, \dots, j_n \geq 0$.
Hence if $R \in [0,b]^n$, then
\[
|c_J|_X r_1^{j_1} \cdots r_n^{j_n} \leq |c_J|_X b^{j_1 + \cdots + j_n};
\]
taking suprema yields $|x|_{X,R} \leq |x|_{X,b}$.
\item[(b)]
Note that the desired inequality holds with equality if
$x = c_J t^J$ is a monomial.
For a general $x = \sum c_J t^J$, we then have
\begin{align*}
|x|_{X,R} &= \sup_J \{ |c_J t^J|_{X,R} \} \\
&= \sup_J \{ |c_J t^J|_{X,A}^c |c_J t^J|_{X,B}^{1-c} \} \\
&\leq \sup_J \{ |c_J t^J|_{X,A} \}^c \sup_J \{|c_J t^J|_{X,B} \}^{1-c} \\
&= |x|_{X,A}^c |x|_{X,B}^{1-c},
\end{align*}
as desired.
\end{enumerate}
\end{proof}
\begin{cor} \label{C:maximum}
For $x \in \calO(X \times A^n_K[a,b])$, the maximum of
$|x|_{X,R}$ over all $R \in [a,b]^n$ is achieved by a tuple
$R \in \{a,b\}^n$.
\end{cor}

\begin{lemma}
For $X$ an affinoid space, $I$ an aligned subinterval of $[0,+\infty)$,
and $A =(a_i) \in I^n \cap (\Gamma^*)^n$,
the norm $|\cdot|_{X,A}$ 
coincides with the supremum seminorm on 
the affinoid space
\[
X \times \{(x_1, \dots, x_n) \in A^n_K(I): |x_i| = a_i \quad (i=1, \dots, n)\}.
\]
\end{lemma}
\begin{proof}
There is no loss of generality in enlarging $K$ so that $a_1, \dots, a_n$ 
land in the image of $|\cdot|$ itself.
Given $\sum c_J t^J \in \calO(X \times A^n_K(I))$, 
the supremum defining $|\sum c_J t^J|_{X,A}$
is achieved by finitely
many tuples $J$. Let $S$ be the set of these tuples;
by enlarging $K$ again, we can ensure that there
exist $x_1, \dots, x_n \in K$, with $|x_i| = a_i$ for each $i$,
such that the evaluation of $\sum_{J \in S} c_J t^J$ at $t_i = x_i$ has norm
equal to $|\sum c_J t^J|_{X,A}$.
\end{proof}
\begin{cor}
For $[a,b]$ aligned,
the affinoid topology on $\calO(X \times A^n_K[0,b])$
coincides with the subspace topology induced by
the affinoid topology on $\calO(X \times A^n_K[a,b])$.
\end{cor}
\begin{cor} \label{C:quasi-Stein}
For any aligned subinterval $I$ of $[0,+\infty)$,
the space $A^n_K(I)$ is a quasi-Stein space. (In particular, if $X$ is
also quasi-Stein, then so is $X \times A^n_K(I)$.)
\end{cor}
\begin{proof}
Let $I_1 \subseteq I_2 \subseteq \cdots$ be a weakly increasing sequence of
closed aligned intervals with union $I$. Then $A^n_K(I)$ is the
union of the $A^n_K(I_j)$; moreover, if $0 \in I$, then the
polynomial ring
$K [t_1,\dots,t_n]$ is dense in each
$\calO(A^n_K(I_j))$, since the Laurent series
$\sum c_J t^J$ is the limit under each $|\cdot|_{X,R}$ of its finite 
partial sums.
By the same token, if $0 \notin I$, then the Laurent polynomial
ring $K [t_1, \dots, t_n, t_1^{-1}, \dots, t_n^{-1}]$ is dense in 
each $\calO(A^n_K(I_j))$. In either case, $A^n_K(I)$ is quasi-Stein.
\end{proof}

We will need a refinement of the argument of Corollary~\ref{C:quasi-Stein}
for quasi-open intervals.
\begin{lemma} \label{L:truncation}
Let $I$ be a quasi-open subinterval of $[0,+\infty)$, and let
$x = \sum_J c_J t^J$
be an element of $\calO(A^n_K(I))$. For $l = 1,2, \dots$, put
\[
x_l = \sum_{J: |j_1|, \dots, |j_n| \geq l} c_J t^J.
\]
Then for any $R \in I^n$, there exists $\eta > 1$ such that
$\lim_{l \to \infty} \eta^l |x_l|_{K,R} = 0$.
\end{lemma}
\begin{proof}
First suppose $I = (a,b)$.
Pick $a', b' \in \Gamma^*$ with $a < a' < r_i < b' < b$ for $i=1,\dots,n$;
then the supremum seminorm of $x_l$ on $A^n_K[a',b']$ tends to 0 as
$l \to \infty$. That is,
\[
\lim_{l \to \infty} \max_{S \in \{a',b'\}^n} \{|x_l|_{K,S}\} = 0.
\]
Given $J$, put
\[
s_i = \begin{cases} b' & j_i \geq 0 \\ a' & j_i < 0 \end{cases}
\]
for $i=1, \dots, n$, and put $S = (s_1, \dots, s_n)$. Then
\[
|c_J t^J|_{X,R} \leq |c_J t^J|_{X,S} \prod_{i=1}^n \max\{a'/r_i, r_i/b'\}^{|j_i|}.
\]
We may thus take $\eta = \prod_{i=1}^n \min\{r_i/a', b'/r_i\} > 1$.

In case $I = [0,b)$, the argument is similar but easier: for
any $b'$ with $r_i < b' < b$ for $i=1,\dots, n$, we have
\[
|c_J t^J|_{X,R} \leq |c_J t^J|_{X,S} \prod_{i=1}^n (r_i/b')^{j_i}.
\]
and so we may take $\eta = \prod_{i=1}^n (b'/r_i) > 1$.
\end{proof}

\subsection{Constant and unipotent connections}
\label{subsec:constant unip}

We now start considering $\nabla$-modules on the product 
of a smooth rigid space with a polyannulus.
For convenience, we encapsulate a running hypothesis.

\begin{hypo} \label{H:setup conn}
Let $f: V \to W$ be a morphism of smooth rigid spaces,
and suppose $x_1, \dots, x_m \in \Gamma(V, \calO)$ 
have zero loci which are
smooth and meet transversely. 
\end{hypo}

\begin{defn} \label{D:LNM}
Under Hypothesis~\ref{H:setup conn},
let $n$ be a positive integer,
and let $X$ be an admissible open subset of $V \times A^n_K[0,1)$.
Define the category $\LNM_{X/W}$ to be the category of
log-$\nabla$-modules over $X$ relative to $W$ 
with respect to $t_1, \dots, t_n,x_1,\dots,x_m$, having nilpotent residues.
\end{defn}

\begin{remark} \label{R:saturated}
In Definition~\ref{D:LNM},
we omit $W$ in case $f$ coincides with the structural morphism
$V \to \Maxspec K$ (which we will more briefly describe hereafter by
saying ``if $W=\Maxspec K$''). If we are in this case and
$m=n=0$, then $\LNM_{X}$ is an abelian category;
this will also turn out to be true if $m>0$ or $n>0$, by virtue of
Lemma~\ref{L:locally free}.
\end{remark}

\begin{convention}
If $I = [0,0]$, we
will regard $\Omega^{1,\log}$ over $A^n_K[0,0]$ as being freely generated
by $\frac{dt_1}{t_1}, \dots, \frac{dt_n}{t_n}$. That is, 
for $U = V \times A^n_K[0,0]$, the elements of
$\LNM_{U/W}$ will be 
log-$\nabla$-modules over $V$ relative to $W$ with respect to
$x_1,\dots, x_m$, equipped with
$n$ commuting  endomorphisms, which for consistency we view
as the actions of the operators $\del_i = t_i \frac{\del}{\del t_i}$ for
$i=1, \dots, n$.
\end{convention}

\begin{defn} \label{D:constant unipotent}
Under Hypothesis~\ref{H:setup conn},
take $\calE \in \LNM_{X/W}$ for 
$X = V \times A^n_K(I)$.
We say that $\calE$ is \emph{constant} if $\calE \cong \pi_1^* \calF$
for some log-$\nabla$-module $\calF$ on $V$ relative to $W$ with
respect to $x_1, \dots, x_m$ (necessarily having nilpotent residues).
Note that if $\calE$ is constant, then for any affinoid subspace $U$ of $V$,
the restriction of $\calE$ to $X \cap (U \times A^n_K[0,1))$
is spanned by finitely many sections which are horizontal
relative to $V$.
We say $\calE$ is \emph{unipotent} if
$\calE$ admits an exhaustive filtration
by log-$\nabla$-submodules 
whose successive quotients are constant; we call such
a filtration a \emph{unipotent filtration}.
Let $\ULNM_{X/W}$ be the subcategory of $\LNM_{X/W}$ consisting of
unipotent objects.
\end{defn}

We will ultimately see (Theorem~\ref{T:unipotent equiv}) that
the following construction produces all unipotent $\nabla$-modules.
\begin{defn} \label{D:unipotent functor}
Under Hypothesis~\ref{H:setup conn},
let $I$ be an aligned subinterval of $[0,+\infty)$.
We define the functor
\[
\calU_I: \LNM_{V \times A^n_K[0,0]/W} \to \LNM_{V \times A^n_K(I)/W}
\]
as follows.
Given a log-$\nabla$-module $\calE$ over $V$ relative to $W$
with respect to $x_1, \dots, x_m$,
equipped with $n$
commuting nilpotent endomorphisms $N_1, \dots, N_n$,
define $\calU_I(\calE)$ to be the sheaf $\pi_1^* \calE$
equipped with the connection
\[
\bv \mapsto \pi_1^*(\nabla) \bv + \sum_{i=1}^n \pi_1^*(N_i)(\bv)
\otimes \frac{dt_i}{t_i}.
\]
This connection is integrable because the $\pi_1^*(N_i)$ commute
with each other and with the action of the connection on the base.
\end{defn}

\begin{remark} \label{R:const unip1}
Suppose $V = W = \Maxspec K$.
In Definition~\ref{D:constant unipotent},
if $X = V \times A^n_K[0,0]$, then
$\calE$ is constant if and only
if the $\del_i$ all act via the zero map, and
$\calE$ is unipotent if and only if the $\del_i$
are all nilpotent. 
In Definition~\ref{D:unipotent functor},
any nilpotent filtration of $\calE$ with respect to 
$N_1,\dots, N_n$ lifts to a unipotent filtration of
$\calU_I(\calE)$, so $\calU_I(\calE)$ is unipotent. 
We will generalize this remark later (Remark~\ref{R:const unip2}),
but beware that it is not true for $V,W$ general.
\end{remark}

In ordinary analysis, $\nabla$-modules on open 
polydiscs are automatically constant,
but this fails in rigid analysis without extra hypotheses; see
Remark~\ref{R:not unipotent}. However, one can at least salvage the following
result.

\begin{lemma} \label{L:small unip1}
Under Hypothesis~\ref{H:setup conn},
suppose that $V$ is affinoid, that $X = V \times A^1_K[0,a]$, and that
$\calE \in \LNM_{X/V}$ is such that
the restriction of $\calE$ to $V \times \{0\}$ is free.
Then there exists $b \in (0,a] 
\cap \Gamma^*$ such that the restriction of $\calE$ to $V \times A^1_K[0,b]$
is in the essential image of
the functor $\calU_{[0,b]}$. 
In particular, if the residue of $\calE$ along $V \times \{0\} = V(t_1)$
vanishes, then the restriction of $\calE$ to $V \times A^1_K[0,b]$ is constant.
\end{lemma}
\begin{proof}
Choose elements $\be_1, \dots, \be_n$ of $\Gamma(X, \calE)$
restricting to a basis of $\calE$ on $V \times \{0\}$.
The locus where these sections fail to be linearly independent
or fail to span $\calE$ is a closed
analytic subspace of $X$ not meeting $V$;
by the maximum modulus principle, the values of $t_1$ on this subspace
are bounded away from zero. Hence by making $a$ smaller, we can ensure that 
$\be_1, \dots, \be_n$ form a basis of $\Gamma(X, \calE)$.

Define the $n \times n$ matrix $N$ over $\calO(X)$ by the formula
\[
\del_1 \be_l = \sum_j N_{jl} \be_j
\]
and formally write $N = \sum_{i=0}^\infty N_i t_1^i$, where each $N_i$
is an $n \times n$ matrix over $\calO(V)$.
We now verify that there is a unique $n \times n$
matrix $M$ over $\calO(V)
\llbracket t_1 \rrbracket$, congruent to the identity matrix modulo $t_1$,
such that $N M + \del_1 M = M N_0$.
Namely, if we write $M = \sum_{i=0}^\infty M_i t_1^i$, for each $i>0$ we then
have
\begin{equation} \label{eq:recursion}
iM_i + N_0 M_i - M_i N_0 = -\sum_{j=0}^{i-1} N_{i-j} M_j.
\end{equation}
Let $e$ be the nilpotency index of $N_0$; then the
map $g$ on the space of $n \times n$ matrices over $\calO(V)$
defined by $g(M_i)= N_0 M_i - M_i N_0$ is itself nilpotent
of index at most $2e-1$.
The map $M_i \mapsto iM_i + N_0 M_i - M_i N_0$ on $n \times n$ matrices
is then the sum of an invertible linear map and a nilpotent linear map,
hence is invertible. Thus
$M_i$ is uniquely determined by $M_0, \dots, M_{i-1}$, proving the existence
and uniqueness of $M$.

We now analyze \eqref{eq:recursion} to show that $M$ converges on
$V \times A^1_K[0,b]$ for some $b$. 
Put $\delta = \max\{1, |N|_{X}\}$; then for all $i$,
\[
|N_i|_V  \leq \delta a^{-i}.
\]
In particular $|N_0|_V \leq \delta$.
We now prove by induction that
\begin{equation} \label{eq:recursion2}
|M_i|_V \leq |i!|^{-2e} a^{-i} \delta^{2ei}.
\end{equation}
For $i=0$, this is merely $1 \leq 1$. Given the result
for all $j<i$, examining the right side of \eqref{eq:recursion} yields
the bound
\begin{equation} \label{eq:recursion3}
|iM_i + N_0 M_i - M_i N_0|_V \leq
|(i-1)!|^{-2e} a^{-i} \delta^{2e(i-1)+1}.
\end{equation}
If $W = iM_i + N_0 M_i - M_i N_0$, we can then write
\[
M_i = \sum_{j=0}^{2e-1} (-1)^j i^{-j-1} g^{(j)}(W),
\]
where $g^{(j)}$ denotes the $j$-fold composition.
In particular, we have
\[
|M_i|_V \leq |i|^{-2e} \delta^{2e-1} |W|_V,
\]
which combines with \eqref{eq:recursion3} to yield \eqref{eq:recursion2}.

Finally, choose $b \in (0,1) \cap \Gamma^*$ with
\[
b < |p|^{2e/(p-1)} a \delta^{-2e}.
\]
By virtue of \eqref{eq:recursion2} and the inequality
$|i!| \leq |p|^{i/(p-1)}$, we have $|M_i|_V b^i < 1$ for all $i>0$.
Hence the matrix $M$ gives rise to an invertible matrix over 
$\calO(V \times A^1_K[0,b])$. 
Define the vectors $\bv_1, 
\dots, \bv_m \in \Gamma(V \times A^1_K[0,b], \calE)$ by
\[
\bv_l = \sum_j M_{jl} \be_j;
\]
then we can write $\calE = \calU_{[0,b]}(\calF)$ for $\calF$ equal to
the $\calO(V)$-span of $\bv_1, \dots, \bv_n$.
\end{proof}

\begin{lemma} \label{L:small unip}
Let $\calE$ be a $\nabla$-module (resp.\ a 
log-$\nabla$-module with nilpotent residues) on $A^n_K[0,a]$
for some $a \in (0,+\infty) \cap \Gamma^*$.
 Then there exists $b \in (0,a] 
\cap \Gamma^*$ such that $\calE$ is constant (resp.\ unipotent) on $A^n_K[0,b]$.
\end{lemma}
\begin{proof}
We proceed by induction on $n$, with vacuous base case $n=0$. Identify
$V = A^{n-1}_K[0,a]$ with the zero locus of $t_n$ in $A^n_K[0,a]$; by the
induction hypothesis, by making $a$ smaller, we can ensure that
the restriction of  $\calE$ to $V$
is constant (resp.\ unipotent). In particular, we can 
choose sections $\be_1, \dots, \be_m \in \Gamma(A^n_K[0,a], \calE)$ 
restricting to sections on $V$ which form
a basis of $\Gamma(V, \calE)$ on which
the $\del_i$ act trivially (resp.\ act via commuting nilpotent matrices
over $K$). The locus where these sections fail to be linearly independent
or fail to span $\calE$ is a closed
analytic subspace  of $A^n_K[0,a]$ not meeting $V$;
by the maximum modulus principle, the values of $t_n$ on this subspace
are bounded away from zero. Hence by making $a$ smaller, we can ensure 
that in fact
$\be_1, \dots, \be_m$ form a basis of $\Gamma(A^n_K[0,a], \calE)$.
We may then apply Lemma~\ref{L:small unip1} to see that for some
$b$, the restriction
of $\calE$ to $A^n_K[0,b]$ can be pulled back from $A^{n-1}_K[0,b]$.
By the induction hypothesis, $\calE$ is in fact constant (resp.\ unipotent).
\end{proof}

Note that Lemma~\ref{L:small unip} has important consequences for 
connections on arbitrary smooth rigid spaces: it gives us a ``very local''
criterion for checking local freeness of a module equipped with a logarithmic
connection. (Here ``very local'' means that the criterion
 can be checked in an affinoid
neighborhood around each point, not just on an admissible affinoid covering.)
\begin{lemma} \label{L:locally free1}
Let $X$ be a rigid space, and let $\calE$ be a coherent sheaf on $X$.
Then $\calE$ is locally free if and only if for each $x \in X$, there is an
affinoid neighborhood of $x$ on which $\calE$ is free.
\end{lemma}
\begin{proof}
By passing to an admissible affinoid cover, it suffices to check this in case
$X = \Maxspec A$ is affinoid. In that case, by 
Kiehl's theorem \cite[Theorem~9.4.3/3]{bgr}, $M = \Gamma(X,\calE)$
is a finitely generated $A$-module and $\calE$ is the coherent sheaf on $X$
associated to $M$. Let $Y$ denote the scheme $\Spec A$; then for each $x
\in X$, we may also regard $x$ as a point of $Y$, and the local ring
$\calO_{X,x}$ is flat over the local ring $\calO_{Y,x}$ 
because both have the same completion \cite[Proposition~7.3.2/3]{bgr}.
Hence the coherent sheaf on $Y$ associated to $M$ has free stalks
at each maximal ideal, and so $M$ is locally free.
\end{proof}

\begin{lemma} \label{L:locally free}
Let $X$ be a smooth rigid space, and suppose the zero loci of $t_1,
\dots, t_n \in \calO(X)$ are smooth and meet transversely. Let
$f: \calE \to \calF$ be a morphism of log-$\nabla$-modules with nilpotent 
residues on $X$ with respect to $t_1, \dots, t_n$.
Then the kernel and cokernel of $f$ are also log-$\nabla$-modules with
nilpotent residues.
\end{lemma}
\begin{proof}
By Lemma~\ref{L:locally free1}, it suffices to check
the local freeness pointwise; clearly the same is true of the nilpotence
of residues. Moreover, there is no harm in enlarging $K$ before checking
these conditions at a given point.
It thus suffices to check that if $x \in X$ is a $K$-rational point, then
the kernel and cokernel of $f$ have free stalks at $x$ and have nilpotent
residues there. 
There is no harm in assuming that $t_1, \dots, t_n$ vanish at $x$
and generate $dt_1, \dots, dt_n$ there. (To get to this case, first drop
the $t_i$ which do not vanish at $x$, then add back additional ones to
fill out a local coordinate system.) That done, by
Remark~\ref{R:polydisc neighborhoods}, we may assume that in fact $X
= A^n_K[0,a]$ for some $a \in (0,+\infty) \cap \Gamma^*$;
then by Lemma~\ref{L:small unip}, we may assume that
$\calE$ and $\calF$ are unipotent.

We proceed by induction on the rank of $\calE \oplus \calF$.
Choose bases $\be_1, \dots, \be_l$ and $\bbf_1, \dots, \bbf_m$
of $\calE$ and $\calF$, respectively, on which each $\del_i$
acts via a nilpotent matrix over $K$. In particular, $\nabla(\be_1) = 0$,
and so $\nabla(f(\be_1)) = 0$. By a formal power series calculation, each
element of the kernel of $\nabla$ belongs to the $K$-span of 
$\bbf_1, \dots, \bbf_m$.
Hence either $f(\be_1) = 0$, or $f(\be_1)$ generates a direct summand of
$\calF$. Quotienting by the spans of $\be_1$ and $f(\be_1)$ and repeating
the argument, we deduce that the kernel and 
image of $f$ are free with nilpotent residues, as desired.
\end{proof}

\begin{remark}
The local freeness in Lemma~\ref{L:locally free} can also be proved
on the level of completed local rings, which does not require the
use of Lemma~\ref{L:small unip}. However, Lemma~\ref{L:small unip} will
come in handy later; see Proposition~\ref{P:saturation}.
\end{remark}

\begin{remark} \label{R:const unip2}
We can now generalize both assertions of Remark~\ref{R:const unip1}
to the case $W =\Maxspec K$ and $V$ arbitrary; it suffices to treat the
first of them, and this can be done as follows.
For $\calE \in \LNM_X$ with $X = V \times A^n_K[0,0]$, 
the map $\del_n: \calE \to \calE$ is a morphism in 
$\LNM_{V \times A^{n-1}_K[0,0]}$, so its kernel is an object in that category.
Repeating the argument, we find that 
$\calE_1 = \cap_i \ker(\del_i)$ is an object in $\LNM_V$, which is nonzero
because the $\del_i$ are nilpotent. By Lemma~\ref{L:locally free},
$\calE/\calE_1 \in \LNM_X$, so we may repeat to conclude that $\calE \in
\ULNM_X$.
\end{remark}

\begin{lemma} \label{L:local filt}
Let $X$ be a smooth rigid space, and suppose the zero loci of $t_1,
\dots, t_n \in \calO(X)$ are smooth and meet transversely. Let
$\calE$ be a coherent $\calO_X$-module equipped with an integrable log-connection
with respect to $t_1, \dots, t_n$. 
Then the following conditions are equivalent.
\begin{enumerate}
\item[(a)]
$\calE$ is locally free (i.e., is a log-$\nabla$-module) and has nilpotent
residues.
\item[(b)]
For each point $x \in X$,
there is an affinoid subdomain of $X$ containing $x$, on which
$\calE$ admits a filtration whose
successive quotients are $\nabla$-modules.
\item[(c)]
For each point $x \in X$,
there is an affinoid subdomain of $X$ containing $x$, on which
$\calE$ admits a filtration whose
successive quotients are \emph{trivial} $\nabla$-modules.
\end{enumerate}
\end{lemma}
\begin{proof}
Note that (c) implies (b) trivially, and (b) implies (a) by
Lemma~\ref{L:locally free1}.
It thus remains to show that (a) implies (c).

Given (a), pick $x \in X$, and let $K'$ be a finite 
Galois extension of $K$ containing
the residue field of $x$. By shrinking $X$, we may further assume 
that
$dt_1, \dots, dt_n$ form a basis of $\Omega^1_{X/K}$ in a neighborhood
of $X$, and that $t_1, \dots, t_n$ all vanish at $x$. 
Then by applying Remark~\ref{R:polydisc neighborhoods} and shrinking
$X$ further, we may reduce to
the case where $X$ is a polydisc, in which case
Lemma~\ref{L:small unip} yields the claim over $K'$. Since the filtration
of $\calE$ can be chosen canonically (by taking the first step to be the
span of all horizontal sections, and so on), it descends from $K'$ to $K$.
Each successive quotient of the result is locally free by Lemma~\ref{L:locally
free}, and becomes trivial over $K'$, hence
is also trivial over $K$: given a spanning set of horizontal sections
defined over $K'$, we can decompose over a basis for $K'$ over
$K$ to get a spanning set of horizontal sections defined over $K$.
\end{proof}

\begin{remark}
Note that the properties of being constant/unipotent are stable under
formation of direct sums, tensor products, and duals; the property of
being unipotent is also stable under extensions. When working with
$\ULNM_X$ (i.e., with $W = \Maxspec K$), one can say more, as follows.
\end{remark}

\begin{lemma} \label{L:h0 lin ind}
Let $\calE$ be a $\nabla$-module over $A^1_K(I)$ for some closed 
aligned interval
$I$. If $\bv_1, \dots, \bv_n \in H^0(A^1_K(I), \calE)$ are linearly
independent over $K$, then they are linearly independent over
$\calO(A^1_K(I))$.
\end{lemma}
\begin{proof}
Suppose the contrary; choose a counterexample with $n$ minimal.
Take $c_1 \bv_1 +\cdots + c_n \bv_n = 0$ 
with $c_1, \dots, c_n \in \calO(A^1_K(I))$; then $c_1, \dots, c_n$ are
all nonzero. Since $\calO(A^1_K(I))$ is a principal ideal domain, we
may divide through to ensure that $c_1, \dots, c_n$ generate the unit
ideal; then they are uniquely determined up to a unit in $\calO(A^1_K(I))$.

Now observe that 
\[
\frac{\del c_1}{\del t_1} \bv_1 + \cdots +
\frac{\del c_n}{\del t_1} \bv_n = 0;
\]
consequently,
$\frac{\del c_1}{\del t_1}, \dots, \frac{\del c_n}{\del t_1}$
must equal $c_1, \dots, c_n$ times an element of $\calO(A^1_K(I))$.
If $c_1$ vanishes anywhere on $A^1_K(I)$, then 
$\frac{\del c_i}{\del t_1}$ vanishes there to lower order,
yielding a contradiction. 

Hence $c_1$ is a unit in $A^1_K(I)$, so we could have
taken $c_1 = 1$ to begin with. But in that case, 
$\frac{\del c_1}{\del t_1}, \dots, \frac{\del c_n}{\del t_1}$
would all vanish, yielding $c_1, \dots, c_n \in K$ 
and contradicting the linear independence of 
$\bv_1, \dots, \bv_n$ over $K$. This proves the claim.
\end{proof}

\begin{prop} \label{P:abelian category}
For any smooth rigid space $X$ over $K$, 
$\LNM_X$ is an abelian tensor category.
If $X = V \times A^n_K(I)$ for $I$ a closed aligned interval, then
$\ULNM_X$ is an abelian tensor subcategory of $\LNM_X$.
\end{prop}
\begin{proof}
The fact that $\LNM_X$ is an abelian category follows from 
Lemma~\ref{L:locally free}. To show that $\ULNM_X$ is an abelian tensor
subcategory, we must check that the property of being constant/unipotent
is preserved by formation of subobjects and quotients within 
$\LNM_X$; in any given situation, it suffices to check one of subobjects
or quotients, as the other will follow by dualizing.
There is no harm in enlarging $K$, so we may assume that there exists
a section $x: V \to X$ of the projection $\pi_1: X \to V$; we will write
$x$ also to mean the image of $x$.

We first check that for $n=1$ and $V = \Maxspec K$,
the property of being constant is stable under taking
quotients.
Suppose $\calE \in \LNM_X$ is constant and $g: \calE \to \calE'$ is a surjection
in $\LNM_X$. Let $\calF'$ be the image of $H^0(X, \calE)$ in $\calE'$;
then the map $\calF' \otimes_K \calO_X \to \calE'$ is surjective by
construction, and injective by Lemma~\ref{L:h0 lin ind}.
Thus $\calE'$ is constant.

We next check that for any $n$ and $V$,
the property of being constant is stable under taking
subobjects.
By induction on $n$, it suffices to check the case $n=1$ for arbitrary $V$.
Let $\calE \in \LNM_X$ be constant, so that there exists $\calF \in \LNM_V$
with $\pi_1^* \calF \cong \calE$; we can identify $\calF$ inside
$\calE$ as the $\pi_1^{-1} \calO_V$-span of the horizontal sections.
Let $g: \calE' \to \calE$ be an injection
in $\LNM_X$. Let $\calF' \in \LNM_V$ be the image of the restriction of 
$g$ to $x$.

We wish to show that $\calE' = \pi_1^* \calF'$, that is, that
the maps $\calE' \to \calE/(\pi_1^* \calF')$ 
and $\pi_1^* \calF' \to \calE/\calE'$ are zero.
Since this is a property that can be checked pointwise on $V$, it is certainly
enough to check on a polydisc around each point of $V$. If $V$ is itself a
polydisc, we may pass to its generic point and check there. We may thus
reduce to the case where $V$ is a point, where we already know that
$\calE$ is constant, and the equality $\calE' = \pi_1^* \calF'$ is thus
straightforward. Hence $\calE' = \pi_1^* \calF'$ in general, so
$\calE'$ is constant.

To conclude, we observe that 
the property of being unipotent is also stable under
subobjects and quotients: we may intersect a unipotent filtration with
a subobject or project it onto a quotient, and the successive quotients
will be constant by the previous paragraph. Hence
$\ULNM_X$ is indeed an abelian tensor subcategory.
\end{proof}

\begin{remark} \label{R:Tannakian}
One could in principle consider the Tannakian category consisting
of the log-$\nabla$-modules over $V \times A^n_K(I)$, 
and reinterpret
the constant/unipotent property for a given log-$\nabla$-module
in terms of the action of the fundamental group 
on that module. However, in order to produce a fibre functor by consideration 
of horizontal sections, it is necessary to restrict to modules with a Frobenius
structure and invoke a suitable form of the $p$-adic local monodromy theorem.
We may address this point in a subsequent paper.
\end{remark}

\subsection{Classification of unipotent log-$\nabla$-modules}

Our next
goal is to give a characterization of unipotent log-$\nabla$-modules
analogous to the pullback definition of constant log-$\nabla$-modules.
For this we will need a relative analogue of
\cite[Proposition~1.1.2]{chiar-lestum}, whose proof is
straightforward.
\begin{lemma} \label{L:ext1}
Under Hypothesis~\ref{H:setup conn}, let $I$ be an
aligned subinterval of $[0,+\infty)$,
put $X = V \times A^n_K(I)$,
and let $\calD_{X/W}$ be the noncommutative ring sheaf of (finite order)
$\calO_W$-linear log-differential operators on $X$.
Let $\calE$ and $\calE'$ be (left) $\calD_{X/W}$-modules on $X$,
with $\calE$ coherent and flat over $\calO_X$.  Then there is a natural
isomorphism
\[
\Ext^i_{\calD_{X/W}}(\calE, \calE') \cong 
\HH^i(X, \calE^\dual \otimes \calE' \otimes \Omega^{\cdot, \log}_{X/W}),
\]
where $\HH$ denotes hypercohomology.
\end{lemma}

Using Lemma~\ref{L:ext1}, we can show that $\calU_I$ commutes with the 
formation of Yoneda $\Ext$ groups.
\begin{lemma} \label{L:ext2}
Under Hypothesis~\ref{H:setup conn}, 
let $I$ be a quasi-open subinterval of $[0,+\infty)$.
Then for any $\calE, \calE' \in \ULNM_{V \times A^n_K[0,0]/W}$, the natural map
\[
\Ext^i(\calE, \calE') \to \Ext^i(\calU_I(\calE), \calU_I(\calE'))
\]
is a bijection.
\end{lemma}
\begin{proof}
If
\[
0 \to \calE_1 \to \calE \to \calE_2 \to 0
\]
is a short exact sequence, then by the long exact sequence for Yoneda
Exts and the five lemma (and the fact that the map in question is functorial),
we can reduce the question of bijectivity from the case of $\calE$ and
$\calE'$ to the cases of $\calE_1$ and $\calE'$, and of $\calE_2$ and
$\calE'$. Of course one has an analogous reduction given a short exact
sequence with $\calE'$ in the middle. We may thus reduce to
the case where $\calE$ and $\calE'$ are constant.

Next, we observe that it suffices to check the case where $W$ is affinoid,
as we may deduce the general case by making an admissible affinoid cover 
of $W$ and using the spectral sequence provided by the corresponding
\v{C}ech complex. Similarly, we may reduce to the case where $V$ is affinoid.

We may now formally imitate the construction of the Katz-Oda spectral sequence
\cite[Theorem~3]{katz-oda} to produce a spectral sequence with
\[
E_2^{pq} = H^p(\Gamma(V, \Omega^{\cdot,\log}_{V/W}) \otimes_{\calO(V)}  
\HH^q(\calE^\dual \otimes \calE' \otimes \Omega^{\cdot,\log}_{X/V}))
\implies 
\HH^{p+q}(X, \calE^\dual \otimes \calE' \otimes \Omega^{q,\log}_{X/W}):
\]
namely, it is the spectral sequence associated to the filtration on the
$\Omega^{\cdot,\log}_{X/W}$ with
\[
\Fil^i(\Omega^{\cdot,\log}_{X/W}) = \image(\Omega^{\cdot - i,\log}_{X/W}
\otimes_{\calO_X} f^*(\Omega^{i,\log}_{V/W}) \to \Omega^{\cdot,\log}_{X/W})
\]
with respect to the derived functors of $\RR^0 \Gamma(X, \cdot)$.

Using the Katz-Oda spectral sequence
(and Lemma~\ref{L:ext1} to translate between Ext groups and
cohomology), we may argue that 
it suffices to prove the desired result in
 the case $V = W$: if each step in the spectral
sequence commutes with the application of $\calU_I$, then so does the 
final result.
Again by passing from $V$ to a suitable cover, we may reduce to the case
where $\calE$ and $\calE'$ are actually free over $\calO$ (and $V=W$);
by arguing again using short exact sequences, we may then
reduce to the case $\calE = \calE' = \calO$.

To summarize, we have so far reduced to consider the case where $\calE = \calE' 
= \calO$ and $V = W$. (Note that since $V = W$, the logarithmic structure
on $V$ no longer intervenes in the calculation.)
Since $V$ is affinoid, $V \times A^n_K(I)$ is a quasi-Stein space by
Corollary~\ref{C:quasi-Stein}, so is acyclic for the cohomology of
coherent sheaves by Kiehl's theorem \cite[Satz~2.4]{kiehl}. Hence
the hypercohomology in
Lemma~\ref{L:ext1} may be computed directly 
on global sections. With this in mind,
we note that the functoriality map
$\Ext^i(\calO, \calO) \to \Ext^i(\calU_I(\calO), \calU_I(\calO))$ 
translates via Lemma~\ref{L:ext1}
into the map on cohomologies induced by the map on complexes
\[
g: \,\Gamma(V \times \{0\}/V, \Omega^{\cdot,\log}_{V \times A^n_K[0,0]/V}) \to
\Gamma(V \times A^n_K(I)/V, \Omega^{\cdot,\log}_{V \times A^n_K(I)/V})
\]
induced by the embedding of $\calO(V \times A^n_K[0,0])$
into $\calO(V \times A^n_K(I))$.

Let $h$ denote the map on complexes
obtained from the ``constant coefficient'' map
$\calO(V \times A^n_K(I)) \to \calO(V)$
(that is, expanding a function as a Laurent series
in $t_1, \dots, t_n$ and extracting the constant coefficient). Then
once we identify $V$ with $V \times \{0\}$,
$h \circ g$ becomes the identity map; we claim that $g \circ h$ is
homotopic to the identity map. One such homotopy can be reconstructed
from the following description on monomials.
Given  the $k$-form
\[
t_1^{i_1} \dots t_n^{i_n}
\frac{dt_{j_1}}{t_{j_1}} \wedge \cdots \wedge \frac{dt_{j_k}}{j_k}
\qquad
(i_1, \dots, i_n \in \ZZ; 1 \leq j_1 < \cdots < j_k \leq n),
\]
pick out the first integer $h$ such that $i_h \neq 0$, and integrate
against $dt_{h}/t_{h}$ (obtaining zero if $h$ is not among
$j_1, \dots, j_k$); this gives a well-defined operation on the complex because
$I$ is quasi-open, so the convergence condition is not disturbed by the
integration.

Since $g$ admits a homotopy inverse, it induces bijections on
cohomology, as desired.
\end{proof}

\begin{remark}
Note that even though only the cases $i=0,1$ of Lemma~\ref{L:ext2} are
needed in what follows,
the higher cases are needed in order to apply the five lemma in the induction
within the proof of Lemma~\ref{L:ext2}.
\end{remark}

\begin{theorem} \label{T:unipotent equiv}
Under Hypothesis~\ref{H:setup conn},
for any nonempty quasi-open subinterval $I$ of $[0,+\infty)$,
the functor $\calU_I: \ULNM_{V \times A^n_K[0,0]/W} 
\to \ULNM_{V \times A^n_K(I)/W}$ is an equivalence of categories.
\end{theorem}
\begin{proof}
By Lemma~\ref{L:ext2} applied in the cases $i=0$ and $i=1$ (as in the proof
of \cite[Proposition~6.7]{crewfin}), $\calU_I$ is an equivalence
whenever $V$ and $W$ are both affinoid. In general,
faithfulness of $\calU_I$ may be checked locally on $V$ and $W$, so
it follows from the affinoid case. Similarly, given faithfulness, 
full faithfulness may be checked locally; given full faithfulness,
essential surjectivity may be checked locally.
\end{proof}
\begin{cor} \label{C:unipotent local}
For $I$ quasi-open and $W = \Maxspec K$, 
the property of an element of $\LNM_{V \times A^n_K(I)/W}$ being
constant/unipotent may be checked locally on $V$.
\end{cor}

\begin{cor} \label{C:horizontal}
Under Hypothesis~\ref{H:setup conn} with $W=V$,
for $\calE \in \ULNM_{V \times A^n_K[0,0]/V}$,
and $I$ an aligned subinterval of $[0,+\infty)$ of positive length,
there is a natural isomorphism
\[
H^0_V(V \times A^n_K[0,0], \calE) \cong H^0_V(V \times
A^n_K(I), \calU_I(\calE)).
\]
\end{cor}
\begin{proof}
Note that elements of $H^0$ can be viewed as homomorphisms from the trivial
log-$\nabla$-module on $V \times A^n_K[0,0]$ (i.e., the sheaf
$\calO_V$ equipped with $n$ endomorphisms all equal to zero).
Hence if $I \subseteq [0,a) \subseteq [0,+\infty)$, 
then by Theorem~\ref{T:unipotent
equiv}, we have a natural isomorphism
\[
H^0_V(V \times A^n_K[0,0], \calE) \cong 
H^0_V(V \times
A^n_K[0,a), \calU_{[0,a)}(\calE)).
\]
inverting the restriction map. We then have another restriction map
\[
H^0_V(V \times
A^n_K[0,a), \calU_{[0,a)}(\calE)) \to 
H^0_V(V \times
A^n_K(I)), \calU_I(\calE));
\]
by Theorem~\ref{T:unipotent equiv}, the composite map
\[
H^0_V(V \times A^n_K[0,0], \calE) \to H^0_V(V \times
A^n_K(I), \calU_I(\calE))
\]
does not depend on the choice of $a$. This composite map is clearly injective;
to see that it is surjective, compose further with the injection
$H^0_V(V \times A^n_K(I), \calU_I(\calE)) \to
H^0_V(V \times A^n_K(J), \calU_J(\calE))$ for any $J \subseteq I$
quasi-open and note that the result is an isomorphism by 
Theorem~\ref{T:unipotent equiv}.
\end{proof}
\begin{remark}
Corollary~\ref{C:horizontal} depends crucially on the nilpotent
residues hypothesis; compare Remark~\ref{R:need nilpotent}.
\end{remark}

As a further consequence of Theorem~\ref{T:unipotent equiv},
we can make an argument that allows us to
ignore hereafter ``logarithmic structure on the base''.
\begin{prop} \label{P:saturation}
Let $X$ be a smooth rigid space, and suppose the zero loci of $t_1,
\dots, t_n \in \Gamma(X, \calO)$ are smooth and meet transversely; let 
$U$ be the complement of these zero loci. Let
$\calE$ be a log-$\nabla$-module with nilpotent 
residues on $X$ with respect to $t_1, \dots, t_n$, and let 
$\calF$ be a $\nabla$-submodule of the restriction of $\calE$ to $U$.
Then $\calF$ extends uniquely to a log-$\nabla$-submodule of $\calE$
with nilpotent residues.
\end{prop}
\begin{proof}
By induction, it suffices to check the following. Let $Z$ be the zero locus
of $t_n$. Suppose that $Z$ is irreducible and that $\calF$
is a log-$\nabla$-submodule with nilpotent residues of the
restriction of $\calE$ to $X \setminus Z$.
Then $\calF$ extends uniquely to a 
log-$\nabla$-submodule of $\calE$ with nilpotent residues.

Since this claim is local (because of the uniqueness assertion), 
we may assume further that $\calE$ and $\calF$ are
free, and that (by imitating the construction of
\cite[Proposition~1.3]{grosse-klonne}, as was done already in
Lemma~\ref{L:small unip}) there exists an admissible subspace
of $X$ containing $Z$ and isomorphic to $Z \times A^1_K[0,a)$ via
a map carrying $Z$ to the zero section. 
Moreover, by Lemma~\ref{L:small unip1}, we may choose $a$ so that
$\calE$ is unipotent on $Z \times A^1_K[0,a)$.
Note that $Z \times A^1_K[0,a)$ and $X \setminus Z$
form an admissible covering of $X$; it thus suffices to exhibit a unique
extension of $\calF$ to $Z \times A^1_K[0,a)$. That is, we may as well
assume outright that $X = Z \times A^1_K[0,a)$ at this point. 

Write $\calE = \calU_{[0,a)}(\calG)$ for some $\calG \in 
\LNM_{Z \times A^1_K[0,0]}$. Then $\calF$ is a subobject 
in $\LNM_{Z \times A^1_K(0,a)}$ of
the restriction of $\calE$ to $\ULNM_{Z \times A^1_K(0,a)}$; by 
Proposition~\ref{P:abelian category}, $\calF$ is itself
unipotent on $Z \times A^1_K(0,a)$. That is, we can write
$\calF = \calU_{(0,a)}(\calH)$ for some $\calH \in \LNM_{Z \times A^1_K[0,0]}$.
By Theorem~\ref{T:unipotent equiv},
the inclusion $\calF \hookrightarrow \left. \calE \right|_{Z \times A^1_K(0,a)}$
is induced by an inclusion $\calH \hookrightarrow \calG$, so
we may take $\calU_{[0,a)}(\calH)$ as the desired extension of $\calF$.
To establish uniqueness of the extension, note that
any two such extensions are both unipotent on some $Z \times A^1_K[0,a)$
by Lemma~\ref{L:small unip1}, so must be isomorphic by
Theorem~\ref{T:unipotent equiv}.
\end{proof}

\subsection{Unipotence and generization}

We now adapt a recipe from \cite[\S 5.3]{me-finite} for iteratively
constructing horizontal elements of a differential module; it shows that
the property of unipotence is ``generic on the base'' in a certain sense.

\begin{lemma} \label{L:gener1}
Let $A$ be an integral
affinoid algebra with $V = \Maxspec A$ smooth over $K$, and let
$L$ be a field containing $A$ which is complete for a norm restricting
to the spectral seminorm on $A$.
(Note that the existence of $L$ forces the reduction of $A$ to be integral.)
Let $I$ be a quasi-open subinterval of $[0,+\infty)$,
take $\calE \in \LNM_{V \times A^n_K(I)/V}$, and let
$\calF$ be the induced element of $\LNM_{A^n_L(I)/\Maxspec L}$.
If $\calF$ is unipotent, then $H^0_V(V \times A^n_K[b,c], \calE) \neq 0$
for any closed aligned subinterval $[b,c]$ of $I$.
\end{lemma}
\begin{proof}
By Theorem~\ref{T:unipotent equiv}, we can express $\calF$ as
$\calU_I(W)$ for some finite dimensional vector space $W$ over $L$
equipped with commuting nilpotent endomorphisms $N_1, \dots, N_n$.
Let $m$ be the minimal length of a unipotent filtration of $W$;
we can then choose $i_1, \dots, i_{m-1} 
\in \{1, \dots, n\}$ such that $N_{i_1} \cdots N_{i_{m-1}} \neq 0$
but $N_{i_1}\cdots N_{i_{m-1}} N_i = 0$ for $i = 1, \dots, n$.

Define the sequence of operators $D_l$ on $\calE$ as follows:
\[
D_l = \prod_{h=1}^{m-1} \left( t_{i_h} \frac{\del}{\del t_{i_h}} \right)
\prod_{i=1}^n \prod_{j=1}^l \left( 1 - 
\frac{t_i}{j} \frac{\del}{\del t_i} \right)^m \left(
1 + \frac{t_i}{j} \frac{\del}{\del t_i} \right)^m.
\]
Pick a closed aligned subinterval $[d,e]$ of $I$ with
$d \leq b$ with strict inequality if $b>0$, and $c < e$.
We claim that for $\bv \in \Gamma(V \times A^n_K[d,e], \calE)$, 
the sequence
$D_l(\bv)$ converges to an element of $H^0_V(V \times A^n_K[b,c], \calE)$.
It suffices to check this in $\Gamma(A^n_L[b,c], \calF)$, where we can 
write $\bv = \sum_J \bv_J t_1^{j_1}\cdots t_n^{j_n}$ for some 
$\bv_J \in W$. In this representation, we have
\begin{equation} \label{eq:dl}
D_l(\bv) = \sum_J t_1^{j_1}\cdots t_n^{j_n} (j_{i_1} + N_{i_1}) 
\cdots (j_{i_{m-1}} + N_{i_{m-1}}) 
\prod_{i=1}^n
\prod_{j=1}^l \left( 1 - \frac{(j_i+N_i)^2}{j^2} \right)^m \bv_J.
\end{equation}
We now analyze the situation following 
\cite[Lemmas~5.3.1 and~5.3.2]{me-finite}.
We may multiply out the summand in \eqref{eq:dl} to get a collection
of terms, each of which consists of $t_1^{j_1} \cdots t_n^{j_n}$ times
a rational number times at most $m-1$ factors from among
$\{N_1, \dots, N_n\}$ (repetitions allowed) times $\bv_J$.
(Remember that the product of any $m$ of the operators $N_1, \dots, N_n$
vanishes, so we can ignore any such product.)
There is a unique term with no $N$'s, in which the rational number factor
is $j_{i_1} \cdots j_{i_{m-1}}$ times the product of the binomial coefficients
$\binom{-j_i-1}{l} \binom{j_i+l}{l}$ for $i=1, \dots, n$; 
in particular, this factor is an integer.
A term with some number $h \leq m-1$ of $N$'s as factors will have a rational
number factor which can be obtained from the integral product we just 
described by multiplying by some integer and then
dividing by $h$ integers, each of absolute value
at most $\max\{|j_1|, \dots, |j_n|\} + l$.

This means that the norms of the terms of the $t$-adic expansion of
$D_l(\bv) - N_{i_1} \cdots N_{i_{m-1}} \bv_0$
are dominated by the norms of the terms of the sum
dominated by
\[
\sum_{J: |j_1|, \dots, |j_n| > l} \sum_N
(\max\{|j_1|, \dots, |j_n|\}+l)^{-m+1} t_1^{j_1} \cdots t_n^{j_n} N \bv_J,
\]
where $N$ runs over the number of products of at most $m-1$
of the operators $N_1, \dots, N_n$ with repetitions allowed. For each
fixed $N$, if we were to consider the sequence (as $l$ varies) of summands
with the factor $(\max\{|j_1|, \dots, |j_n|\}+l)^{-m+1}$ removed,
then Lemma~\ref{L:truncation} would force the sequence to be
$\eta$-null over $X \times A^n_K[b,c]$ for some $\eta > 1$.
Putting the factor back in, we obtain the same conclusion by replacing $\eta$
by any smaller value, since $|\max\{|j_1|, \dots, |j_n|\}+l|^{-m+1}$ 
is dominated by $\rho^l$ for any $\rho>1$.

We conclude that $\{D_l(\bv) - N_{i_1} \cdots N_{i_{m-1}} \bv_0\}_{l=0}^\infty$ 
is $\eta$-null over $V \times A^n_K[b,c]$ for some $\eta > 1$, and
so in particular is convergent to zero. Hence
the $D_l(\bv)$ converge to an element of $H^0_{V}(V \times A^n_K[b,c], \calE)$,
and the limit is nonzero if and only if $N_{i_1} 
\cdots N_{i_{m-1}} \bv_0 \neq 0$
(since the map $\Gamma(V \times A^n_K[b,c], \calE) \to
\Gamma(A^n_L[b,c], \calF)$ is injective).
Let $\bv_1, \dots, \bv_k$ be a set of
generators of $\Gamma(V \times A^n_K[d,e], \calE)$.
If $0 \in I$,
take $S = \{\bv_1, \dots, \bv_k\}$; otherwise, let
$S$ be the set consisting of elements of 
$\Gamma(V \times A^n_K[d,e], \calE)$ of the form $t^J \bv_l$ for 
$J \in \ZZ^n$ and $l \in \{1, \dots, k\}$. Then as $\bv$ runs over $S$,
the resulting values of $\bv_0$ must span $W$ over $L$; in particular, we can choose
$\bv$ so that $N_{i_1} \cdots N_{i_{m-1}} \bv_0 \neq 0$, 
and so the limit of the
$D_l(\bv)$ is a nonzero element of $H^0_V(V \times A^n_K[b,c], \calE)$.
\end{proof}

\begin{prop} \label{P:generization}
Let $A$ be an integral
affinoid algebra with $V = \Maxspec A$ smooth over $K$, 
take $x_1, \dots, x_m \in A$ whose zero loci are smooth 
and meet transversely, and let
$L$ be a field containing $A$ which is complete for a norm restricting
to the spectral seminorm on $A$.
Let $I$ be a quasi-open subinterval of $[0,1)$,
take $\calE \in \LNM_{V \times A^n_K(I)}$, and let
$\calF$ be the induced element of $\LNM_{A^n_L(I)}$.
Then $\calE$ is constant (resp.\ unipotent) if and only if $\calF$
is constant (resp.\ unipotent).
\end{prop}
\begin{proof}
If $\calE$ is constant (resp.\ unipotent), then clearly $\calF$
is constant (resp.\ unipotent). We prove the converse by induction on 
the rank of $\calE$.

Let $[b,c]$ be any closed aligned subinterval of $I$ of positive length.
By Lemma~\ref{L:gener1}, $H^0_V(V \times A^n_K[b,c], \calE)$ is nonzero;
if we let $V'$ be the complement on $V$ of the zero loci of $x_1, \dots, x_m$,
then it follows that $H^0_{V'}(V' \times A^n_K[b,c], \calE)$ is also
nonzero.
By Proposition~\ref{P:saturation}, the $\calO_{V' \times A^n_K[b,c]}$-span of
$H^0_{V'}(V' \times A^n_K[b,c], \calE)$ extends to a subobject $\calG$
of $\calE$ in
$\LNM_{V \times A^n_K[b,c]}$; we will show that $\calG$ is constant.

Let $\calH \in \LNM_{A^n_L[b,c]}$ be induced by $\calG$.
Since $\calF$ is unipotent and $\calH$ injects into $\calF$,
$\calH$ is unipotent by the proof of Proposition~\ref{P:abelian category}.
On the other hand, $\calH$ is also generated by global sections, namely those
coming from $H^0_{V'}(V' \times A^n_K[b,c], \calE)$, so $\calH$ is constant.
In particular, $H^0(A^n_L[b,c], \calH)$ is a finite-dimensional $L$-vector
space and, writing $\pi_L$ for the structure map $A^n_L[b,c] \to \Maxspec
L$, the natural map $\pi_L^* H^0(A^n_L[b,c],
\calH) \to \calH$ is an isomorphism.

For any finitely generated $\calO_V$-submodule $M$ of
$H^0_V(V \times A^n_K[b,c], \calG)$, we have a commuting diagram
\[
\xymatrix{ 
\pi_1^* M \ar[r] \ar[d] & \calG \ar[d] \\
\pi_L^* H^0(A^n_L[b,c], \calH) \ar[r] & \calH,
}
\]
in which the vertical arrows are visibly injective. We showed above that
the lower horizontal arrow is an isomorphism, so the upper horizontal arrow
is also injective. 
Since $\calG$ is a finitely generated module over
the noetherian sheaf of rings $\calO_{V \times A^n_K[b,c]}$, we can 
choose some $M$ as above, which we call $M_1$, such that
$\pi_1^* M_1$ is maximal among the $\pi_1^* M$.
On the other hand, 
if $M_2$ were a $\calO_V$-submodule of $H^0_V(V \times A^n_K[b,c],
\calG)$ strictly containing $M_1$, then $\pi_1^* M_2$ would strictly contain 
$\pi_1^* M_1$. We conclude that $H^0_V(V \times A^n_K[b,c], \calG) = M_1$
is finitely generated over $\calO_V$, and that
$\pi_1^* H^0_V(V \times A^n_K[b, c], \calG) = \pi_1^* M_1
\to \calG$ is injective.

We next prove that the map 
$\pi_1^* H^0_V(V \times A^n_K[b,c], \calG) \to \calG$ is surjective.
With notation as in the proof of Lemma~\ref{L:gener1},
let $f(\bv)$ denote the limit of the $D_l(\bv)$.
Then for any $\bv \in \Gamma(V \times A^n_K[d,e], \calG)$,
we have $f(\bv) \in H^0_V(V \times A^n_K[b,c], \calG)$.

Suppose that $b \neq 0$. Then
\[
\bv = \sum_{J \in \ZZ^n} t_1^{j_1} \cdots t_n^{j_n} 
f(t_1^{-j_1} \cdots t_n^{-j_n} \bv)
\]
as an equality of sections of $\calG$ on $V \times A^n_K[b,c]$;
this implies that $\bv \in \pi_1^* H^0_V(V \times A^n_K[b,c], \calG)$.
Since $\Gamma(V \times A^n_K[d,e], \calG)$ is dense in
$\Gamma(V \times A^n_K[b,c], \calG)$, this yields the desired surjectivity.

Suppose now that $b = 0$.
Before proceeding further, 
we verify that for $\alpha \in \Gamma^*$ and $J \in \ZZ^n_{\geq 0}$,
an element $\bx \in \Gamma(V \times A^n_K[0,\alpha], \calG)$
is divisible by a monomial $t^J$ if and only if the restriction
$\bx_L$ of $\bx$ to $\Gamma(A^n_L[0,\alpha], \calH)$ is divisible by $t^J$;
that is, we can check divisibility by $t^J$ from the expansion of $\bx$
as a formal series in $t_1,\dots,t_n$.
By induction on the sum of the entries of $J$,
it suffices to check the claim for $t^J = t_1$. Write $\iota, \iota_L$
for the inclusions
\[
V \times A^{n-1}_K[0,\alpha] \to V \times A^n_K[0,\alpha],
\qquad
A^{n-1}_L[0,\alpha] \to A^n_L[0,\alpha]
\]
into the locus $t_1 = 0$.
Write $\iota^*, \iota_L^*$ for the induced morphisms
\[
\Gamma(V \times A^n_K[0,\alpha], \calG) \to
\Gamma(V \times A^{n-1}_K[0,\alpha], \calG), \qquad
\Gamma(A^n_L[0,\alpha], \calH) \to \Gamma(A^{n-1}_L[0,\alpha], \calH).
\]
Then $\bx$ is divisible by $t_1$ if and only if $\iota^*(\bx) = 0$,
which is equivalent to $\iota^*_L(\bx_L) = 0$ because the restrictions
$\Gamma(V \times A^*_K[0,\alpha], \calG) \to \Gamma(A^*_L[0,\alpha], \calH)$
are injective for $* \in \{n-1,n\}$.
The latter is equivalent to $\bx_L$ being divisible by $t_1$, as desired.

For $J, J'$, write $J \leq J'$ if $J'$ is componentwise greater than or
equal to $J$.
Let $J_0, J_1, \dots$ be a total ordering of 
$\ZZ_{\geq 0}^n$ refining the partial ordering $\leq$;
write $J_j = (a_{j,1}, \dots, a_{j,n})$.
Choose a decreasing sequence of aligned intervals
\[
[d,e] = [0, e_0] \supset [0, e_1] \supset \cdots
\]
satisfying $\cap_j [0, e_j] \supseteq [0, c] = [b,c]$.
For $\alpha$ equal to one of the $e_j$, and $k \in \{1, \dots, n\}$,
write
\[
\hat{\pi}_k: V \times A^n_K[0,\alpha] \to V \times A^{n-1}_K[0,\alpha]
\]
for the projection omitting the $k$-th coordinate of $A^n_K[0,\alpha]$.
Let $f_k(\bv)$ denote the limit of the $D_l(\bv)$ when computed for the
projection $\hat{\pi}_k$; then $f_k$ defines a map 
$\Gamma(V \times A^n_K[0,e_j], \calG) \to
\Gamma(V \times A^n_K[0, e_{j+1}], \calG)$.
For $\bv \in \Gamma(V \times A^n_K[0,\alpha], \calG)$, write 
$\bv = \sum_J t^J \bv_J$ for the series expansion of $\bv$
over $A^n_L[0,\alpha]$; in terms of these series, $f_k$ acts as
\[
\sum_{j \geq 0} t^{J_j} \bv_j \mapsto \sum_{j \geq 0, a_{j,k} = 0} 
t^{J_j} \bv_j.
\]

We check by induction on $j$ that $\sum_{J_j \leq J} t^J
\bv_J \in \Gamma(V \times A^n_K[0,e_j], \calG)$. This is given for
$j=0$; if $j > 0$, we can choose $k \in \{1,\dots,n\}$ such that
$a_{j,k} > 0$, and there is an index $j' < j$ such that
\[
J_{j'} = (a_{j,1}, \dots, a_{j,k} - 1, \dots, a_{j,n}).
\]
By the induction hypothesis, 
\[
\sum_{J_{j'} \leq J} t^J \bv_J
\in \Gamma(V \times A^n_K[0,e_{j'}], \calG)
\subseteq \Gamma(V \times A^n_K[0,e_{j-1}], \calG)
\]
and $\sum_{J_{j'} \leq J} t^J \bv_J$ is divisible by 
$t_k^{a_{j,k}-1}$, since we showed above that we
can check this divisibility  on the level of formal series.
We now have
\[
\sum_{J_j \leq J} t^J \bv_J =
\sum_{J_{j'} \leq J} t^J \bv_J - 
t_k^{a_{j,k}-1} f_k (t_k^{-a_{j,k}+1}
\sum_{J_{j'} \leq J} t^J \bv_J)
\in \Gamma(V \times A^n_K[0,e_j], \calG),
\]
completing the induction.

By the previous paragraph, for each $j$,
$\sum_{J_j \leq J} t^J \bv_J \in \Gamma(V \times A^n_K[0,e_j],
\calG)$, and 
$\sum_{J_j \leq J} t^J \bv_J$ is divisible by $t^{J_j}$
because again we can check this divisibility on the level
of formal series. We then have
\[
\bv_{J_j} = f(t^{-J_j} \sum_{J_j \leq J} t^J \bv_J)
\in H^0_V(V \times A^n_K[0,c],\calG).
\]
Because the norm on $L$ is compatible with that on $V$, the sum
$\sum_J t^J \bv_J$ converges to $\bv$, and so
$\bv \in \pi_1^* H^0_V(V
\times A^n_K[0,c], \calG)$ as in the case $b \neq 0$.
Again because the restriction map $H^0_V(V \times A^n_K[0,e],\calG)
\to H^0_V(V \times A^n_K[0,c], \calG)$ has dense image,
this yields the desired surjectivity.

In either of the cases $b\neq 0$ or $b=0$, we now see that the map
$\pi_1^* H^0_V(V \times A^n_K[b,c], \calG) \to \calG$ is surjective;
since we already showed injectivity, the map is an isomorphism.
At this point, there is no harm in replacing $K$ by a finite extension,
as what we are checking is local freeness and nilpotence of residues
for $H^0_V(V \times A^n_K[b,c], \calE)$. In particular,
we may assume that $A^n_K[b,c]$ contains a $K$-rational point $x$.

Writing $i$ for the injection
$V \times \{x\} \to V \times A^n_K[b,c]$, we obtain an isomorphism
\[
H^0_V(V \times A^n_K[b,c], \calG) = i^* \pi_1^* H^0_V(V 
\times A^n_K[b,c], \calG) \cong i^* \calG.
\]
Conesquently, $H^0_V(V \times A^n_K[b,c], \calG)$ defines an
object in $\LNM_V$, and $\calG \cong \pi_1^* H^0_V(V \times A^n_K[b,c],
\calG)$ is constant over $V$. 

By the induction hypothesis, we may deduce that the restriction of
$\calE$ to $V \times A^n_K(b,c)$ in case $b > 0$,
or $V \times A^n_K[0,c)$ in case $b=0$,
is unipotent over $V$.
Since $[b,c]$ was an arbitrary
closed aligned subinterval of $I$, we deduce by 
Theorem~\ref{T:unipotent equiv} that
$\calE$ is unipotent over all of $V \times A^n_K(I)$,
as desired. (If $\calF$ is constant, then $\calE$ is constant by comparison
of residues.)
\end{proof}
\begin{remark}
Proposition~\ref{P:generization} even makes a nontrivial assertion when
$V = \Maxspec K$, as we may take $L$ to be any
extension of $K$ complete under some extension of $|\cdot|$. The assertion
is that unipotence can be tested after making an arbitrary base field
extension. (As everything involved is $K$-linear, this should not
be surprising; a special case of this was already proved in 
\cite[Proposition~6.11]{me-local} using this linearity.)
\end{remark}

\begin{cor} \label{C:generization}
Let $P$ be a smooth affine formal scheme of finite type over $\gotho_K$,
suppose $x_1,\dots,x_m \in \Gamma(P, \calO)$ have zero loci on $P_K$
which are smooth and meet transversely, and
let $X$ be an open dense subscheme of $P_k$.
Given a quasi-open subinterval $I$ of $[0,1)$ and
an object $\calE \in \LNM_{P_K \times A^n_K(I)}$,
suppose that the restriction of $\calE$ to
$]X[ \times A^n_K(I)$ is constant/unipotent.
Then $\calE$ is constant/unipotent.
\end{cor}
\begin{proof}
By shrinking $X$ further, we may reduce to the case where
$X = P_k \setminus V(g)$ for some $g \in \Gamma(P, \calO)$.
Then $]X[$ is the affinoid space associated to the affinoid algebra
$\Gamma(P_K, \calO) \langle g^{-1} \rangle$; in particular,
$\Gamma(P_K, \calO)$ and $\Gamma(]X[, \calO)$ have the same completed
fraction field $L$. We may thus apply
Proposition~\ref{P:generization} to deduce that $\calE$
induces a constant/unipotent $\nabla$-module over $A^n_L(I)$,
and then that $\calE$ is constant/unipotent
over $P_K \times A^n_K(I)$.
\end{proof}

\subsection{Unipotence and overconvergent generization}

We will also need a variant of the
construction of Proposition~\ref{P:generization}
in which we allow an ``overconvergent'' base. We start with a Gr\"obner
basis calculation derived from \cite[Section~2.4]{me-finite}, but modified
to avoid relying on discreteness of $K$.

\begin{lemma} \label{L:grobner basis}
For $\lambda \in [1, \infty) \cap \Gamma^*$, 
let $R_\lambda$ be the (affinoid) ring of rigid analytic functions on the subspace
\[
|x_1|\leq 1, \quad \dots, \quad |x_{n-1}| \leq 1, \quad |x_n| \leq \lambda
\]
of the rigid affine $n$-space over $K$, and write $|\cdot|_\lambda$
for the supremum norm on $R_\lambda$.
Let $\gotha$ be an ideal of $R_\delta$ for some $\delta \in
(1, \infty) \cap \Gamma^*$. Then there exists $\rho_0 \in (1, \delta]
\cap \Gamma^*$ such that for any $\rho \in (1, \rho_0] \cap \Gamma^*$
and any $y,z \in R_{\delta}$ with $y-z \in \gotha$, one can find 
$u \in R_\delta$ with 
\[
u-z \in \gotha, \quad |u|_1 \leq |y|_1, \quad |u|_\rho \leq |z|_\rho.
\]
\end{lemma}
\begin{proof}
If $y = 0$, we may take $u = 0$, so we assume instead that $y \neq 0$.
Choose a total ordering $\leq$ on $\ZZ^n_{\geq 0}$ extending the partial order
$\preceq$
by termwise comparison and the partial order by comparison only in the last 
component. The former partial order is well-founded, so the total ordering
is a well ordering.

For $y = \sum y_I x^I \in R_\delta$ and $\lambda \in [1,\delta] \cap
 \Gamma^*$, define
the \emph{$\lambda$-leading term} of $y$ to be the expression $y_I x^I$
for $I$ the largest tuple under $\leq$ which maximizes
$|y_I x^I|_\lambda = |y_I| \lambda^{i_n}$; 
such a tuple exists because there only finitely many
tuples achieving the maximum.

We claim that for each $y \in R_\delta$, 
the $1$-leading term of $y$ coincides with the
$\rho$-leading term for each sufficiently small $\rho \in (1,\delta]
\cap \Gamma^*$ (depending on $y$). To see this, let $y_I x^I$ be the 1-leading
 term of $y$.
For each tuple $J$, we then have either 
\begin{enumerate}
\item[(a)] $|y_J| < |y_I|$, or 
\item[(b)]
$|y_J| = |y_I|$ and $J \leq I$; in this case we have $j_n \leq i_n$.
\end{enumerate}
If $|y_J x^J|_\delta \leq |y_I x^I|_\delta$, then in case (a),
we have $|y_J x^J|_\rho < |y_I x^I|_\rho$ for all
$\rho \in [1, \delta) \cap \Gamma^*$;
in case (b), we have $|y_J x^J|_\rho \leq |y_I x^I|_\rho$ and $J \leq I$.
So these terms are all okay for any $\rho$; in fact, because $y \in R_\delta$,
there are only finitely many tuples $J$ with
$|y_J x^J|_\delta > |y_I x^I|_\delta$.
For each such $J$, we must be in case (a), so $|y_J x^J|_\rho 
< |y_I x^I|_\rho$ for $\rho \in (1,\delta]$
sufficiently small. This yields the claim.

Define elements $a_1, a_2, \dots$ of $\gotha$ as follows.
Given $a_1, \dots, a_{i-1}$, choose $a_i$ if possible to be an element
of $\gotha$ whose $1$-leading term is not a multiple of the
$1$-leading term of $a_j$ for any $j<i$,
otherwise stop. By the well-foundedness of
$\preceq$, this process must eventually stop; at that point, every 
$1$-leading term of every element of $\gotha$ is a multiple of 
the $1$-leading term of some $a_i$. Let $A$ be the finite set consisting
of the $a_i$ just constructed.

As shown above, we can choose $\rho_0 \in (1, \delta] \cap \Gamma^*$
such that for $\rho \in [1, \rho_0) \cap \Gamma^*$, the $1$-leading
term and $\rho$-leading term of each $a \in A$ coincide. Moreover, we can
choose $\epsilon \in (0,1)$ such that for each $a \in A$, if $y_I x^I$ is the
1-leading term of $a$, then for each $J$ in case (a) above, we
actually have $|y_J| \leq \epsilon |y_I|$. (Namely, for any particular
$\epsilon$, there are only finitely $J$ contradicting this inequality;
by making $\epsilon$ large enough, we can eliminate all of these.)

We construct a sequence $\{c_j\}$ of monomials and a sequence
$\{d_j\}$ of elements of $A$ as follows. Given the sequences up to
$c_j$ and $d_j$, put $z_j = z - c_1 d_1 - \cdots - c_j d_j$ (or $z_0 = z$ initially). If
$|z_j|_1 \leq |y|_1$, then stop. Otherwise, let $e_I x^I$ be the 1-leading
term of $z_j - y$. By the construction of $A$, we can find
a monomial $c_{j+1}$ and some $d_{j+1} \in A$ such that $c_{j+1} d_{j+1}$
has 1-leading term, and hence $\rho$-leading term, equal to $e_I x^I$.

{}From the construction, we clearly have $|z_j|_\rho \leq |z|_\rho$.
On the other hand, if the process were never to terminate, we could show
that $|z_j|_1 \to 0$ as $j \to \infty$ as follows. It would suffice to
show that eventually $|z_j|_1 \leq \epsilon |z|_1$, as this argument could
then be iterated. Let $s_j$ be the set of monomials
of $z_j$ of 1-norm greater than $\epsilon |z|_1$. If $s_j$ is nonempty,
then $s_{j+1}$
is obtained from $s_j$ by taking out a term of maximal 1-norm and possibly
adding back in some other terms of the same 1-norm which are smaller
under $\leq$. In particular, the set of all possible 1-norms of elements
of the $s_j$ is finite; moreover, 
since $\leq$ is a well-ordering, we must eventually
run out of terms of any particular 1-norm. Hence eventually $s_j$ becomes 
empty, and so
$|z_j|_1 \leq \epsilon |z|_1$.

Again assuming that the process does not terminate, the previous paragraph
would imply that
$|z_j|_1 \to 0$ as $j \to \infty$. But since
we stop whenever $|z_j|_1 \leq |y|_1$, this can only happen if 
$y = 0$, which contradicts an earlier assumption.
Thus the process terminates at some $z_j$, and we
may take $u=z_j$.
\end{proof}

\begin{prop} \label{P:gen Hadamard}
Let $X$ be a reduced affinoid space, and take $f \in \calO(X)$ with
$|f|_X = \delta > 1$. For $\lambda \in [1,\delta] \cap \Gamma^*$, put
$U_\lambda = \{x \in X: |f(x)| \leq \lambda\}$. Suppose that $U_1 \neq \emptyset$.
Then for each 
$c \in (0,1) \cap \QQ$, there exists $\lambda \in (1,\delta] \cap \Gamma^*$ such that
for all $g \in \calO(X)$,
\[
|g|_{U_\lambda} \leq |g|_{U_1}^c |g|_{X}^{1-c}. 
\]
\end{prop}
\begin{proof}
With notation as in Lemma~\ref{L:grobner basis}, we can choose a
closed immersion $\phi: X \hookrightarrow \Maxspec R_\delta$ which pulls
$x_n$ back to $f$; then $U_\lambda = \phi^{-1}(\Maxspec R_\lambda)$.
We then choose $\rho_0$ as in Lemma~\ref{L:grobner basis}.

For each $\lambda \in [1,\delta] \cap \Gamma^*$,
the supremum norm on $U_\lambda$ is equivalent to the quotient norm
induced from $R_\lambda$. We can thus choose $\epsilon > 1$ such that
for any $g \in \calO(X)$, there exist
$y, z \in R_\delta$ with
\[
\phi^*(y) = \phi^*(z) = g,
\quad
|y|_1 \leq \epsilon |g|_{U_1},
\quad
|z|_{\rho_0} \leq \epsilon |g|_{U_{\rho_0}}.
\]
By Lemma~\ref{L:grobner basis}, we can choose $u \in R_\delta$
with
\[
\phi^*(u) = g, \quad |u|_1 \leq |y|_1, \quad |u|_{\rho_0} \leq |z|_{\rho_0}.
\]
Now put $\lambda = \rho_0^{1-c}$; by Lemma~\ref{L:Hadamard}, we have
\begin{align*}
|g|_{U_\lambda} &\leq |u|_\lambda \\
&\leq |u|_1^c |u|_{\rho_0}^{1-c} \\
&\leq |y|_1^c |z|_{\rho_0}^{1-c} \\
&\leq \epsilon |g|_{U_1}^c |g|_{U_{\rho_0}}^{1-c} \\
&\leq \epsilon |g|_{U_1}^c |g|_{X}^{1-c}.
\end{align*}
Since supremum norms are multiplicative, applying the same argument to
$g^n$ instead of $g$ yields
\[
|g|_{U_\lambda} \leq \epsilon^{1/n} |g|_{U_1}^c |g|_{X}^{1-c},
\]
and the desired result now follows by taking the limit as $n \to \infty$.
\end{proof}

\begin{prop} \label{P:generization2}
Let $P$ be an affine formal scheme of finite type over $\gotho_K$, and let
$X$ be an open dense subscheme of $P_k$ such that $P$ is smooth in a 
neighborhood of $X$.
Take $x_1, \dots, x_m \in \Gamma(P, \calO)$ whose
zero loci on $P_K$ are smooth and meet transversely.
Let $I$ be a quasi-open subinterval of $[0,1)$,
let $V$ be a strict neighborhood of $]X[$ in $P_K$, and suppose
that $\calE \in \LNM_{V \times A^n_K(I)}$ becomes
constant/unipotent on $]X[ \times A^n_K(I)$. 
Then for any closed aligned subinterval $[b,c] \subset I$ of positive length, 
there exists a strict neighborhood
$V'$ of $]X[$ in $P_K$ such that 
$\calE$ is constant/unipotent over $V' \times A^n_K[b,c]$.
\end{prop}
\begin{proof}
We may assume without loss of generality that $V$ is affinoid.
Let $[d,e] \subset I$ 
be a closed aligned subinterval
with $[b,c] \subseteq [d,e)$, and with $d < b$ unless
$b = 0$.
As in the proof of Lemma~\ref{L:gener1},
we can choose $\bv \in \Gamma(V \times A^n_K[d,e], \calE)$ such that
the sequence $\{ D_l(\bv) \}$ converges to a nonzero element of
$H^0_{]X[}(]X[ \times A^n_K[b,c], \calE)$.
Moreover, from the construction in 
Lemma~\ref{L:gener1},
we see that there exists $\eta>1$ so that the sequence
$\{D_{l+1}(\bv)- D_l(\bv)\}$ is $\eta$-null over $]X[ \times A^n_K[b,c]$.

Suppose $W$ is a connected affinoid subdomain of $V \times A^n_K[d,e]$
over which $\calE$ becomes free. 
Choose a basis $\be_1, \dots, \be_r$ of $\Gamma(W, \calE)$,
and for $i=1, \dots, r$,
let $A_i$ be the matrix via which $t_i \frac{\del}{\del t_i}$
acts on the basis $\be_1, \dots, \be_r$.
Define a system $V_\lambda$ of strict neighborhoods of $]X[$ in
$P_K$ as in Lemma~\ref{L:strict neighborhood}, and let
$g_i(\lambda)$ denote the
maximum supremum seminorm of any entry of $A_i$ over
$W \cap (V_{\lambda} \times A^n_K[d,e])$.
Then we see directly from the definition of $D_l$ that the sequence
$\{D_{l+1}(\bv) - D_l(\bv)\}$ is $\rho$-null over
$W \cap (V_{\lambda} \times A^n_K[d,e])$ for some $\rho > 0$, e.g.,
\[
\rho = (\max\{1,g_1(\lambda)\} \cdots \max\{1,g_n(\lambda)\} |p|^{-1/(p-1)})^{-2m}.
\]

If $W$ has nonempty intersection with $]X[ \times A^n_K[d,e]$,
we may apply Proposition~\ref{P:gen Hadamard} to deduce that
the sequence $\{D_{l+1}(\bv) - D_l(\bv)\}$ is 1-null over 
$W \cap (V_{\lambda} \times A^n_K[b,c])$ for some
$\lambda \in (0,1) \cap \Gamma^*$.
If on the other hand $W$ has empty intersection with $]X[ \times A^n_K[d,e]$,
then by the maximum modulus principle, $W$ also has empty intersection
with $V_\lambda \times A^n_K[d,e]$ for some $\lambda \in (0,1) \cap \Gamma^*$,
so there is nothing to check in this case.

Note that we can cover $V
\times A^n_K[d,e]$ with finitely many affinoid subdomains $W$,
over each of which $\calE$ becomes free.
Hence we can choose $\lambda
 \in (0, 1) \cap \Gamma^*$
such that the limit of the $D_l(\bv)$ exists over $V_{\lambda}
\times A^n_K[b,c]$.
Thus $H^0_{V_{\lambda}}(V_{\lambda} \times A^n_K[b,c], \calE) \neq 0$ for 
some $\lambda$. As in the proof of Proposition~\ref{P:generization},
we may obtain a nonzero constant 
log-$\nabla$-submodule of $\calE$, quotient by it, and repeat to
obtain the desired result.
(The role of $L$ in the proof of Proposition~\ref{P:generization}
is played by a complete field containing $\calO(]X[)$ whose norm
is compatible with the norm on $\calO(V_\lambda)$.)
\end{proof}

\subsection{Convergence and unipotence}
\label{subsec:conv unip}

Contrary to what one's intuition from real analysis would suggest,
a log-$\nabla$-module over $V \times A^n_K[0,1)$ with nilpotent
residues need not be unipotent; see Remark~\ref{R:not unipotent}
below. 
What distinguishes unipotent log-$\nabla$-modules is $\eta$-convergence
(see Definition~\ref{D:eta convergent}), in the following fashion.

\begin{lemma} \label{L:easy convergent}
For any smooth affinoid space $X$, any $a,b \in (0,1) \cap \Gamma^*$ with
$a \leq b$, and any $\calE \in \ULNM_{X \times A^n_K[a,b]/X}$, $\calE$
is $\eta$-convergent with respect to $t_1, \dots, t_n$ (relative to $X$)
for any $\eta < a$. Moreover, if 
$\calE \in \ULNM_{X \times A^n_K[a,b]}$ and 
there exists a point $x \in A^n_K[a,b]$ such that the restriction of
$\calE$ to $X \times \{x\}$ is $\eta$-convergent with respect to some
coordinate system $z_1, \dots, z_l$ on $X$ and some $\eta < a$, then 
$\calE$ is $\eta$-convergent with respect to $t_1, \dots, t_n, z_1, \dots, z_l$.
\end{lemma}
\begin{proof}
First note that the question is local on $X$, so we may reduce to the case
where $\calE$ admits a filtration whose successive quotients are constant
and pulled back from free $\calO_X$-modules.
By Remark~\ref{R:eta conv exact}, we may assume that $\calE$ itself is constant.

Note that the claim in the first instance
holds for $\calE = \calO$ by direct calculation:
for any $x \in \calO(X \times A^n_K[a,b]$), any tuple $R  = (r_1, \dots, r_n)
\in [a,b]^n$, and any tuple $I = (i_1, \dots, i_n)$ of nonnegative integers,
one has
\[
\left| \frac{1}{I!} \frac{\del^{i_1}}{\del t_1^{i_1}} \cdots
\frac{\del^{i_n}}{\del t_n^{i_n}} x\right|_R 
\leq r_1^{-i_1} \cdots r_n^{-i_n} |x|_R,
\]
yielding the $\eta$-convergence.
In particular, $t_1, \dots, t_n$ form an $\eta$-admissible coordinate
system on $X \times A^n_K[a,b]$ relative to $X$.

In the second instance, $\calE$ is obtained by pullback from a
log-$\nabla$-module $\calF$ on $X$, which by the given hypothesis is
$\eta$-convergent with respect to $z_1, \dots, z_l$. The
$\eta$-convergence of $\calE$ follows by the same calculation as in
the previous paragraph.
\end{proof}

\begin{lemma} \label{L:convergent log}
Let $X$ be a smooth affinoid space,
and take $\calE \in \LNM_{X \times A^n_K[0,b]}$ for some $b \in (0,1)
\cap \Gamma^*$. Suppose that the restriction of $\calE$ to 
$X \times A^n_K[a,b]$ is $\eta$-convergent with respect to $t_1, \dots, t_n$
(relative to $X$)
for some $a \in (0,b)
\cap \Gamma^*$ and some $\eta \in (0,a) \cap \Gamma^*$. Then 
$\calE$ is unipotent on $X \times A^n_K[0,\eta)$.
Moreover, if all of the residues are zero, then
$\calE$ is constant on $X \times A^n_K[0,\eta)$.
\end{lemma}
\begin{proof}
We proceed by induction on $n$.
Write $Y = (X \times A^{n-1}_K[0,b]) \times A^1_K[0,\eta]$.
Suppose $\calF \in \ULNM_{Y}$ 
is a (possibly zero) proper
subobject of the restriction of $\calE$.
Let $d$ be the length of the shortest
unipotent filtration of the restriction of the residue of $\calE/\calF$
along $t_n = 0$.
Let $P_j(x)$ denote the $j$-th binomial polynomial, i.e.,
\[
P_j(x) = \frac{x(x-1)\cdots(x-j+1)}{j!} \qquad (j = 1,2,\dots).
\]
Then an exercise in elementary
number theory shows that the $\ZZ$-module of 
polynomials with rational coefficients
carrying $\ZZ$ into itself is freely generated by the $P_n(x)$.
Moreover, if $Q$ is a polynomial carrying $\ZZ$ into itself and
$Q(0) = \cdots = Q(j-1) = 0$, then $Q$ is an integer linear
combination of $P_j, P_{j+1}, \dots, P_{\deg Q}$. (Evaluating at 0 shows that
the coefficient of $P_0$ vanishes; then evaluating at 1 shows that
the coefficient of $P_1$ vanishes, and so on.)
In particular, if we set
\[
Q_j(x) = x^{d-1} \left( \frac{(1-x)\cdots(j-x)}{j!} \right)^d,
\]
then $Q_{j+1}(x) - Q_j(x)$
is an integer linear combination of $P_{j+1}(x), \dots, P_{dj+d-1}(x)$.

By computing on formal power series in $t_n$
(with which we can formally construct a basis of sections killed by
$(t_n \frac{\del}{\del t_n})^d$) or invoking Lemma~\ref{L:small unip1}, 
we see that
\[
(Q_{j+1}- Q_j)\left( t_n \frac{\del}{\del t_n} \right)
\]
carries any element of $\Gamma(Y, \calE/\calF)$
to a multiple of $t_n^{j+1}$ in the same module. That is,
\[
\frac{1}{t_n^{j+1}}(Q_{j+1} - Q_j) 
\left( t_n \frac{\del}{\del t_n} \right)
\]
is a well-defined operator on $\calE/\calF$.
As we saw above, $\frac{1}{t_n^{j+1}} (Q_{j+1}-Q_j)(t_n \frac{\del}{\del t_n})$
is a $\ZZ$-linear combination of
\[
\frac{1}{t_n^{j+1}} P_l\left(t_n \frac{\del}{\del t_n}\right) \qquad (l = j+1, \dots, dj+d-1),
\]
and hence is a $\Gamma(Y, \gotho)$-linear 
combination of the $\frac{1}{t_n^l} P_l(t_n \frac{\del}{\del t_n})$
for $l=j+1, \dots, dj+d-1$.

However,
\[
\frac{1}{t_n^l}
P_l\left(t_n \frac{\del}{\del t_n}\right) = 
\frac{1}{l!} \frac{\del^l}{\del t_n^l};
\]
by the $\eta$-convergence condition,
for any $\bw \in \Gamma(A^n_K[a,b], \calE)$,
the sequence $\{\frac{1}{j!} \frac{\del^j}{\del t_n^j} 
\bw\}_{j=1}^\infty$ is $\eta$-null on $A^n_K[a,b]$.
If we choose 
$\bw \in \Gamma(A^n_K[0,b], \calE)$, then
Lemma~\ref{L:Hadamard} implies that 
$\{\frac{1}{j!} \frac{\del^j}{\del t_n^j} 
\bw\}_{j=1}^\infty$ is also $\eta$-null on $A^n_K[0,b]$,
so in particular on $Y$. In particular, if $\bv$
denotes the image of $\bw$ in $\Gamma(Y, \calE/\calF)$, then
the sequence
\begin{equation} \label{eq:unipotent sequence}
\{t_n^{-j-1} (Q_{j+1} - Q_j)\left(t_n 
\frac{\del}{\del t_n}\right)\bv\}_{j=1}^\infty
\end{equation}
is $\eta$-null on $Y$;
that means that the sequence
$\{(Q_{j+1} - Q_j)(t_n \frac{\del}{\del t_n}) \bv\}$ is 
$1$-null on $Y$.
That is, the limit
\[
f(\bv) = \lim_{j \to \infty} Q_j\left(t_n \frac{\del}{\del t_n}\right) \bv
\]
exists in $\Gamma(Y, \calE/\calF)$.

Again from the formal power series computation, we see that 
$f(\bv)$ is killed by
$\frac{\del}{\del t_n}$; that is, the kernel of 
$\frac{\del}{\del t_n}$ is nonempty. We may now repeat the proof of
Proposition~\ref{P:generization}, using this last result to replace
Lemma~\ref{L:gener1} (and inspecting its proof similarly)
to produce a nonzero constant subobject $\calG$
of $\calE/\calF$. 
(The role of $L$ in the proof of Proposition~\ref{P:generization}
is played by the completed fraction field of $\calO(X \times A^{n-1}_K[0,b])$.)
Repeating the argument with $\calF$ replaced by the preimage
of $\calG$ in $\calE$, we eventually deduce that $\calE \in \ULNM_Y$.

To summarize, we have shown that 
$\calE$ is unipotent on $Y = (X \times A^{n-1}_K[0,b]) \times A^1_K[0,\eta)$
relative to $X \times A^{n-1}_K[0,b]$.
Since the restriction of $\calE$ to $X \times A^{n-1}_K[0,b] \times \{0\}$
again satisfies the convergence hypothesis (by Lemma~\ref{L:Hadamard} again),
we may invoke the induction hypothesis to obtain the desired result.
\end{proof}
\begin{remark}
The subtlety in the above proof is that the application of
Lemma~\ref{L:Hadamard} must be to a sequence without poles;
this is why we must apply it to \eqref{eq:unipotent sequence}
rather than to the sequence $\{\frac{1}{j!} \frac{\del^j}{\del t_n^j}
\bv\}$ directly.
\end{remark}

\begin{remark} \label{R:not unipotent}
We have already seen (in Lemma~\ref{L:small unip} for $X = K$;
apply Proposition~\ref{P:generization} to deduce the general case) that
without the convergence hypothesis, one can only prove that
$\calE$ is unipotent over $X \times A^n_K[0,a)$ for some
$a \in [0,1]$. Indeed, simple examples show that the
stronger conclusion of unipotence over $X \times A^n_K[0,1)$
cannot be achieved; for instance,
the log-$\nabla$-module of rank $1$ on $A^1_K[0,1)$
with generator $\bv$ satisfying
\[
\frac{\del}{\del t} \bv = \bv
\]
is only unipotent on $A^1_K[0,|p|^{1/(p-1)})$. (Its
horizontal sections are the scalar multiples of $\exp(-t) \bv$,
and the exponential only converges on the smaller disc.)
\end{remark}

\begin{defn}
Let $X$ be a smooth rigid space,
and let $\calE$ be a log-$\nabla$-module on
$X \times A^n_K[a,1)$  or
$X \times A^n_K(a,1)$
for some $a \in [0,1) \cap \Gamma^*$. We say
$\calE$ is \emph{convergent} if for any
$\eta \in (0,1)$, there exists $b \in (a,1) \cap \Gamma^*$ such that
for all $c \in [b,1) \cap \Gamma^*$, 
$\calE$ is $\eta$-convergent with respect to $t_1, \dots, t_n$ on
$X \times A^n_K[b,c]$ (relative to $X$).
\end{defn}
\begin{example}
If $\calE$ is constant, then it is convergent by Lemma~\ref{L:easy convergent}.
It follows (from the fact that $\eta$-convergence is stable under
formation of extensions) that any unipotent log-$\nabla$-module is also
convergent. It also follows that $t_1, \dots, t_n$ is an
$\eta$-convergent coordinate system on $X \times 
A^n_K[b,c]$ (relative to $X$), so we may
check $\eta$-convergence of $\calE$ on $X$
by just checking $\eta$-convergence at a set of generators.
\end{example}

\begin{remark}
If $\calE$ is the $\nabla$-module over $A^1_K[a,1)$ associated to
a finite free module $M$ over $\Gamma(A^1_K[a,1), \calO)$, then
$\calE$ is convergent if and only if $M$ is ``soluble at 1''
in the terminology of \cite[4.1-1]{chris-meb3}. (See also
\cite[\S 2.3]{christol-dwork}, where the notion of ``generic
radius of convergence'' used in \cite{chris-meb3} is introduced.)
\end{remark}

Putting Lemma~\ref{L:convergent log} together with 
Theorem~\ref{T:unipotent equiv} gives us the following characterization
of constant/unipotent $\nabla$-modules.

\begin{prop} \label{P:characterize unipotent}
Under Hypothesis~\ref{H:setup conn} with $W = \Maxspec K$,
take $a \in [0,1) \cap \Gamma^*$ and
suppose that
$\calE \in \LNM_{V \times A^n_K(a,1)}$ is convergent.
Then $\calE$ is unipotent if and only if
$\calE$ extends to a a log-$\nabla$-module
with nilpotent residues on $V \times A^n_K[0,1)$. Moreover, this
extension is unique if it exists, and $\calE$ is constant if and only
if the residues of $\frac{\del}{\del t_i}$ are all zero.
\end{prop}
\begin{proof}
If $\calE$ is unipotent, then the desired extension
exists and is unique thanks to Theorem~\ref{T:unipotent equiv}.
Conversely, if $\calE$ extends, then the extension is 
unipotent by Lemma~\ref{L:convergent log}.
\end{proof}

\begin{remark}
In case $\calE$ is already known to be isomorphic as an $\calO$-module
to the pullback of a coherent locally free $\calO$-module on $V$, one
may invoke \cite[Corollary~6.5.2]{balda-chiar} to give an
alternate derivation of Proposition~\ref{P:characterize unipotent}.
\end{remark}

\section{Monodromy of isocrystals}
\label{sec:mono}

In this section, we explain what it means for an isocrystal on
a smooth variety to have
``constant/unipotent monodromy'' along a divisor, and show that
one can ``fill in'' an overconvergent isocrystal along a divisor of
constant monodromy.

\subsection{Partial compactifications}

\begin{defn}
Let $X$ be a $k$-variety. By a \emph{partial compactification} of
$X$, we will mean a pair $(Y, j)$, where $Y$ is a $k$-variety and
$j: X \hookrightarrow Y$ is an open immersion. We do not require that
$j$ have dense image, though we will see soon (Remark~\ref{R:closed
immersion}) that this permissiveness is not so critical.
If $X$ is closed in $Y$ (e.g., if
$Y = X$ and $j = \id_X$), we say $(Y,j)$ is a \emph{trivial 
compactification}.
If the closure of $X$ in $Y$ is proper over $k$
(e.g., if $Y$ is proper over $k$), 
we say $(Y,j)$ is a \emph{full compactification}.
\end{defn}

\begin{defn} \label{D:equivalent compactifications}
Given a $k$-variety $X$ and two partial compactifications $(Y_i,j_i)$
of $X$ ($j=1,2$), put $Y_3 = Y_1 \times_k Y_2$; then $j_1$ and $j_2$ induce
an open immersion $j: X \hookrightarrow Y_3$.
Let $\overline{X}_i$ denote the Zariski closure of $X$ within
$Y_i$ for $i=1,2,3$. We write $(Y_1,j_1) \geq (Y_2,j_2)$ if
the map $\overline{X}_3 \to \overline{X}_2$ is proper; clearly this
relation is a reflexive partial ordering. In particular, we say 
that $(Y_1,j_1)$
and $(Y_2,j_2)$ are \emph{equivalent} if they are mutually comparable
under $\geq$.
Note that this does indeed give an equivalence relation;
moreover, a compactification is trivial/full if and only if it is
minimal/maximal under $\geq$.
\end{defn}

In practice, instead of checking the definition of equivalence directly,
we use the following result.
\begin{lemma} \label{L:check equivalence}
With notation as in Definition~\ref{D:equivalent compactifications},
suppose that there exists a proper map $\phi: Y_1 \to Y_2$ such that
$j_2 = \phi \circ j_1$. Then $(Y_1,j_1)$ and $(Y_2, j_2)$ are equivalent.
\end{lemma}
\begin{proof}
The map $\id_{Y_1} \times \phi: Y_1 \to Y_3$ is proper and sections
the projection $\pi_1: Y_3 \to Y_1$; we thus have regular maps
$\overline{X_3} \to \overline{X_1}$ and $\overline{X_1} \to \overline{X_3}$, 
induced by $\pi_1$ and $\id_{Y_1} 
\times \phi$,
respectively, which compose both ways
to give maps which restrict to the identity map on $X$. Since $X$ is
dense in both $\overline{X_1}$ and $\overline{X_3}$, 
the compositions really are the identity
maps; that is, the induced maps $\overline{X_3} \to \overline{X_1}$ and 
$\overline{X_1} \to \overline{X_3}$ are
isomorphisms.

In particular,
$\pi_1: \overline{X_3} \to \overline{X_1}$ is proper; since 
$\pi_2: \overline{X_3} \to \overline{X_2}$ factors as
$\phi \circ \pi_1$, it is also proper. This yields the desired
equivalence.
\end{proof}

\begin{remark} \label{R:closed immersion}
In particular, if $(Y,j)$ is a partial compactification and
$\overline{X}$ is the Zariski closure of $X$ within $Y$, then
$(Y,j)$ and $(\overline{X},j)$ are equivalent, because a closed 
immersion is proper.
\end{remark}

\begin{remark}
We have observed previously (Definition~\ref{defn:invimage}) that if
$(Y_1, j_1)$ and $(Y_2, j_2)$ are equivalent partial compactifications, and
$Y_3 = Y_1 \times_k Y_2$, then the inverse image functors
$\Isoc^\dagger(X,Y_1/K) \to \Isoc^\dagger(X,Y_3/K)$ and
$\Isoc^\dagger(X,Y_2/K) \to \Isoc^\dagger(X,Y_3/K)$ are equivalences
of categories. In other words, the category of isocrystals on $X$
overconvergent along $Y \setminus X$ depends only on the equivalence
class of the partial compactification $(Y,j)$.
\end{remark}

Since any variety can be covered by open subvarieties which are
affine and hence quasi-projective, it will be helpful to know something
similar for partial compactifications; the following lemma is a step
in this direction.
\begin{lemma} \label{L:gruson-raynaud}
Let $X$ be a quasi-projective $k$-variety. Then for any partial
compactification $(Y,j)$ of $X$, there exists a partial compactification
$(Y',j')$ with $Y'$ quasi-projective and a proper map $\phi: Y' \to Y$
such that $j = \phi \circ j'$. In particular, the two partial
compactifications $(Y,j)$ and $(Y',j')$ are equivalent.
\end{lemma}
\begin{proof}
This is precisely the statement (restricted from algebraic spaces to
varieties) of the quantitative Chow's lemma of Gruson-Raynaud
\cite[Corollaire~5.7.14]{gruson}.
\end{proof}

\subsection{Smooth varieties and small frames}
\label{subsec:small frames}

We now focus attention on isocrystals on smooth varieties; it will be
convenient to handle them using a special sort of frame.
\begin{defn}
A \emph{small frame} is a frame $(X,Y,P,i,j)$ in which $Y = P_k$,
the map $i$ is the identity,
and $Y \setminus X$ is the zero locus of some regular function on $Y$.
We will drop $Y$ and $i$ from the
notation for a small frame, denoting it by $(X,P,j)$. Note that
in any small frame, $X$ must be smooth, since $X$ is open in $P_k$
and $P$ is smooth in a neighborhood of $X$.
\end{defn}

In order to make much use of small frames, we need the following lemma.
\begin{lemma} \label{L:enough small frames}
Let $j: X \hookrightarrow Y$ be an open immersion of $k$-varieties,
with $X$ dense in $Y$.
Then there exists a blowup $Y' \to Y$ centered in $Y \setminus X$,
an open cover $U_1, \dots, U_n$ of $Y'$, and
for $i=1,\dots,n$, a partial compactification
$(Y_i, j_i)$ of $X \cap U_i$, enclosed by a small frame, such that 
$Y_i$ admits a proper morphism $\phi_i$ to 
$Y' \cap U_i$ with $j = \phi_i \circ j_i$ on $X \cap U_i$.
In particular, $(Y_i,j_i)$ is equivalent to $(Y' \cap U_i,j)$.
\end{lemma}
\begin{proof}
By blowing up in $Y \setminus X$, we may reduce to the case where all components
of $Y \setminus X$ have codimension 1 in $Y$.
By then passing to open affine covers,
we may reduce to the case where $X$ and $Y$ are affine (and 
$Y \setminus X$ is still a divisor).
By a theorem of Arabia \cite[Th\'eor\`eme~1.3.1]{arabia}
(generalizing a theorem of Elkik \cite{elkik} in the case of
$K$ discretely valued), there exists a smooth affine
scheme $\tilde{X}$ over $\gotho_K$ with $\tilde{X} \times_{\gotho_K} k
\cong X$. Choose an embedding of $\tilde{X}$ into a projective space
$\PP^n_{\gotho_K}$ and let $P$ be the formal completion of the
projective closure of $\tilde{X}$ in $\PP^n_{\gotho_K}$.

Choose a closed immersion $Y \hookrightarrow \AAA^l_k$, where the latter
has coordinates $x_1, \dots, x_l$. Then along the rational map
$P_k \dashrightarrow Y \hookrightarrow \AAA^l_k$ induced by the
isomorphism between the two copies of $X$, each of
$x_1, \dots, x_l$ pulls back to a rational function $f_1, \dots, f_l$
on $P_k$. For some $m>0$, these functions can be written as quotients
of homogeneous polynomials of degree $m$ (i.e., sections of
$\calO(m)$); lift these polynomials to homogeneous polynomials of degree
$m$ over $\gotho_K$. The resulting rational functions define a rational
map $P \dashrightarrow \widehat{\AAA^l_{\gotho_K}}$; let $P'$ denote 
the closure of the graph of this rational map. Then $P'_k$ is a
partial compactification of $X$ admitting a proper map to $Y$,
and the complement $P'_k \setminus X$ is the zero locus of a regular
function;
we can cover $P'$ with affines to obtain the desired small frames.
\end{proof}

\begin{remark}
Lemma~\ref{L:enough small frames} may be interpreted as saying that
any isocrystal can be described entirely using small frames. However, 
this does not assert by itself that one can reconstruct the whole theory
of isocrystals using only small frames, since functoriality is defined
by passing to a restriction from a product frame, which is not small.
One could get around this using sophisticated ``lifting lemmas'' of the
sort given in \cite{arabia}; this would amount to giving a development
of isocrystals from the point of view of Monsky and Washnitzer's
``formal cohomology'' (see \cite{mw} for the construction, and
\cite[Section~2.5]{ber2} for its relationship to Berthelot's
construction). We will not give such a development here.
\end{remark}

\subsection{Monodromy: a restricted definition}
\label{subsec:monodromy}

\begin{lemma} \label{L:cohen}
Let $A$ be a noetherian ring, such that $A$
is complete
with respect to the $x$-adic topology for some $x \in A$ not a zero divisor, 
and let
$R$ be a subring of $A$.
Suppose that $B = A/xA$ is formally smooth over $R$. Then there is
an isomorphism $A \cong B\llbracket x \rrbracket$ sending $x$ to $x$,
whose composition with the quotient $B \llbracket x \rrbracket
\to B \llbracket x \rrbracket / x B \llbracket x \rrbracket
\cong B$ gives the quotient map $A \to A/xA \cong B$.
\end{lemma}
\begin{proof}
The proof is as in \cite[Lemma~II.1.2]{hartshorne}, 
except there $R$ is taken to be a field (but the argument does not
change). See also \cite[Expos\'e~III, 5.6]{sga1}.
\end{proof}

\begin{hypo} \label{H:mono setup}
Let $X \hookrightarrow Y$ be an open immersion
of smooth affine $k$-varieties, with $X$ dense in $Y$ and $Z = Y \setminus X$
also smooth.
Suppose that there exists a small frame $(X,P,j)$ enclosing $Y$,
and that there exists $f \in
\Gamma(P, \calO_P)$ which cuts out $Z$ within $Y$, such that
$df$ generates a direct summand of $\Omega^1$
in a neighborhood of $Z$.
Let $Q$ be the zero locus of $f$ on $P$.
\end{hypo}

\begin{lemma} \label{L:tube isom}
Under Hypothesis~\ref{H:mono setup}, there exists an isomorphism
$\phi: ]Z[_Q \times A^1_K[0,1) \to ]Z[_P$.
\end{lemma}
\begin{proof}
Apply Lemma~\ref{L:cohen} to produce an isomorphism
$\Gamma(]Z[_Q, \gotho) \llbracket t_1 \rrbracket
\cong \Gamma(]Z[_P, \gotho)$.
This yields the desired map. (This can also be proved using the strong
fibration theorem; compare the proof of Lemma~\ref{L:direct extension}.)
\end{proof}

\begin{defn} \label{D:monodromy 1}
Under Hypothesis~\ref{H:mono setup}, let $\calE$ be an isocrystal
on $X$ overconvergent along $Y \setminus X$. We confound $\calE$ with
its realization on the small frame $F = (X,P,j)$; the latter is a
$\nabla$-module on a strict neighborhood $V$ of $]X[_P$ in $]Y[_P$.
Since $]Y[_P = P_K$ is an affinoid space, by 
Lemma~\ref{L:strict neighborhood}, $V \cap ]Z[_P$ contains a
subspace of the form
\[
\{y \in P_K: |f(y)| \geq \lambda\}
\]
for some $\lambda \in (0,1) \cap \Gamma^*$. Under $\phi^{-1}$, 
such a space maps to
$]Z[_Q \times A^1_K[\lambda,1)$,
so $\calE$ restricts to a $\nabla$-module on $]Z[_Q \times A^1_K[\lambda,1)$,
which is convergent thanks to Proposition~\ref{P:taylor} (applied with 
$g = f$). We say that $\calE$ has
\emph{constant/unipotent monodromy along $Z$}
(with respect to $f, \phi$) if 
$\calE$ is constant/unipotent over $]Z[_Q 
\times A^1_K[\lambda,1)$ for some $\lambda \in (0,1) \cap \Gamma^*$.
\end{defn}

So far, the definition of the phrase
``$\calE$ has constant/unipotent monodromy along $Z$'' depends on
the choices of the frame $(X,P,j)$, the map $\phi$, and the function $f$.
To eliminate these dependencies, we make the usual argument of passing to
a product frame, but since the latter is not a small frame,
some care is required.

\begin{prop}
Under Hypothesis~\ref{H:mono setup},
let $(X, P',j')$ be another small frame satisfying the same hypotheses
(with corresponding objects denoted by primes). 
Let $\calE'$ be the realization of $\calE$ on $(X,P',j')$. Then
$\calE$ has constant/unipotent monodromy along $Z$ if and only
if $\calE'$ has constant/unipotent monodromy along $Z$.
\end{prop}
\begin{proof}
We first note that 
by Proposition~\ref{P:characterize unipotent}, $\calE$ has constant monodromy
along $Z$ if and only if $\calE$ extends from some $]Z[_Q \times 
A^1_K[\lambda,1)$ to a $\nabla$-module on $]Z[_Q \times A^1_K[0,1)$.
Similarly, $\calE$ has unipotent monodromy along $Z$ if and only if
$\calE$ admits a filtration $0 = \calE_0 \subset \calE_1 \subset
\cdots \subset \calE_l = \calE$ whose successive quotients extend to
$\nabla$-modules on $]Z[_Q \times A^1_K[0,1) = ]Z[_P$.

Suppose now that $\calE$ has unipotent monodromy along $Z$.
By passing to an affine cover, we may assume that there
exist $x_1,\dots,x_m \in \Gamma(P, \calO_P)$ and
$x'_1, \dots, x'_m \in \Gamma(P', \calO_{P'})$ whose differentials
generate $\Omega^1$ on $P$ and $P'$, respectively,
such that $x_i \equiv x'_i$ as elements of $\Gamma(Y, \calO) = 
\Gamma(P_k, \calO)
= \Gamma(P'_k, \calO)$.
Put $P'' = P \times P'$, 
put $j'' = j \times j'$, 
put $t_i = x_i - x'_i \in \Gamma(P'', \calO_{P''})$, 
and let $\calE''$ be the realization of $\calE$ on $(Y,P'',j'')$.
On one hand, $\calE''$ is isomorphic to the pullback $\pi_1^* \calE$ along
the projection $P'' \to P$; so on the intersection of $]Z[_{P''}$ with
some strict neighborhood of $]X[_{P''}$ in $]Y[_{P''}$,
$\calE''$ admits a filtration $0 = \calE''_0 \subset \calE''_1 \subset
\cdots \subset \calE''_l = \calE''$ whose successive quotients extend to
$\nabla$-modules on $]Z[_{P''}$.
On the other hand, $\calE''$ is also isomorphic to the pullback
$\pi_2^* \calE'$, and in fact we can recover $\calE'$ from $\calE''$ by
restricting to a component of 
the subspace $t_1 = \cdots = t_m = 0$ of $P''$. In particular,
we obtain a filtration $0 = \calE'_0 \subset \calE'_1 \subset
\cdots \subset \calE'_l = \calE'$ whose successive quotients extend to
$\nabla$-modules on $]Z[_{P'} = ]Z[_{Q'} \times A^1_K[0,1)$. Hence
$\calE'$ also has unipotent monodromy along $Z$. Moreover, if $\calE$
actually has constant monodromy along $Z$, then we can take the
filtration of $\calE$ to be the trivial one $0 = \calE_0 \subset \calE_1 = 
\calE$, move it through the above argument, and deduce that $\calE'$
has constant monodromy along $Z$.
\end{proof}

\begin{remark} \label{R:constant monodromy}
If $\calE$ extends to a convergent isocrystal on $Y$, then $\calE$ has
constant monodromy along $Z$ by 
Proposition~\ref{P:characterize unipotent}. We will prove a converse
of this observation; see Theorem~\ref{T:no log extension}.
\end{remark}

\begin{remark}
As noted in Remark~\ref{R:Tannakian}, one could in principle
construct a local monodromy representation (along $Y \setminus X$)
for an isocrystal on $X$ overconvergent along $Y \setminus X$.
We will defer doing so to a subsequent paper.
\end{remark}

\subsection{Monodromy: a general definition}

We now wish to extend the definition of constant/unipotent monodromy;
first we make some comments about the existing definition.
\begin{remark} \label{R:reductions}
Under Hypothesis~\ref{H:mono setup}, let $\calE$ be the realization,
on a fixed small frame $F$, of an isocrystal on
$X$ overconvergent along $Z = Y \setminus X$. Then the following are true.
\begin{itemize}
\item
Let $U_1, \dots, U_n$ be an open cover of $Y$. Then
$\calE$ has constant/unipotent monodromy along $Z$ if and only if for
$i=1, \dots, n$, the restriction of $\calE$ to $U_i \cap X$ has
constant/unipotent monodromy
along $U_i \cap Z$; this follows from Corollary~\ref{C:unipotent
local} applied to the admissible cover $\{]U_i \cap Z[\}$ of $]Z[$.
\item
Let $K'$ be a field containing $K$ which is complete under an extension
of $|\cdot|$. Then $\calE$ has constant/unipotent monodromy along
$Z$ if and only if this is true after changing the base field to $K'$;
this follows from Proposition~\ref{P:generization}.
\item
Let $U$ be an open subscheme of $Y$ such that $U \cap Z$ is dense in $Z$.
Then $\calE$ has constant/unipotent monodromy along $Z$ if and only
if the restriction of $\calE$, to an isocrystal on $U \cap X$
overconvergent along $U \cap Z$, has constant/unipotent monodromy along
$U \cap Z$; this also follows from Proposition~\ref{P:generization},
or more precisely from Corollary~\ref{C:generization}.
\item
If $\calE$ extends to a convergent isocrystal on $Y$, then $\calE$
has constant monodromy along $Z$, by Proposition~\ref{P:characterize
unipotent}.
\end{itemize}
\end{remark}

\begin{defn} \label{D:monodromy 2}
Let $X \hookrightarrow Y$ be an open immersion of smooth $k$-varieties,
and let $\calE$ be an isocrystal on $X$ overconvergent along $Z = Y
\setminus X$. 
We say $\calE$ has \emph{constant/unipotent monodromy along $Z$} if for
any extension field $k'$ of $k$, any field $K'$ containing $K$ which is
complete under an extension of $|\cdot|$ with residue field $k'$,
and any small frame $(U,P,j)$ over $K$ enclosing an open subset $V = P_k$ of $Y$
(with $U = V \cap X$)
which satisfies Hypothesis~\ref{H:mono setup} (i.e., $V \setminus U$
is smooth and is the zero locus of some $f \in
\Gamma(P, \calO_P)$),
the realization of $\calE$ on $(U,P,j)$ has constant/unipotent monodromy along
$V \setminus U$. By virtue of Remark~\ref{R:reductions}, this agrees
with Definition~\ref{D:monodromy 1} when they both apply; also, the analogue
of Remark~\ref{R:reductions} holds for this expanded definition.
\end{defn}
\begin{remark}
The checking over extension fields is only necessary when $k$ is imperfect:
when $k$ is perfect, $Z$ (being reduced, thanks to our running hypothesis that
all $k$-varieties are reduced) is generically smooth, so we may sample
on a suitable open subset of $Y$ without enlarging $k$. However, if
$k$ is imperfect, then $Z$ may fail to be geometrically reduced, and 
one must extend $k$ in order to guarantee that the underlying reduced
subscheme is generically smooth. This will require us to do a bit
of work in the case of $k$ imperfect in order to complete the proof
of the extension theorem (Theorem~\ref{T:no log extension}).
\end{remark}

An important property of the definition of constant/unipotent monodromy
is its ``codimension 1 nature''.
\begin{prop} \label{P:codimension 1}
Let $U \hookrightarrow X \hookrightarrow Y$ be open immersions
of smooth $k$-varieties, such that $Y \setminus X$ has codimension at least $2$
in $Y$. Let $\calE$ be an isocrystal on $U$ overconvergent
along $Y \setminus U$. Then $\calE$ has constant/unipotent monodromy
along $Y \setminus U$ if and only if 
$\calE$ has constant/unipotent
monodromy along $X \setminus U$.
\end{prop}
\begin{proof}
There is no harm in shrinking $U$ so that $Y \setminus U$ becomes purely
of codimension 1, as $\calE$ automatically has constant monodromy along any
added component. In this case, $X \setminus U$ is dense in $Y \setminus U$,
so we obtain the desired equivalence
as in Remark~\ref{R:reductions}.
\end{proof}

\section{Monodromy and extensions}
\label{sec:ext}

In this section, we clarify the relationship between extendability of
an isocrystal and the property of having constant monodromy along
some boundary variety.

\subsection{An extension lemma}

We now prove a lemma about extending $\nabla$-modules 
in a key geometric setting. To avoid having to repeat effort, we set up
the lemma so that it also handles log-$\nabla$-modules with nilpotent
residues; hence the somewhat complicated statement.
\begin{lemma} \label{L:direct extension}
Let $V \hookrightarrow U \hookrightarrow 
X \hookrightarrow Y$ be open immersions of 
$k$-varieties such that $X$ is smooth, $V$ is dense in $Y$, $X \setminus V$
is a strict normal crossings divisor on $X$, and $X \setminus U$ is a single
component of $X \setminus V$.
Suppose further that there exist:
\begin{itemize}
\item
a small frame
$F = (X,P,j)$ enclosing $(X,Y)$;
\item
functions $f_1, \dots, f_r \in \Gamma(P, \calO_P)$
whose zero loci cut out the components of the 
closure of $X \setminus V$ in $Y$, with $f_1$ cutting out $X \setminus U$;
\item
functions $f_{r+1}, \dots, f_n \in \Gamma(P, \calO_P)$ such that
$df_1, \dots, df_n$
freely generate $\Omega^1$ in a neighborhood of $X$;
\item 
a function $g \in \Gamma(P, \calO_P)$ whose zero locus cuts out
$Y \setminus X$ within $Y$.
\end{itemize}
Then the following results hold.
\begin{enumerate}
\item[(a)]
Let $\calE$ be a $\nabla$-module on a strict neighborhood of $]U[_P$ in
$]Y[_P = P_K$ representing an isocrystal on $U$
overconvergent along $Y \setminus U$.
Then $\calE$ has
constant monodromy along $X \setminus U$ if and only if $\calE$
extends to an isocrystal on $X$ overconvergent along $Y \setminus X$.
\item[(b)]
Let $\calE$ be a log-$\nabla$-module with nilpotent residues
on a strict neighborhood of $]U[_P$
in $P_K$ with respect to $f_1, \dots, f_r$, whose restriction to a strict
neighborhood of $]V[_P$ in $P_K$ represents an isocrystal on $V$
overconvergent along $Y \setminus V$.
Then $\calE$ has unipotent
monodromy along $X \setminus U$ if and only if $\calE$
extends to a log-$\nabla$-module with nilpotent residues
on a strict neighborhood of $]X[_P$
in $P_K$ with respect to $f_1, \dots, f_r$.
\item[(c)]
In both (a) and (b), the implied restriction functor is fully faithful:
that is, morphisms between $\calE$ and $\calE'$ always uniquely induce
morphisms
on their extensions.
\end{enumerate}
\end{lemma}
\begin{proof}
Let $P'$ be the zero locus of $f_1$ on $P$.
Let $F'$ be the frame $(X \setminus U,P',j')$,
and let $f'_2, \dots, f'_n$ be the restrictions of $f_2,\dots, f_n$
to $P'$. 
Put $Z = Y \setminus U$.
By the strong fibration
theorem (Proposition~\ref{P:strong fibration}), there exists a strict
neighborhood of $]X \setminus U[_{P \times P'}$ in $]Z[_{P \times P'}$
isomorphic on one hand to a strict neighborhood $V_1$ of
$]X \setminus U[_{P \times \widehat{\AAA^{n-1}}} 
\cong ]X \setminus U[_P \times A^{n-1}_K[0,1)$ in
$]Z[_{P \times \widehat{\AAA^{n-1}}} = ]Z[_P \times A^{n-1}_K[0,1)$ 
via the functions
$f_2-f'_2, \dots, f_n - f'_n$, and on the other hand to a strict 
neighborhood $V_2$ of 
$]X \setminus U[_{P' \times \widehat{\AAA^n}} =
]X \setminus U[_{P'} \times A^1_K[0,1) \times A^{n-1}_K[0,1)$ in
$]Z[_{P' \times \widehat{\AAA^n}} = ]Z[_P \times A^1_K[0,1) \times
A^{n-1}_K[0,1)$ via the functions
$f_1, f_2 - f'_2, \dots, f_n - f'_n$.
If we restrict the resulting isomorphism $V_1 \to V_2$ to the
inverse image of $0 \in A^{n-1}_K[0,1)$ in both factors,
we get an isomorphism between a strict neighborhood of
$]X \setminus U[_P$ in $]Z[_P$ with a strict neighborhood of
$]X \setminus U[_{P'} \times A^1_K[0,1)$ in
$]Z[_{P'} \times A^1_K[0,1)$, whose composition with the projection
$]Z[_{P'} \times A^1_K[0,1) \to A^1_K[0,1)$ is precisely $f_1$.

By assumption, $\calE$ is defined on some subset of
$P_K$ of the form
\[
V_\lambda = \{x \in P_K: |f_1(x)| \geq \lambda, |g(x)| \geq \lambda\}
\]
with $\lambda \in (0,1) \cap \Gamma^*$,
and 
its restriction to $V_\lambda \cap ]X \setminus U[_P$ is in case (a)
a constant $\nabla$-module and in case (b) a unipotent log-$\nabla$-module. 
Now pass $\calE$ over to a strict neighborhood of
$]X \setminus U[_{P'} \times A^1_K[0,1)$ in
$]Z[_{P'} \times A^1_K[0,1)$; then for each closed subinterval $[a,b]$
of $(\lambda, 1)$, $\calE$ is defined on $V_0 \times A^1_K[a,b]$
for some strict neighborhood $V_0$ of $]X \setminus U[_{P'}$ in
$]Z[_{P'}$. By Proposition~\ref{P:generization2}, there exists another
strict neighborhood $V_1$ of $]X \setminus U[_{P'}$ in $]Z[_{P'}$
such that $\calE$ becomes constant/unipotent on $V_1 \times A^1_K[a,b]$.
By Theorem~\ref{T:unipotent equiv}, this restriction of $\calE$
extends in case (a) to a $\nabla$-module, or in case (b) to a
log-$\nabla$-module with nilpotent residues,
on $V_1 \times A^1_K[0,b]$, which
we may glue with the original $\calE$ to extend it to a strict neighborhood
of $]X[_P$ in $]Y[_P$. 
The assertion of (c) follows from Corollary~\ref{C:horizontal}.

Finally, we check the overconvergence of the extension in the case (a),
by verifying the condition of Proposition~\ref{P:taylor}; that is,
we claim that our extension is $\eta$-convergent with respect to
$f_1, f_2, \dots, f_n$ on some affinoid strict neighborhood of 
$]X[_P$ in $]Y[_P$ (which may depend on $\eta$). 
We need only verify the $\eta$-convergence condition for each of a 
set of generating sections; by Proposition~\ref{P:taylor},
we already know this on some $V_\lambda$.
Now run the aforementioned construction for a choice of $[a,b]$ with
$\eta < a$. Then the fact that $\calE$ is constant on $V_1 \times A^1_K[a,b]$
means (by Lemma~\ref{L:easy convergent}) that the extension of 
$\calE$ to $V_1 \times A^1_K[0,b]$ is $\eta$-convergent.
This yields $\eta$-convergence of the extension of $\calE$ to a strict
neighborhood of $]X[_P$ in $]Y[_P$, as desired.
\end{proof}

\subsection{Extension of overconvergent isocrystals}

With Lemma~\ref{L:direct extension} in hand, we can now prove a 
definitive theorem about extending overconvergent isocrystals.
\begin{theorem} \label{T:no log extension}
Let $U \hookrightarrow X \hookrightarrow Y$ be open immersions
of $k$-varieties, such that $X$ is smooth and $U$ is dense in $Y$.
Let $\calE$ be an isocrystal on $U$ overconvergent along
$Y \setminus U$. Then $\calE$
has constant monodromy along $X \setminus U$
if and only if $\calE$ extends to an isocrystal on 
$X$ overconvergent along $Y \setminus X$.
Moreover, the functor $\Isoc^\dagger(X,Y/K) \to \Isoc^\dagger(U,Y/K)$ is fully
faithful, so the extension is unique if it exists.
\end{theorem}
\begin{proof}
As in Remark~\ref{R:constant monodromy}, if $\calE$ extends, it
must have constant monodromy along $X \setminus U$. We will prove
the converse and the full faithfulness
under several sets of hypotheses, culminating in the
unrestricted form. 

To begin with, suppose that
$X \setminus U$ is a smooth divisor on $X$. By applying 
Lemma~\ref{L:enough small frames} (allowing $Y$ to be replaced by a blowup
centered in $Y \setminus X$), then passing to an open cover of $Y$
and replacing each open subset of $Y$ 
by an equivalent partial compactification (of the subset of $X$ it contains),
we may reduce the desired assertion to a collection of instances
of  Lemma~\ref{L:direct extension}, in which we fill in one component of
$X \setminus U$ at a time. (Note that part (c) of the lemma ensures
that the extensions produced can be glued back together.)

Next, suppose that $k$ is perfect but $U,X,Y$ are not
further restricted. If $X \setminus U$ is nonempty,
we can find a smooth closed point $x$ on (the reduced
subscheme underlying) $X \setminus U$, since the latter is also
geometrically reduced. Let $Z$ be the unique component of $X \setminus U$
passing through $x$, and let $D$ be an irreducible divisor of $X$ containing
$Z$ which is smooth in a neighborhood $V$ of $x$. (For instance,
choose functions $t_1, \dots, t_r$ cutting out $Z$ within $X$
whose differentials form part of a basis of $\Omega^1$ in 
a neighborhood of $x$, then take $D$ to be the component of
the zero locus of $t_1$ passing through $x$.)
Then either $D = Z$, or $D \setminus Z$ is dense in $D$. In either case,
the restriction of $\calE$ to $X \setminus D$ has constant monodromy
along $D$: in the former case this is by hypothesis, whereas in the 
latter case this is automatic.

Let $Z'$ be the union of the components of $X \setminus U$ other than
$Z$, together with the nonsmooth locus of $D$. 
By the previously treated case, $\calE$ extends to an isocrystal on
$V \setminus Z'$ overconvergent along $Y \setminus (V \setminus Z')$,
and the corresponding restriction functor is fully faithful.
Since $x \in V \setminus Z'$, we may glue to obtain an extension of
$\calE$ to an open subset of $X$ which is strictly larger than $U$.
By noetherian induction, repeating this process eventually yields
an extension of $\calE$ to $X$ and the full faithfulness of the
restriction functor.

Finally, suppose that $k$ is arbitrary.
In this case, we can still run the previous
argument at the expense of replacing $k$ by a finite radicial extension.
It thus suffices to show the following: suppose that $K' = K(y^{1/p})$
for some $y \in \gotho_K$ whose image in $k$ is not
a $p$-th power, and that the assertion of the theorem holds
for $U,X,Y$ over $K'$. Then it also holds for $U,X,Y$ over $K$.
(Namely, with this result in hand,
we can enlarge the residue field
from $k$ to any desired finite radicial extension by a sequence of
such extensions of $K$, then back down the tower to deduce the
theorem.)

Since everything under consideration is local,
we may assume thanks to Lemma~\ref{L:enough small frames}
that $(X,Y)$ is enclosed by a small frame $(X,P,j)$.
Take $\calE \in \Isoc^\dagger(U,Y/K)$ with constant monodromy
along $X \setminus U$.
For $V$ an affinoid strict neighborhood of $]U[$ in $]Y[$, put
$A_V = \Gamma(V, \calO)$ and $M_V = \Gamma(V, \calE)$.
For $W$ an affinoid strict neighborhood of $]X[$ in $]Y[$, put
$B_W = \Gamma(W, \calO)$.
For everything in sight, insert a prime to denote tensoring
with $K'$ over $K$.
We have (by applying the theorem over $K'$)
that for some affinoid strict neighborhood
$V$ of $]U[$ in $]Y[$, there exists an affinoid strict
neighborhood $W$ of $]X[$ in $]Y[$ containing $V$ and
a finitely generated $B'_W$-submodule $N'_W$ of 
$M_V'$, stable under
$\nabla$ and satisfying $N'_W \otimes_{B'_W} A'_V = M'_V$.
By the full faithfulness
of restriction from $X$ to $U$ over $K'$, $N_W'$ is uniquely determined by
these conditions.

Put $N_W = N'_W \cap M_V$; then $N_W$ is
a $B_W$-submodule of $M_V$ which is stable under $\nabla$.
We will show that $N_W$ is finitely generated and that
$N_W \otimes_{B_W} A_V = M_V$.
It suffices to check this after enlarging $K$
and $K'$ to contain a primitive $p$-th root of unity $\zeta_p$ (since
$K(\zeta_p)$ and $K'$ are linearly disjoint over $K$, by the
hypothesis on $y$). In this case, $K'$ becomes Galois with group
$G = \Gal(K'/K)$, which we identify with $\ZZ/p\ZZ$ by
declaring that $e \in \ZZ/p\ZZ$ carries $y^{1/p}$
to $\zeta_p^e y^{1/p}$.

Thanks to Proposition~\ref{P:pullback0} and the fact that $G$ acts
trivially modulo $\gothm_K$, we obtain a canonical action of $G$
on $M'_V$ with invariants $M_V$ (at least after shrinking $V$, which
is harmless). By the uniqueness of $N'_W$,
$N'_W$ also carries an action of $G$. For $i=0, \dots, p-1$ and
$\bv \in M'_V$, set
\[
f_i(\bv) = (y^{1/p})^{-i} \sum_{e \in \ZZ/p\ZZ} \zeta_p^{-ei} \bv^e.
\]
Then each $f_i$ carries $M'_V$ into $M_V$, and so carries
$N'_W$ into $N_W$. 

It is clear that the natural map $N_W \otimes_K K' \to N'_W$ is injective.
On the other hand, for $\bv = \sum_{l=0}^{p-1} (y^{1/p})^l \bv_l \in 
N'_W$, with each $\bv_l \in M_V$, we have $\bv_l = f_l(\bv) \in N_W$ as in
the previous paragraph. Hence $N_W \otimes_K K' \to N'_W$ is also surjective,
so
\[
(N_W \otimes_{B_W} A_V) \otimes_K K' = N'_W \otimes_{B'_W} A'_V = M'_V 
= M_V \otimes_K K'
\]
and so $N_W \otimes_{B_W} A_V = M_V$ by Galois descent.

Moreover, if $\bw \in M_V$ and $\bv_j
\in M'_V$ satisfy $\sum b_j \bv_j = \bw$ for some $b_j \in B'_W$,
write $b_j = \sum_{l=0}^{p-1} b_{j,l} (y^{1/p})^{-l}$ with
$b_{j,l} \in A_V$ (resp.\ $b_{j,l} \in B_W$); we then have
\begin{align*}
p\bw &= f_0(\bw) \\
&= \sum_j \sum_{e \in \ZZ/p\ZZ} b_j^e \bv_j^e \\
&= \sum_{j} \sum_{e \in \ZZ/p\ZZ} \sum_{l=0}^{p-1}
b_{j,l} \zeta_p^{-el} (y^{1/p})^{-l}  \bv_j^e \\
&= \sum_{j} \sum_{l=0}^{p-1} b_{j,l} f_l(\bv_j).
\end{align*}
That is, $\bw$ is also a $B_W$-linear
combination of the $f_l(\bv_j)$.
Consequently, given any finite set of generators of $N'_W$ over 
$B'_W$ which also generate $M'_V$ over $A'_V$, their images under
all of the $f_i$ generate $N_W$ over $B_W$.

Since $N_W$ is finitely generated and $N_W \otimes_{B_W} A_W = M_V$,
we can extend $\calE$ to a $\nabla$-module on $W$;
its overconvergence can be
checked after tensoring with $K'$. Thus $\calE$ extends to an
element of $\Isoc^\dagger(X,Y/K)$.

To obtain the extension of horizontal sections, suppose
$\bv \in M_V$ is horizontal. Then on one hand $\bv \in N'_W$
by the assertion of the theorem over $K'$; on the other hand,
$\bv$ is $G$-invariant. Hence $\bv \in N_W$, i.e., $\bv$ extends
to $X$ as desired.
\end{proof}

\begin{remark}
The full faithfulness of restriction to an open subscheme generalizes
a result of \'Etesse \cite[Th\'eor\`eme 4]{etesse}, by eliminating the
restrictions that $K$ be discretely valued and that the isocrystals
carry Frobenius structures. On the other hand, the extension criterion
seems to be new in essentially all cases except perhaps on curves (where
it is straightforward).
\end{remark}

\subsection{Consequences of overconvergent extension}

Before proceeding to the logarithmic situation, we pause to
record some consequences of Theorem~\ref{T:no log extension}. Some
of these may be of independent interest.

We first give a result about extending sub-isocrystals.
\begin{prop} \label{P:extend subobject}
Let $U \hookrightarrow X \hookrightarrow Y$ be open immersions
of $k$-varieties, such that $X$ is smooth and $U$ is dense in $Y$.
Let $\calE$ be an isocrystal on $X$ overconvergent along $Y \setminus X$, 
and let
$\calF$ be a sub-isocrystal of $\calE$ over $U$ overconvergent
along $Y \setminus U$. Then $\calF$
is the restriction to $U$ of a sub-isocrystal of $\calE$ over $X$
overconvergent along $Y \setminus X$.
\end{prop}
\begin{proof}
By Theorem~\ref{T:no log extension}, $\calE$ has constant monodromy along
$X \setminus U$, as then does $\calF$ by
Proposition~\ref{P:abelian category}, so $\calF$ extends to
an isocrystal $\calG$ on $X$ overconvergent along $Y \setminus X$.
By the full faithfulness component of Theorem~\ref{T:no log extension},
the inclusion $\calG \hookrightarrow \calE$ extends from $U$ to $X$.
This yields the desired result.
\end{proof}
\begin{remark} \label{R:legendre}
This situation should be contrasted with the situation that arises
when proving that the forgetful functor from overconvergent to
convergent $F$-isocrystals (isocrystals with Frobenius structures;
see Definition~\ref{D:Frobenius})
is fully faithful, as in \cite{me-full}.
There one does not have an analogue of Proposition~\ref{P:extend subobject},
as an overconvergent $F$-isocrystal can have
nonconstant
convergent subcrystals that do not descend to the overconvergent category.
For instance, if $f: X \to B$ is the Legendre family of elliptic curves
minus the supersingular fibres,
then $R^1 f_* \calO_X$ is a rank two overconvergent $F$-isocrystal on $B$
which has a unit-root subobject in the convergent category, but not in
the overconvergent category. (If it had a unit-root subobject in the
overconvergent category, then by Proposition~\ref{P:extend subobject}, 
it would also have
a unit-root subobject even if the supersingular fibres were not excluded,
which is absurd.)
\end{remark}

We next observe that isocrystals extend across holes of codimension
at least 2.
\begin{prop} \label{P:codimension 2}
Let $U \hookrightarrow X \hookrightarrow Y$ be open immersions
of $k$-varieties, such that $X$ is smooth, $U$ is dense in $Y$,
and $X \setminus U$ has codimension at least $2$ in $X$.
Then the restriction functor $\Isoc^\dagger(X,Y/K) \to \Isoc^\dagger(U,Y/K)$ is
an equivalence of categories.
\end{prop}
\begin{proof}
The restriction functor is fully faithful by Theorem~\ref{T:no log
extension}, so we must show that it is essentially surjective.
Let $\calE$ be an isocrystal on $U$ overconvergent along $Y \setminus U$.
Then applying Proposition~\ref{P:codimension 1} shows that
$\calE$ has constant monodromy along $X \setminus U$ if and only if it
has constant monodromy along the empty scheme. The latter is
vacuously true, so $\calE$ extends to $X$. This yields
the desired essential surjectivity.
\end{proof}
\begin{remark}
The restriction that $X$ be smooth is critical, just as
the regularity
restriction is critical in the Zariski-Nagata purity theorem;
one can construct counterexamples in the nonsmooth case much as
in the algebraic de Rham setting, e.g., by taking the rank 1
$\nabla$-module defined by $\nabla(\bv) = \bv \otimes \frac{dx}{2x}$ on the
surface $z^2 = xy$ away from $x=y=z=0$.
On the other hand, \cite[Expos\'e~X, Th\'eor\`eme~3.1]{sga1}
gives another form of the purity theorem which we are unable
to analogize using our techniques; we leave it as a question.
\end{remark}

\begin{question}
Let $U \hookrightarrow X \hookrightarrow Y$ be open immersions
of $k$-varieties, such that $X$ is a local complete intersection,
$U$ is dense in $Y$, and $X \setminus U$ has codimension at least
$3$ in $X$. Is the restriction functor
$\Isoc^\dagger(X,Y/K) \to \Isoc^\dagger(U,Y/K)$ an equivalence of categories?
This has been verified explicitly in some special cases by
Tsuzuki (private communication).
\end{question}

Using Proposition~\ref{P:codimension 2}, we can analogize the
invariance of the algebraic fundamental group under a blowup.
\begin{prop}
Let $f: Y \to X$ be a proper birational morphism of smooth
$k$-varieties, and let $\calE$ be an overconvergent isocrystal on $Y$.
Then there exists an overconvergent isocrystal $\calF$ on $X$
such that $\calE \cong f^* \calF$.
\end{prop}
\begin{proof}
Since $f$ is birational, there is an open subset $U$ of $X$,
whose complement has codimension at least 2 in $X$, on which $f$ is
an isomorphism. The restriction of $\calE$ to $U$ extends to
an overconvergent isocrystal $\calF$ on $X$ by
Proposition~\ref{P:codimension 2}; the isomorphism
$\calE \cong f^* \calF$ over $U$ extends to $X$ by
the full faithfulness aspect of Theorem~\ref{T:no log extension}.
\end{proof}

Finally, we give a result to the effect that ``overconvergence is 
contagious''.
\begin{prop} \label{P:contagion}
Let $U \hookrightarrow X \hookrightarrow Y$ be open immersions of 
$k$-varieties, such that $X$ is smooth and $U$ is dense in $Y$.
Let $\calE$ be a convergent isocrystal on $X$
whose restriction to $\Isoc^\dagger(U,X/K)$ is isomorphic to the restriction
of an isocrystal on $U$ overconvergent along $Y \setminus U$.
Then $\calE$ itself is the restriction to $\Isoc(X/K)$ of an isocrystal
on $X$ overconvergent along $Y \setminus X$.
\end{prop}
\begin{proof}
Let $\calF$ be an isocrystal on $U$ overconvergent along $Y \setminus U$
whose restriction to $\Isoc^\dagger(U,X/K)$ is isomorphic to
the restriction of $\calE$. Then $\calF$ has constant monodromy along
$X \setminus U$, so
by Theorem~\ref{T:no log extension} it extends to an isocrystal
$\calG$ on $X$ overconvergent along $Y \setminus X$.
If we compare $\calE$ and the restriction of $\calG$ to
$\Isoc(X/K)$, we see that they become isomorphic in
$\Isoc^\dagger(U,X/K)$; by the full faithfulness aspect of
Theorem~\ref{T:no log extension}, they are isomorphic in
$\Isoc(X/K)$. This yields the desired result.
\end{proof}
\begin{remark}
Proposition~\ref{P:contagion} seems tantalizing close to, but distinct
from, a result of Matsuda and Trihan \cite[Theorem~1]{matsuda-trihan}.
The latter says (with more restrictive hypotheses, namely
discreteness of $K$ and presence of a Frobenius structure)
that on a curve, whether a convergent isocrystal is overconvergent
can be checked locally. It would be interesting
to give a higher-dimensional analogue of the result of Matsuda-Trihan;
our methodology is unsuited for this,
as we must have some sort of global overconvergence
in order to make any monodromy constructions.
\end{remark}

\begin{remark} \label{R:contagion}
If one knew that restriction from $\Isoc^\dagger(U,X/K)$ to $\Isoc(U/K)$ were
fully faithful, one could perform the comparison in 
Proposition~\ref{P:contagion} in $\Isoc(U/K)$ instead.
By \cite[Theorem~1.1]{me-full}
this full faithfulness is known under some additional restrictions:
$K$ must be discretely valued, $X$ must be proper (so that
$\Isoc^\dagger(U,X/K) = \Isoc^\dagger(U/K)$), and 
one must consider isocrystals with Frobenius structures.
(Strictly speaking, \cite[Theorem~1.1]{me-full} only extends morphisms
which commute with the Frobenius structures, but it is not difficult
to remove that restriction.)
\end{remark}

\section{Logarithmic extensions}
\label{sec:logext}

We now turn to the problem of extending isocrystals into log-isocrystals.
The context in which we will do this is the work of Shiho \cite{shiho1,
shiho2}, which constructs categories of ``convergent log-isocrystals''
analogous to the convergent isocrystals of Berthelot-Ogus;
indeed, the bulk of this section will be spent reviewing foundational
aspects of logarithmic structures on schemes, then making explicit
one of Shiho's constructions for a smooth pair (a smooth variety
equipped with a strict normal crossings divisor).

In principle, our methods can also be used to construct
``overconvergent log-isocrystals''; the trouble is that there is no
analogue of Shiho's work to use as the foundation. Since building
such a foundation is somewhat orthogonal to our present purposes, we
will not do so here; see Remark~\ref{R:no overconvergent} for further
discussion.

\setcounter{equation}{0}

\begin{convention} \label{conv:Zariski}
We continue to assume that the field $K$ has characteristic
0 and residue field $k$.
However, throughout this section, we also assume that $K$ is
complete with respect to a \emph{discrete} valuation; this is in
order to invoke Shiho's results.
Also, ``locally'' on a scheme or formal scheme (e.g., in the notion of a sheaf)
will always mean locally for the Zariski
topology; note that though some of the constructions can be made using the
\'etale topology (as in \cite{kato} or \cite{shiho1}), the
relevant constructions in \cite{shiho2} require working Zariski locally.
Finally, all monoids to which we refer will be
commutative, and for $M$ a monoid,
$M^{\gp}$ will denote the group generated by $M$.
\end{convention}

\begin{convention}
For definitions and notations regarding log-schemes,
see \cite{kato} and \cite[Section~2.1]{shiho1}.
We follow the following convention from \cite{shiho1}:
if $(X, \calM)$ is a $p$-adic log formal scheme, we get an ordinary
log scheme by reduction modulo $p^n$; we call the result $(X_n, \calM_n)$.
Ditto for morphisms between $p$-adic log formal schemes.
\end{convention}

\begin{remark}
It was explained to us by Shiho that the results of this section
can be extended to the case of nondiscrete $K$. We omit this verification
here, since it requires repeating a fair bit of \cite{shiho1} in restricted
generality, it being not completely clear
whether one can redo \cite{shiho1} at full strength for nondiscrete $K$.
\end{remark}

\subsection{Convergent log-isocrystals}

In the process of introducing Fontaine-Illusie logarithmic structures,
Kato constructed the category of crystals on a log-scheme and checked some
of its basic properties.  The analogue of the Berthelot-Ogus
constructions of convergent isocrystals in the logarithmic setting
is the work of Shiho \cite{shiho1, shiho2}. We will not recall Shiho's
original definition here; rather, we will use the alternate description
in the case of interest provided by \cite[Proposition~2.2.7]{shiho2}.

\begin{hypo} \label{H:shiho}
Let $(X,\calM)$ be a fine log scheme over $k$,
and let $i: (X, \calM) \hookrightarrow (P, \calL)$ be a closed immersion
of $(X, \calM)$ into a noetherian fine log formal scheme $(P, \calL)$
over $\Spf \gotho_K$ whose underlying scheme is of finite type over $k$.
Assume also that there exists a factorization of $i$
of the form
\begin{equation} \label{eq:shiho factor}
(X, \calM) \stackrel{i'}{\to} (P', \calL') \stackrel{f'}{\to} (P, \calL),
\end{equation}
in which $i'$ is an exact closed immersion and $f'$ is a formally log
\'etale morphism.
\end{hypo}

By \cite[Lemma~2.2.2]{shiho2}, one has the following.
\begin{lemma} \label{L:log tube}
Under Hypothesis~\ref{H:shiho},
let $\widehat{P}'$ be the completion of $P'$ along $X$.
Then the rigid analytic space $\widehat{P}'_K$ is independent of
the choice of the factorization, up to canonical isomorphism.
\end{lemma}

\begin{defn} \label{D:log tubular}
Under Hypothesis~\ref{H:shiho}, 
we write $](X,\calM)[_{(P,\calL)}$ for the space
$\widehat{P}'_K$ defined in Lemma~\ref{L:log tube}; for brevity,
we also notate it by $]X[_P^{\log}$ if the sheaves of monoids are to
be understood. Define the specialization map
$\speci: ]X[_P^{\log} \to X$ as the
composite of the ordinary specialization map $\speci: \widehat{P'}_K \to 
\widehat{P}'_k$ with the map $\hat{f}'_k$.
\end{defn}

\begin{remark} \label{R:Zariski}
It is shown in \cite[(2.2.1)]{shiho2}
that Hypothesis~\ref{H:shiho} and
Definition~\ref{D:log tubular} admit a natural sheafification
for the Zariski topology, but it is not clear whether this is true
for the \'etale topology. Shiho handles this by hypothesizing that
$(X, \calM)$ and $(P, \calL)$ are of ``Zariski type'',
i.e., is Zariski locally associated to a finitely generated monoid;
given Convention~\ref{conv:Zariski},  this is automatic for us.
\end{remark}

\begin{hypo} \label{H:shiho2}
Suppose
\[
\xymatrix{
(X, \calM) \ar^i[r] \ar_f[d] & (P, \calL) \ar_g[d] \\
\Spec k \ar^{\iota}[r] & \Spf \gotho_K
}
\]
is a commuting diagram, where the top row satisfies Hypothesis~\ref{H:shiho},
the log structures on the bottom row are trivial,
and $g$ is formally log smooth. 
For $j \in \NN$, let
$(P(j), \calL(j))$ denote the $(j+1)$-st fibre product of
$(P,\calL(j))$ over $\Spf \gotho_K$, and let
$i(j): (X, \calM) \to (P(j), \calL(j))$ be the locally closed immersion
induced by $i$ (and the diagonal $X \mapsto X(j)$). It can be
shown \cite[Proposition~2.2.4]{shiho2} that
each of the $i(j)$ also satisfies Hypothesis~\ref{H:shiho} Zariski locally.
\end{hypo}

\begin{defn} \label{D:convergent log iso}
Under Hypothesis~\ref{H:shiho2}, 
define a \emph{convergent log-isocrystal} on $(X, \calM)$
(with respect to $i$) to be a pair $(\calE, \epsilon)$, where
$\calE$ is a coherent $\calO_{]X[^{\log}_P}$-module and
$\epsilon: \pi_2^* (\calE) \stackrel{\sim}{\to} \pi_1^* (\calE)$ is an
isomorphism of $\calO_{]X[^{\log}_{P(1)}}$-modules
such that $\Delta^*(\epsilon) = \id_{\calE}$, where $\Delta:
(P, \calL) \to (P(1), \calL(1))$ is the diagonal, and
the cocycle condition $\pi_{12}^*(\epsilon) \circ \pi_{23}^* (\epsilon)
= \pi_{13}^*(\epsilon)$ holds on $]X[^{\log}_{P(2)}$.
Then by \cite[Proposition~2.2.7]{shiho2}, the category of
convergent log-isocrystals on $(X, \calM)$
in this sense is equivalent to the
category of convergent log-isocrystals on $(X, \calM)$ in Shiho's sense;
in particular, the former is canonically independent of the
choice of $i$.
\end{defn}

\begin{remark}
The specific analogue of \cite[Proposition~2.2.7]{shiho2} in the 
nonlogarithmic
case is the combination of Ogus's description of convergent
isocrystals in terms of a canonical sequence of enlargements
\cite[Proposition~2.11]{ogus1} and Berthelot's reinterpretation
of Ogus's description in terms of rigid analytic geometry
\cite[(2.3.4)]{ber2}.
\end{remark}

\subsection{Log-$\nabla$-modules and Shiho's construction}

We now clarify how to construct a convergent log-isocrystal,
in the sense of Definition~\ref{D:convergent log iso},
from a log-$\nabla$-module arising as an extension of an
overconvergent isocrystal.

\begin{hypo} \label{H:log extension}
Let $F = (X,P,j)$ be a small frame with $X = P_k$, and suppose that the
differentials of $t_1, \dots, t_n \in \Gamma(P, \calO_P)$ freely
generate $\Omega^1$. Choose $m \leq n$, let $Q$ denote the zero
locus of $t_1 \cdots t_m$ on $P$, and put $Z = Q_k$ and
$U = X \setminus Z$. Since $Z$ and $Q$ are (relative) strict 
normal crossings divisors on $X$ and $P$, respectively,
we obtain log structures $(X, \calM)$ and $(P, \calN)$
and a morphism $i: (X, \calM) \to (P, \calN)$
satisfying Hypothesis~\ref{H:shiho2}.
Define $X(j)$ and $P(j)$ accordingly, and put
\begin{align*}
Z(j) &= \pi_1^{-1}(Z) \cup \cdots \cup \pi_j^{-1}(Z) \subset X(j) \\
Q(j) &= \pi_1^{-1}(Q) \cup \cdots \cup \pi_j^{-1}(Q) \subset P(j).
\end{align*}
\end{hypo}

In order to apply Definition~\ref{D:convergent log iso}, we 
need to explicitly identify the spaces $]X[^{\log}_{P(j)}$ for $j=1,2$.
\begin{defn}
Under Hypothesis~\ref{H:log extension},
for $i=1, \dots, m$ and $l=1, \dots, j$, put $t^{(l)}_i = \pi_l^* (t_i)$.
Let $\widehat{\AAA^{mj^2}_{\gotho_K}}$ be the completion of the
affine space with coordinates
$u^{(l,l')}_i$ for $i=1,\dots, m$ and $l,l' = 1, \dots, j$.
Let $P'(j)$ be the closure in $P(j) \times \widehat{\AAA^{mj^2}_{m,\gotho_K}}$
of the graph of the map
$P(j)^{\triv} \to P(j) 
\times \widehat{\AAA^{m j^2}_{\gotho_K}}$, induced by the functions
$t^{(l)}_i/t^{(l')}_i$ for $i=1, \dots, m$ and $l,l' = 1, \dots, j$.
Let $i'(j): X \to P'(j)$ be the map induced by composing
$i(j): X \to P(j)$ with the rational map $P(j) \dashrightarrow P'(j)$; note that
$i'(j)$ is a \emph{regular} map, not just a rational map.
Let $f'(j): P'(j) \to P(j)$ 
be the map
obtained by composing the injection $P'(j) \hookrightarrow
P(j) \times \widehat{\AAA^{mj^2}_{\gotho_K}}$ with the first projection
from $P(j) \times \widehat{\AAA^{mj^2}_{\gotho_K}}$; then
$i(j) = f'(j) \circ i'(j)$.
\end{defn}

\begin{lemma} \label{L:blowing up}
Let $X = \Spec A \to S = \Spec B$ be a morphism of integral affine schemes,
and suppose that for
some $n \geq 2$, the differentials of $t_1, \dots, t_n
\in A$ freely generate $\Omega^1_{X/S}$. Put
$A' = A [t_1/t_2, t_2/t_1]$ and $X' = \Spec A'$.
Then $\Omega^1_{X'/S}$ is freely generated by
the differentials of $t_1/t_2, t_2, \dots, t_n$.
\end{lemma}
\begin{proof}
Given $f \in A[t_1/t_2]$, we can write $f = (t_1/t_2)^l a$ with $a \in A$
for some $l \in \NN$.
Then $df = l a (t_1/t_2)^{l-1} d(t_1/t_2) + (t_1/t_2)^l da$ is a linear 
combination
of $d(t_1/t_2), dt_2, \dots, dt_n$,
since $dt_1 = t_2 d(t_1/t_2) + (t_1/t_2) dt_2$ can be reexpressed
in terms of $d(t_1/t_2)$ and $dt_2$.
The same is true if $f \in A[t_2/t_1]$.
Finally, any element of $A[t_1/t_2,t_2/t_1]$ can be written as the sum
of an element of $A[t_1/t_2]$ and an element of $A[t_2/t_1]$, so
$\Omega^1_{X'/S}$ is indeed generated by $d(t_1/t_2), dt_2, \dots, dt_n$.

On the other hand,
suppose that $f d(t_1/t_2) + e_2 dt_2 + \cdots + e_n dt_n = 0$
in $\Omega^1_{X'/S}$ for some $e_2,
\dots, e_n,f \in A[t_1/t_2,t_2/t_1]$. 
By multiplying through by a power of
$t_1 t_2$, 
we may reduce to the case where $e_2, \dots, e_n \in A$
and $f \in t_2^2 A$. Then
\[
0 = (f/t_2) dt_1 + (e_2 - t_1 f/t_2^2) dt_2 + e_3 dt_3 + \cdots + e_n dt_n
\]
(using the fact that $X$ is integral, so the division $f/t_2^2$ makes sense),
so we must have $e_3 = \cdots = e_n = f/t_2 = 0$, so that $f=0$, and then
$e_2 - t_1 f/t_2^2 = 0$, so that $e_2 = 0$. Thus
$d(t_1/t_2), dt_2, \dots, dt_n$ freely generate
$\Omega^1_{X'/S}$, as desired.
\end{proof}
\begin{remark}
All that Lemma~\ref{L:blowing up} is doing is blowing up the smooth
$S$-scheme $X$ along the smooth $S$-subscheme $t_1 = t_2 = 0$.
\end{remark}
\begin{cor} \label{C:relative strict}
The sheaf $\Omega^1_{P'(j)/\gotho_K}$ is freely generated by
the differentials of the regular functions
\[
t^{(1)}_i, t^{(2)}_i/t^{(1)}_i, \dots, t^{(j)}_i/t^{(j-1)}_i \quad
(i=1, \dots, m), \qquad
t^{(1)}_i, t^{(2)}_i, \dots, t^{(j)}_i \quad (i=m+1, \dots, n).
\]
In particular, 
the divisor $f'(j)^{-1}(Q(j))$ is a relative strict normal crossings
divisor on $P'(j)$ (relative to $\gotho_K$).
\end{cor}

\begin{remark} \label{R:inverse image}
In fact, $f'(j)^{-1}(Q(j))$ is quite simple: for $i=1,\dots,m$,
the zero locus of $t^{(1)}_i$ on $P'(j)$ is isomorphic to the zero locus of
$t_i$ on $P$ via the first projection from $P(j)$, and the union of these
loci is all of $f'(j)^{-1}(Q(j))$ since the functions
$t^{(2)}_i/t^{(1)}_i, \dots, t^{(j)}_i/t^{(j-1)}_i$ are all invertible
on $P'(j)$.
\end{remark}

\begin{defn}
Let $\calL'(j)$ be the canonical log-structure on $P'(j)$ associated
to $f'(j)^{-1}(Q(j))$, which is a relative strict normal
crossings divisor by Corollary~\ref{C:relative strict}.
Then $f'(j)$ gives rise to a natural morphism
$(P'(j), \calL'(j)) \to (P(j), \calL(j))$.
On the other hand, since $i'(j)^{-1}(f'(j)^{-1}(Q(j))) = Z$,
$i'(j)$ extends to a morphism $(X, \calM) \to
(P'(j), \calL'(j))$ of log formal schemes,
and the composition $f'(j) \circ i'(j)$ coincides with $i(j)$ as 
a map of log formal schemes.
\end{defn}

\begin{remark} \label{R:strict normal2}
Suppose that $(X,Z)$ and $(X',Z')$ are (formal) smooth pairs,
and $i:X \to X'$ is a closed immersion such that $i^{-1}(Z') \subseteq Z$
as (formal) schemes.
Then $i$ induces a morphism
between the canonical log schemes $\calM$ and 
$\calM'$ corresponding to $(X,Z)$ and $(X',Z')$.
On an open subset $U$ of $X$, $(i^* \calM')/\calO_X^*
\cong i^{-1} (\calM'/\calO_X^*)$ is generated by the components
of $Z'$ meeting $i(U)$, whereas $\calM/\calO_X^*$ is generated by the
components of $Z$ meeting $U$. 
Hence a sufficient (but not necessary) condition for the map
$i^* (\calM') \to \calM$ to be an isomorphism is that each
component of $Z'$ that meets $X$ does so in a single component of $Z$,
and no two components of $Z'$ meet $X$ along the same component of $Z$.
\end{remark}

\begin{lemma}
The factorization
\[
(X, \calM) \stackrel{i'(j)}{\to} (P'(j), \calL'(j)) \stackrel{f'(j)}{\to} 
(P(j), \calL(j))
\]
of $i(j)$ satisfies Hypothesis~\ref{H:shiho}.
\end{lemma}
\begin{proof}
We first check that the map $i'(j)$ is an exact closed immersion
using the criterion from Remark~\ref{R:strict normal2}.
Namely, each component of $f'(j)^{-1}(Q(j))$ is a component
of the zero locus of $t^{(1)}_i$ for some $i \in \{1, \dots, m\}$,
which meets $X$ in the corresponding component of the zero locus of
$t_i$. In particular, each component $f'(j)^{-1}(Q(j))$ meeting
$X$ does so in a single component and no two of these intersections
coincide. Hence the map $i'(j)^* \calL'(j) \to \calM$ is an isomorphism,
and $i'(j)$ is an exact closed immersion.

We next check that the map $f'(j)$ is formally log \'etale.
The structural map 
$(P'(j), \calL'(j)) \to \Spf \gotho_K$ is formally log smooth; by
the formal analogue of \cite[Proposition~3.12]{kato},
it then suffices to show that
the map
$f'(j)^* (\Omega^1_{P(j)/K}) \to \Omega^1_{P'(j)/K}$ is an isomorphism.
But this is a straightforward consequence of the fact that
\[
d\log(t^{(l)}_i/t^{(l')}_i)
= d\log(t^{(l)}_i) - d\log(t^{(l')}_i):
\]
as we adjoin each fraction $t^{(l)}_i/t^{(l')}_i$,
we do not change $\Omega^1$.
\end{proof}

We now have the tools with which to construct convergent log-isocrystals
on the log schemes associated to strict normal crossings divisors on 
smooth $k$-varieties. Before doing so, we must collect a bit of
information about log-$\nabla$-modules.

\subsection{Log-$\nabla$-modules and unipotent monodromy}

\begin{defn}
Under Hypothesis~\ref{H:log extension}, let $\calE$ be a
log-$\nabla$-module on $]X[$ with respect to $t_1, \dots, t_n$.
We say $\calE$ is \emph{convergent} if the restriction of
$\calE$ to a strict neighborhood of $]U[$ in $]X[$ is overconvergent
along $Z$.
\end{defn}

We now have the following limited logarithmic analogue of 
Theorem~\ref{T:no log extension}.
(Note however that the work has been done already in the proof
of Lemma~\ref{L:direct extension}.)
\begin{prop} \label{P:log direct extension}
Under Hypothesis~\ref{H:log extension},
let $\calE$ be a $\nabla$-module on a strict neighborhood
on $]U[$ in $]X[$ which is 
overconvergent along $Z$. Then $\calE$ has unipotent monodromy
along $Z$ if and only if $\calE$ extends to a 
convergent log-$\nabla$-module
on $]X[$ with nilpotent residues. Moreover, the restriction functor, from
convergent log-$\nabla$-modules with nilpotent residues on $]X[$
to isocrystals on $U$ overconvergent along $Z$, is fully faithful.
\end{prop}
\begin{proof}
By covering $X$ with affines,
we may reduce to the case where 
we may repeatedly apply Lemma~\ref{L:direct extension} to obtain the desired
result.
\end{proof}
\begin{remark} \label{R:need nilpotent}
The full faithfulness assertion in Proposition~\ref{P:log direct extension}
depends crucially on the nilpotent residues hypothesis. 
This is analogous to the situation in
\cite[II.5]{deligne}, where logarithmic extensions with nilpotent residue are
``canonical'' and logarithmic extensions with arbitrary residue are not;
indeed, one of the simplest examples in that setting is relevant here also.
Namely, put $P = \Spf K \langle t \rangle$, $X = P_k = \AAA^1_k$,
and $U = \AAA^1_k \setminus \{0\}$,
let $n$ be a positive integer,
and let $\calE$ be the $\nabla$-module on $P_K$
generated by a single element $\bv$ such that $\nabla \bv = n \bv \otimes
\frac{dt}{t}$ for some $n \in \NN$.
Then one easily verifies that
$\nabla$ is overconvergent along $X \setminus U$,
and the kernel of $\nabla$ on
$]X[$ is trivial, but the kernel of $\nabla$ on any strict 
neighborhood of $]U[$ in $]X[$ not containing the point $t=0$
includes the section $t^{-n} \bv$.
(A similar point arises in \cite{lestum-trihan}, which is concerned
with the passage from a log-$F$-crystal to an isocrystal on
the log-trivial subscheme overconvergent along the complement.)
\end{remark}

\begin{lemma} \label{L:log convergent}
Under Hypothesis~\ref{H:log extension},
let $\calE$ be a convergent log-$\nabla$-module on $P_K$.
Then for any $\bv \in \Gamma(]X[, \calE)$
and any $\eta \in (0,1)$, the multisequence
\[
\frac{1}{i_1! \cdots i_n!} \left( \prod_{j=1}^n \prod_{l=0}^{i_j-1}
\left( t_j \frac{\del}{\del t_j} - l \right) \right) \bv
\]
is $\eta$-null.
\end{lemma}
\begin{proof}
Since $\calE$ restricts to a convergent isocrystal on $U$,
the multisequence
\[
\frac{1}{i_1! \cdots i_n!} \left( \prod_{j=1}^n
\frac{\del^{i_j}}{\del t_j^{i_j}} \right) \bv
\]
is $\eta$-null on $]U[$ by the definition of $\eta$-convergence
plus Proposition~\ref{P:taylor}.
Since $|t_i| \leq 1$ for each $i$, the multisequence
\[
\frac{t_1^{i_1} \cdots t_n^{i_n}}{i_1! \cdots i_n!} \left( \prod_{j=1}^n
\frac{\del^{i_j}}{\del t_j^{i_j}} \right) \bv
\]
is also $\eta$-null on $]U[$. However, this is precisely the desired
multisequence, and the fact that it is $\eta$-null on $]U[$ implies
the fact that it is $\eta$-null on $]X[$. Namely, this follows
from the fact that the spectral
seminorm on $\calO(]U[)$ restricts to the spectral
seminorm on $\calO(]X[)$, which is true 
because $U$ is open dense in $X$.
\end{proof}

\subsection{Convergent log-isocrystals and log-$\nabla$-modules}

With the constructions of the previous subsection in hand, we can now
explicitly describe convergent log-isocrystals, in the case of
the log structure associated to a smooth pair, in terms
of log-$\nabla$-modules.

\begin{theorem} \label{T:convergent log}
Under Hypothesis~\ref{H:log extension}, there is an equivalence between
the category of
convergent log-isocrystals on $(X,Z)$ 
and the category of convergent log-$\nabla$-modules on $P_K$.
\end{theorem}
\begin{proof}
Suppose $\calE$ is a convergent log-isocrystal on $(X,Z)$
in the sense of Definition~\ref{D:convergent log iso}.
Then $\calE$ restricts to an isocrystal on $X$ overconvergent
along $Z$, and hence to an overconvergent $\nabla$-module on
some strict neighborhood $V$ of $]X[_P$. Moreover, by
\cite[Proposition~1.2.7]{shiho2}, the
isomorphism $\epsilon: \pi_2^* (\calE) \to \pi_1^* (\calE)$ on
the second infinitesimal neighborhood of $X$ in $P'(1)$
defines a log-connection $\nabla: \calE \to \calE \otimes 
\Omega^{1,\log}_{P_K/K}$ extending the connection on $V$.
This yields the data of a convergent log-$\nabla$-module on $P_K$.

Conversely, suppose that $\calF$ is a convergent
log-$\nabla$-module on $P_K$.
Write $u_i$ for the function $t^{(2)}_i/t^{(1)}_i$ on $P'(1)$.
Following \cite[(6.7.1)]{kato}, we observe that
the isomorphism $\epsilon: \pi_2^*(\calF) \stackrel{\sim}{\to}
 \pi_1^* (\calF)$ over a
suitable strict neighborhood of $]U[_{P'(1)}$ in $]X[_{P'(1)}$
induced by $\nabla$ can be written in the form
\[
1 \otimes \bv \mapsto \sum_{i_1,\dots,i_n=0}^\infty
\left( \prod_{j=1}^n \frac{(u_j - 1)^{i_j}}{i_j!}
\right) \otimes \left( \prod_{j=1}^n \prod_{l=0}^{i_j-1}
\left( t_i \frac{\del}{\del t_i} - l \right) (\bv) \right).
\]
By Lemma~\ref{L:log convergent}, this series converges uniformly
on any affinoid subspace of $]X[_{P'(1)}$ of the form
$\max_j\{|u_j-1|\}\leq \lambda$ for $\lambda \in (0,1) \cap \Gamma^*$. Hence
$\epsilon$ is defined on all of $]X[_{P'(1)} = ]X[^{\log}_{P(1)}$.

We now have an isomorphism $\epsilon: \pi_2^* (\calF) \stackrel{\sim}{\to}
 \pi_1^* (\calF)$
on $]X[^{\log}_{P(1)}$ satisfying $\Delta^*(\epsilon) = \id$.
It is straightforward to check that
the cocycle condition $\pi_{12}^*(\epsilon) \circ \pi_{23}^* (\epsilon)
= \pi_{13}^*(\epsilon)$ holds on $]X[^{\log}_{P(2)}$ from the formula,
but easier is to deduce it by restricting to a strict 
neighborhood of $]U[_{P'(2)}$, where it holds because of the equivalence
of categories between ordinary overconvergent
isocrystals and overconvergent $\nabla$-modules.

We conclude that every convergent log-$\nabla$-module
on $P_K$ does indeed give rise to a convergent
log-isocrystal. This establishes the desired equivalence.
\end{proof}
\begin{remark}
Note that the equivalence in Theorem~\ref{T:convergent log}
is compatible with restriction to an open subscheme, so in
principle its statement can be ``sheafified''. 
\end{remark}
\begin{remark} \label{R:no overconvergent}
While Lemma~\ref{L:direct extension} can also be applied with
$Y \neq X$ to 
construct ``overconvergent log-$\nabla$-modules'', their interpretation
in the Grothendieckian sense (i.e., as isomorphisms between two pullbacks
to the diagonal) seems subtle. Probably the right thing to do is
to globally replace tubes with strict neighborhoods throughout
the proof of Theorem~\ref{T:convergent log}; however, in the
absence of a ``reference category'' of overconvergent log-isocrystals,
one then has to check all the relevant compatibilities by hand.
The main problem is that we do not presently have an 
``overconvergent topos'' analogizing \cite{ogus2}; however,
the ongoing
work of le~Stum mentioned earlier \cite{lestum, lestum2} seems to be
heading in the right direction, and it is possible it will
ultimately be adapted to include logarithmic structures.
In the meantime, however, we will stick to convergent log-isocrystals.
\end{remark}

\begin{defn}
Under Hypothesis~\ref{H:log extension}, we say a
convergent log-isocrystal on $(X,Z)$ \emph{has nilpotent
residues} if its image under the functor of Theorem~\ref{T:convergent log}
is a log-$\nabla$-module with nilpotent residues. More generally,
if $X$ is a smooth $k$-variety and $Z$ is a strict normal crossings
divisor on $X$, we say a convergent log-isocrystal $\calE$ on $(X,Z)$
has nilpotent residues if there is an open cover $U_1, \dots, U_n$
of $X$ such that each pair $(U_i \cap X, U_i \cap Z)$ satisfies
Hypothesis~\ref{H:log extension}, and the restriction of
$\calE$ to $U_i \cap X$ has nilpotent residues. The same is then
true on any open cover.
\end{defn}

{}From Theorem~\ref{T:convergent log}, we obtain the following.
\begin{theorem} \label{T:extension}
Let $U \hookrightarrow X$ be an open immersion of smooth
$k$-varieties such that $Z = X \setminus U$ is a strict
normal crossings divisor on $X$. Let $\calE$ be an isocrystal
on $U$ overconvergent along $Z$. Then $\calE$ has unipotent monodromy
along $Z$ if and only if $\calE$ extends to
a convergent log-isocrystal with nilpotent residues
on $(X,Z)$. Moreover, the restriction functor, from
convergent log-isocrystals with nilpotent residues on
$(X,Z)$ to isocrystals on $U$ overconvergent along $Z$,
is fully faithful.
\end{theorem}
\begin{proof}
Everything being asserted is Zariski local, so
we may reduce to the case where
Hypothesis~\ref{H:log extension} holds. In this case,
Proposition~\ref{P:log direct extension} and
Theorem~\ref{T:convergent log} together yield the claim.
\end{proof}

\begin{remark}
The word ``strict'' is probably not necessary in
Theorem~\ref{T:extension}; removing it would require performing
an appropriate \'etale descent (but beware of some technical
problems, as in Remark~\ref{R:Zariski}). However, in the desired application
to semistable reduction, we can always get to the strict normal
crossings situation using an alteration, in the manner of
de Jong \cite{dejong1}.
\end{remark}

\begin{remark}
It should be possible to improve the full faithfulness conclusion
of Theorem~\ref{T:extension} to allow restriction
all the way to the category of convergent isocrystals
on $U$. In fact, this is possible under additional hypotheses;
see Remark~\ref{R:contagion}.
\end{remark}

\begin{remark}
In some cases, one may want to apply Theorem~\ref{T:extension} to
construct logarithmic extensions of crystals in coherent
$\calO$-modules, rather than isocrystals. This should be
a straightforward consequence of the fact that
isocrystals can be viewed as elements of the isogeny category of
crystals (as in \cite{ogus1}), but we have not checked any details.
\end{remark}

\subsection{Extension classes of log-isocrystals}

In the logarithmic setting, one can show that restriction to the
log-trivial subscheme preserves extension classes.

\begin{prop}
Let $(X,Z)$ be a smooth pair, and let $\calE_1$ and $\calE_2$ be
convergent log-isocrystals with nilpotent residues on $(X,Z)$. Then
$\Ext^1(\calE_1, \calE_2)$ is the same whether computed in the 
category of convergent log-isocrystals on $(X,Z)$ or in the 
category of isocrystals on $U = X \setminus Z$ overconvergent along $Z$.
\end{prop}
\begin{proof}
Recall that the Yoneda Ext group
$\Ext^1(\calE_1,\calE_2)$ classifies short exact sequences
\[
0 \to \calE_1 \to \calF \to \calE_2 \to 0.
\]
Let $\Ext_X$ and $\Ext_U$ denote the group $\Ext^1(\calE_1, \calE_2)$
computed in the category of convergent log-isocrystals on $(X,Z)$
and in the category of isocrystals on $U$ overconvergent along $Z$,
respectively; then restriction gives a map $\Ext_X \to \Ext_U$.
Note that this map is injective thanks to full faithfulness of
restriction (Theorem~\ref{T:extension}): any isomorphism over $U$
between two short exact sequences over $X$ extends to $X$.

To see that $\Ext_X \to \Ext_U$ is surjective, note that if
$\calF$ fits into a sequence over $U$, then $\calF$ has unipotent
monodromy along $Z$, because $\calE_1$ and $\calE_2$ both do.
Hence $\calF$ extends to a convergent log-isocrystal on $(X,Z)$,
as do the maps $\calE_1\to \calF$ and $\calF \to \calE_2$
by Theorem~\ref{T:extension}. Hence $\Ext_X \to \Ext_U$ is 
surjective, and thus is a bijection as desired.
\end{proof}

\section{Conclusion: a look ahead}
\label{sec:finale}

We conclude by cataloging some of the questions we will be discussing later
in this series of papers, in the terminology we have established.
Note that this section is intended as a ``pre-introduction'' to
the subsequent papers, and so statements here have not been made in a
precise fashion; they will be articulated properly in due course.

\subsection{Semistable reduction: Shiho's conjecture}
\label{subsec:shiho}

We give the statement of Shiho's conjecture \cite[Conjecture~3.1.8]{shiho2},
or in our terminology, the ``semistable reduction problem''.
First, we must recall the notion of a Frobenius structure on an isocrystal.
\begin{defn} \label{D:Frobenius}
Suppose that $\sigma_K: \gotho_K \to \gotho_K$ is an endomorphism
lifting the $p^a$-power Frobenius map on $k$, for some positive integer
$a$. Let $X \hookrightarrow Y$ be an open immersion of $k$-varieties.
A \emph{Frobenius structure (of order $a$)} on an isocrystal $\calE$ on $X$
overconvergent along $Y \setminus X$ is an isomorphism 
$F^*_X \sigma^*_K \calE \stackrel{\sim}{\to} \calE$, where
$F_X$ is the relative $p^a$-power Frobenius.
An isocrystal equipped with a Frobenius structure of order $a$
is called an \emph{$F^a$-isocrystal}.
\end{defn}

\begin{conj}[Shiho] \label{conj:Shiho}
Assume the field $k$ is perfect.
Let $X$ be a smooth $k$-variety and let $\calE$ be an overconvergent
$F^a$-isocrystal on $X$. Then there exists a proper, surjective,
generically \'etale morphism $f: X_1 \to X$, and an open immersion
$j: X_1 \hookrightarrow \overline{X_1}$ of $X_1$ into a smooth
projective $k$-variety in which the complement $D = \overline{X_1}
\setminus X_1$ is a strict normal crossings divisor, 
such that $f^* \calE$ extends to a
convergent $F^a$-log-isocrystal $\calF$ on $(\overline{X_1},D)$.
\end{conj}
\begin{remark}
Absent the isocrystal, the existence of the maps $f$ and $j$ is 
the content of de Jong's alterations theorem \cite[Theorem~4.1]{dejong1};
indeed, the map $f$ is precisely an alteration in de Jong's sense.
\end{remark}

\begin{remark}
Note that it is actually enough to show that $f^* \calE$ extends as
a convergent log-isocrystal; then the Frobenius structure will extend
from $X_1$ to $\overline{X_1}$ thanks to the full faithfulness
aspect of Theorem~\ref{T:extension}. Note also that a convergent
log-$F$-isocrystal necessarily has nilpotent residues.
\end{remark}

\begin{remark}
Shiho's conjecture is a higher-dimensional version of
de Jong's formulation of Crew's conjecture \cite{dejong3};
the case where $X$ is a curve is known to follow from
the $p$-adic local monodromy theorem \cite{me-semicurve}.
As noted in the introduction, its resolution is expected to have
various consequences for the theory of rigid cohomology,
especially in the relative setting, and perhaps for the theory
of arithmetic $\mathcal{D}$-modules, which are to the isocrystals
considered here as constructible sheaves are to lisse sheaves in
\'etale cohomology.
\end{remark}

\subsection{Monodromy of exceptional components}

The $p$-adic local monodromy theorem of Andr\'e \cite{andre},
Mebkhout \cite{mebkhout}, and the present author \cite{me-local}
implies a strong statement in the direction of Conjecture~\ref{conj:Shiho}.
(We will describe the exact statement of the $p$-adic local monodromy
theorem and the nature of its application here more thoroughly
later in the series.)
Namely, if one starts with a compactification $X \hookrightarrow \overline{X}$
such that $(\overline{X}, \overline{X} \setminus X)$ is a smooth pair
(which one may do without loss of generality by pulling back 
along an alteration, thanks to de Jong's theorem),
one can construct the maps $f$ and $j$ so that $f$ extends to
a map $\overline{X_1} \to \overline{X}$, and $\calE$ has
unipotent monodromy along each component of $\overline{X_1}
\setminus X_1$
which dominates a component of $\overline{X} \setminus{X}$.

Unfortunately, this statement together with Theorem~\ref{T:extension}
do not suffice to imply Conjecture~\ref{conj:Shiho}, because there
may be components of $\overline{X_1} \setminus X_1$ which do not
dominate any component of $\overline{X} \setminus X$. In order
to deduce Conjecture~\ref{conj:Shiho} along these lines, one must
somehow gain control of the monodromy of these ``exceptional''
divisors. Otherwise, one is forced to alter again, possibly introduce
more exceptional divisors, and perhaps repeat \emph{ad infinitum}
without reaching the desired conclusion.

The control of exceptional divisors will be accomplished by considering
monodromy also along certain ``fake annuli'', corresponding to irrational
valuations on the function field $K(X)$. These form a compact space
(an example of a Gelfand spectrum, as in Berkovich's foundations 
of rigid analytic geometry \cite{berkovich-icm}), so one can 
prove a global quasi-unipotence theorem ``topologically'', by verifying
it on an open neighborhood of each valuation. This is most easily done
for surfaces, so we will focus on that case initially.

It must be stressed that the presence of the exceptional divisors
is not an artifact of the use of de Jong's theorem in lieu
of the as-yet-unknown resolution of singularities in positive
characteristic. That is because the underlying finite cover
given by the $p$-adic local monodromy theorem is typically unavoidably
singular,
due to wild ramification; contrast this situation to what happens in the
complex analytic setting, where one can locally avoid introducing any
singularities by making the right toroidal cover.

\end{document}